\documentclass[a4paper]{article}

\usepackage[a4paper]{geometry}
\usepackage[T1]{fontenc}
\usepackage[utf8]{inputenc}
\usepackage[english]{babel}
\usepackage{amsmath,amsthm,amsfonts,amssymb}
\usepackage{mathtools,graphicx,subfig,float,xcolor}
\usepackage{indentfirst,afterpage,mathrsfs,url}
\usepackage{soul}

\begin{document}


\theoremstyle{definition}
\newtheorem{remark}{Remark}[section]
\theoremstyle{plain}
\DeclarePairedDelimiter{\abs}{\lvert}{ ert}

\title{Numerical simulations of the Gatenby-Gawlinski model\\
with heterogeneous acid diffusion in one space dimension}

\author{
Chiara Simeoni -- \url{chiara.simeoni@univ-cotedazur.fr}\\
Laboratoire de Math\'ematiques J.A. Dieudonn\'e CNRS UMR 7351\\
Universit\'e C\^ote D'Azur, Parc Valrose - 06108 Nice Cedex 2 (France)\\ \\
Elisa Scanu -- \url{e.scanu@qmul.ac.uk}\\
Barts Cancer Institute, Queen Mary University of London\\
Charterhouse Square, EC1M 6AU (United Kingdom)\\ \\
Donato Pera -- \url{donato.pera@univaq.it}\\
Dipartimento di Ingegneria e Scienze dell'Informazione e Matematica\\
Universit\`a degli studi dell'Aquila, via Vetoio 1 - 67100 Coppito, L'Aquila (Italy)\\ \\
Corrado Mascia\footnote{Corresponding author}\quad -- \url{corrado.mascia@uniroma1.it}\\
Dipartimento di Matematica G. Castelnuovo, Sapienza, Universit\`a di Roma\\
piazzale Aldo Moro 2 - 00185 Roma (Italy)\\ \\
}

\date{}
\maketitle

\rightline{\footnotesize{\it Chiara Simeoni passed away untimely from cancer in 2025.}}
\rightline{\footnotesize{\it We will miss her deeply for her wit and strength of spirit.}}
\rightline{\footnotesize{\it And, of course, her extraordinary eloquence...}}
\vskip.25cm

\abstract{In this work, we introduce a variant of the \emph{Gatenby--Gawlinski model} for 
acid-mediated tumor invasion in the one-dimensional experimental setting, accounting for 
heterogeneous diffusion of the lactic acid across the surrounding healthy tissues. 
Numerical simulations are performed by employing finite volume schemes on staggered 
cartesian grids, together with explicit time discretization. 
The effectiveness of such approach is proven by reproducing biologically relevant results like 
the formation of propagating invasion fronts and the emergence of an interstitial gap between
normal and cancerous cells, which is driven by the pH lowering strategy and depends 
significantly on the non-homogenous diffusion rates. 
By means of a comparison analysis, we infer that a homogenization phenomenon arises in 
the long run for appropriate values of the physical parameters, thus paving the way for complex 
applications to interface diffusion problems of invasive processes.




\section{Introduction}

Nowadays, \emph{cancer research} is one of the most active and interdisciplinary investigation fields 
 {because of its evident relevance to our lives, \cite{McKiEtAl06, SBEB21, Wang10, Wein06, 
weinberg1, weinberg2}. 
In particular, in the past few years, the role of mathematical modeling}
has been crucial in supporting theoretical and experimental studies 
 {\cite{AltrLiuMich15, BeckKareAdle20, BellDeAnPrez03, CrisLoweNie03, KuznClaiVolp21, LoweEtAl10}.
Since such models have often several limitations in particular, with respect to the possibility
of proving rigorous results the implementation of numerical simulations and computations 
(the so-called \emph{in silico} experiments) are indispensable for designing new therapies and
prescribing preventive measures, \cite{CrisLoweNie03, Fasano 2, Prez03, SBMA03}.}

Here, we focus on the so-called \emph{Warburg effect} \cite{Warburg1,Warburg2}, and its mathematical 
modeling by means of the \emph{acid-mediated invasion hypothesis}, namely the typical strategy of increasing 
the acidity against the local environment operated by tumors to enhance their growth 
 {at the expense of healthy tissue, \cite{SmalGateMain08, astanin}.}
This phenomenon concerns the intrinsic metabolism 
 {and regulation mechanisms operated by} 
cancer cells, essentially providing their glucose uptake rates. 
Indeed, as  {first} observed by Otto Warburg in the 1920s (and afterwards confirmed through extensive 
 {and comprehensive} 
experiments), tumor cells tend to rely on glycolytic metabolism even in presence of large oxygen amounts 
 {enhancing their capability to degradate the surrounding (healthy) tissue}.
From a strictly biochemical point of view, normal cells undergo glucose metabolism by employing oxidative
phosphorylation pathways, which is the most effective process in terms of Adenosine Triphosphate Production 
(ATP) requiring oxygen as principal resource. 
Tumor cells  {behaviour} seems to forbear this conventional pathway and appeals instead 
 {to a different mechanism, called \emph{glycolysis}, which, in short, induces}
lactic acid fermentation, a product
 {typically}
released in hypoxia regimes (see Fig.\ref{figura1}).
\begin{figure}[bht]
	\centering
	\includegraphics[width=.50\textwidth]{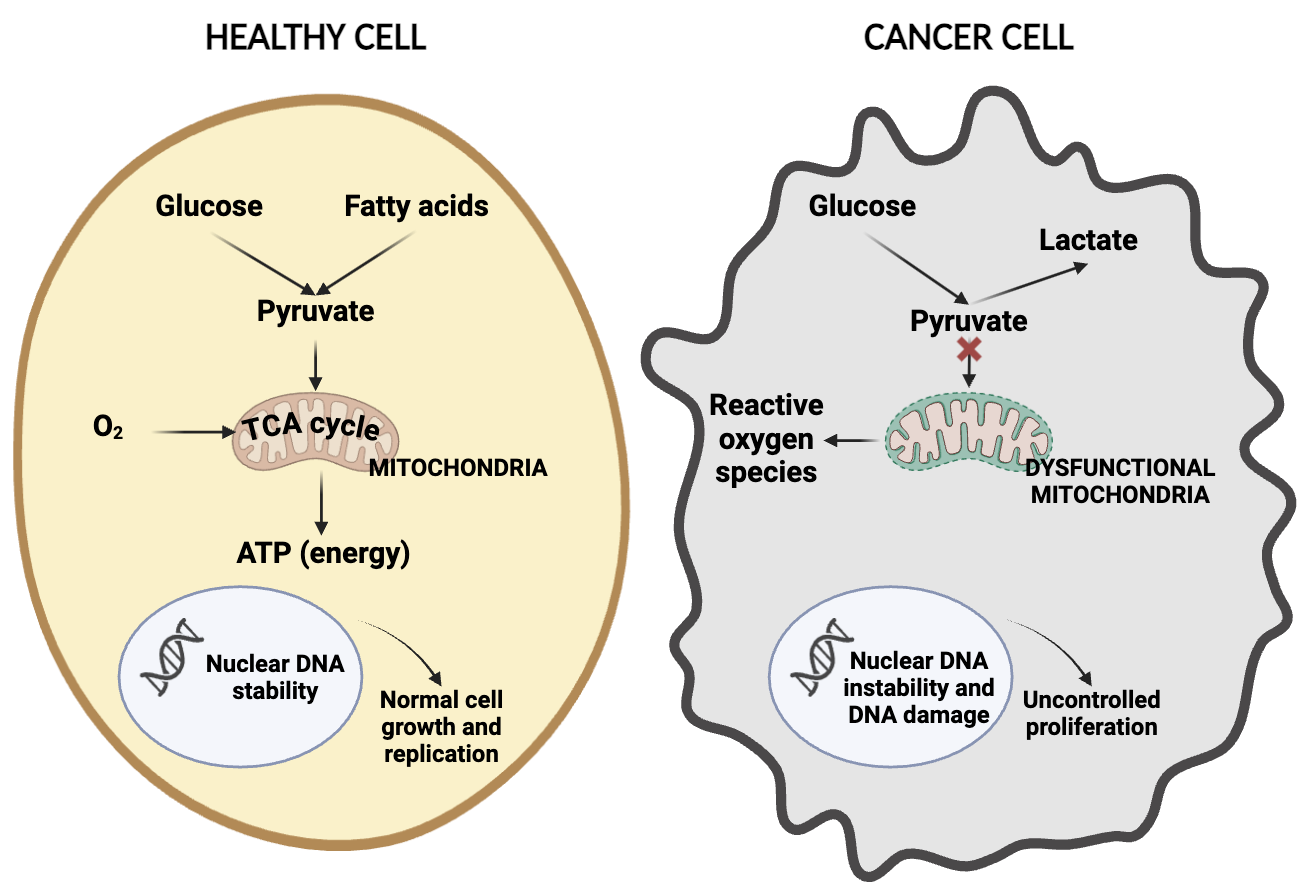}
	\caption{\footnotesize{Pathways of healthy (left) and cancer (right) cells metabolism.}}
	\label{figura1}
\end{figure}
Although the Warburg effect has been intensely studied during the past twenty years
 {(among others, see \cite{Albe23, SmalGateMain08} and reference therein)},
how precisely this phenomenon happens and affects cancer proliferation is  {still} an open problem.
For example, the crucial passage from normal to glycolytic metabolism operated by cancer cells is
still object of  {research} \cite{porporato, sattler}. 
However, by exploiting a combination of modeling and in vitro experiments, it has been 
 {established}
that cancer metabolic changes  {can} define a 
 {local micro-environment, acting on the production of some appropriate \emph{biochemical intermediaries}},
where the better adapted malignant cells overwhelm the others and,  {as a consequence}, 
non-trivial spatial structures are created \cite{archetti, porporato}.
 {These intermediaries have yet to be precisely identified. 
Here, inspired by the Warburg effect, our attention is focused on the molecules of lactic acid pH
and their local fluctuations.}

In the present work, we suppose that the  {beginning of} carcinogenesis has already happened. 
The early growth stages of a primary tumour are a consequence of successive divisions of its initial cells.
When a critical size (of about $10^4$ cells) is reached, further growth requires aggressive actions 
 {exerted} on the surrounding healthy tissue.
 {This can be accomplished in many different ways,}
one of which typically consists  {of} coopting blood vessels to provide
oxygen and nutrients to the expanding colony \cite{bertolini}. 
Another fundamental strategy consists  {of} releasing products, such as lactic acid, 
which favour the degradation of healthy cells. 
 {Therefore, inspired by the Warburg effect, Gatenby and Gawlinsky proposed a simplified model
with three unknowns  representing the healthy cells, the tumour cells and the level of lactic acid concentration
in the tissue, respectively  with the form of a partial differential equations of reaction-diffusion type with a 
coupling which will be detailed in a while, \cite{GateGawl96}.
By means of numerical simulations, they have shown the computational evidence of existence of invasive
fronts that could be the mathematical counterpart corresponding to {\it in vitro} experiments.
In particular, in the presence of a very aggressive tumor, the model predicts the presence of an \emph{interstitial gap}, 
i.e. an intermediate region between the healthy area and the tumor area, in which the main contribution is 
composed only of lactic acid, consistently with what is seen in the simulations carried out {\it in vitro}.
}

Therefore, we make the assumption that acidification, namely the increase caused by lactic acid production, 
is advantageous for the spread of cancer cells, whilst it destroys healthy cells. 
Hence, the original modelling based on reaction-diffusion equations developed by Gatenby and Gawlinski
\cite{GateGawl96, GateGawl03} is suitable  {for performing} numerical investigation and 
 {analysing} travelling wave solutions, because it takes into account the previously explained 
behaviour of both the healthy and cancer cells in an acidic 
 {micro-environment.}

The novelty introduced in the present work consists in taking into account a non-homogeneous
function for the diffusion rates of the lactic acid equation, which translates the crucial biophysical 
properties of a variable microenvironment.
The results obtained from the numerical simulations can be crucial for diagnostic and therapeutic 
applications, in particular for cancer forecasting and prevention.

The manuscript is organised as follows. 
In Section \ref{sect:GG}, we present the Gatenby--Gawlinski model, including its non-dimensionalised 
version and the analysis of invasion fronts, and we illustrate the passages leading to an approximation for
various biological regimes. 
In Section \ref{numerical algorithm}, we introduce the numerical algorithm 
 {based on a finite volume scheme, explaining also why we have chosen such approach 
(see Subsection \ref{finitevolumes}).
We end the Section by introducing the \emph{LeVeque--Yee formula}, i.e.,} 
a theoretical formula for approximating the propagation speed of the 
 {heteroclinic} traveling wave solutions. 
In Section \ref{homogeneous}, we simulate many different cases depending on the biological context.
Then, we make considerations on a possible homogenisation in the long run 
and we provide a comparison between analytical and numerical approximations.
 {Additional computations can be found in  {Appendix \ref{appendix}, at the very end of the paper.}
The paper is concluded by Section \ref{sect:conclusions} which contains a selection of possible 
further explorations and future perspectives.}

The {\tt SciLab} codes for reproducing the numerical simulations presented in this work
are available upon request to the authors.


\section{ {The Gatenby--Gawlinski model with heterogenous acid diffusion}}\label{sect:GG}

 {The Gatenby--Gawlinski model} is developed in order to reproduce cancer cells invasion within healthy tissue, 
starting from a stage in which the carcinogenesis has already happened. 
The attention is on the interaction between malignant and healthy cells occurring at the tumour-host interface, 
where a significant role is played by the lactic acid production and spreading, because of glycolytic metabolism
exploited by the tumour cells, see Fig.\ref{figura1}.

From a mathematical point of view, we consider the following system:
\begin{equation}\label{sistema GG}
	\left\{\begin{aligned}
	U_t &= \rho_1 U \big(1-{U}/{\kappa_1} \big)-\delta_1 U W\\
	V_t &= \rho_2 V \big(1-{V}/{\kappa_2}\big) + D_2 \big[\big(1-{U}/{\kappa_1}\big) V_x  \big]_x\\
	W_t &= \rho_3 V - \delta_3 W + \big[D_3(x) W_x \big]_x
	\end{aligned}\right.
\end{equation}
where $D_3=D_3(x)>0$ is a heterogeneous (positive) diffusion function. 
The boundary and initial conditions will be specified later on.

System \eqref{sistema GG} is a mathematical model of partial differential equations  {of evolutive type}
for the growth, diffusion and chemical action of tumour cells against the surrounding environment, 
where $U$ and $V$ are the healthy and tumour tissue concentrations with carrying capacities 
$\kappa_1$ and $\kappa_2$, respectively, and $W$ represents the excess of $H^+$ ions concentration 
induced by the tumour cells metabolism, 
 {translating the mechanism widely discussed in \cite{Warburg1, Warburg2}.}
A logistic-type growth is assumed in the first and second equations, with steady states $U=0$ and $V=0$ 
which are dynamically unstable, so that small perturbations drive the concentrations towards the stable
states $U=\kappa_1$ and $V=\kappa_2$ with growth rates $\rho_1$ and $\rho_2$, respectively, 
 {in the absence of other species.}

The diffusion rate for lactic acid in the third equation depends on the heterogeneous function  {$D_3=D_3(x)$} 
and, from a biological point of view, this accounts for the different physical composition of tissues
(see Fig.\ref{figura3}), which are not structurally homogeneous   {but} rather diversified according to 
their function and anatomical position \cite{histology}. 
On the other hand, it has been experimentally found that the diffusion rate for tumours typically
depends on the concentration of healthy cells in the surrounding  {micro-environment,} 
which is translated into the nonlinear (possibly degenerate) diffusion term in the second equation 
of model \eqref{sistema GG},  {see \cite{MaruAlmePoly12}}.

The destructive effects of acidity on the healthy tissue are described by the reaction term in the first equation, 
through a physical parameter $\delta_1$ whose values are crucial for the experimental analysis.
Finally, the (linear) production of lactic acid constant rate $\rho_3$, and its loss due to
deactivation kinetics with constant rate $\delta_3$, are included  {in} the third equation 
 {of system \eqref{sistema GG}.}

\begin{figure}\centering
	\subfloat[][epithelial tissue]{\includegraphics[width=.35\textwidth]{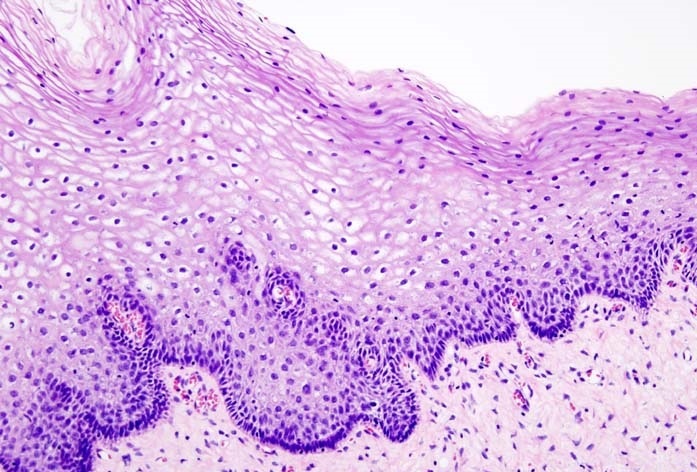}}\qquad
	\subfloat[][connective tissue]{\includegraphics[width=.35\textwidth]{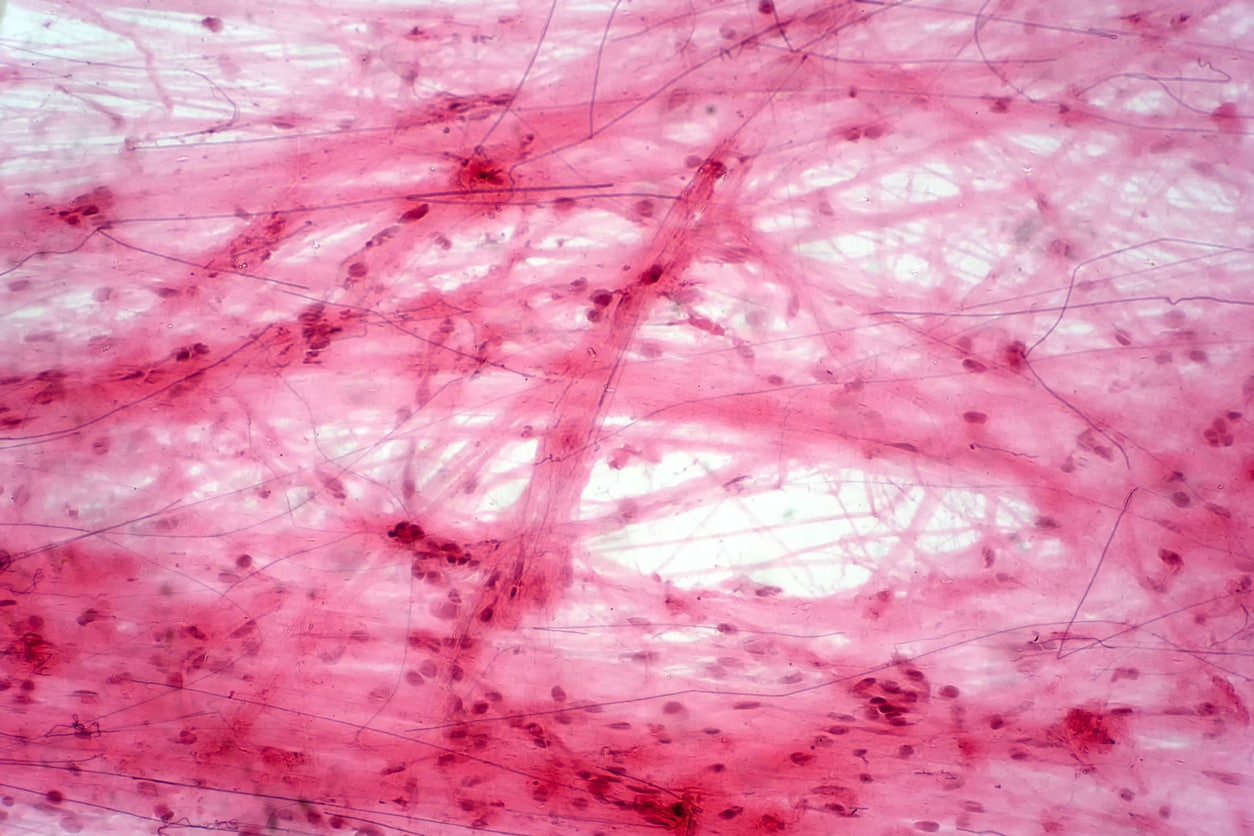}}\\
	\subfloat[][muscular tissue]{\includegraphics[width=.35\textwidth]{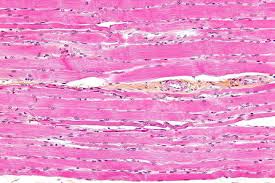}}\qquad
	\subfloat[][nervous tissue]{\includegraphics[width=.35\textwidth]{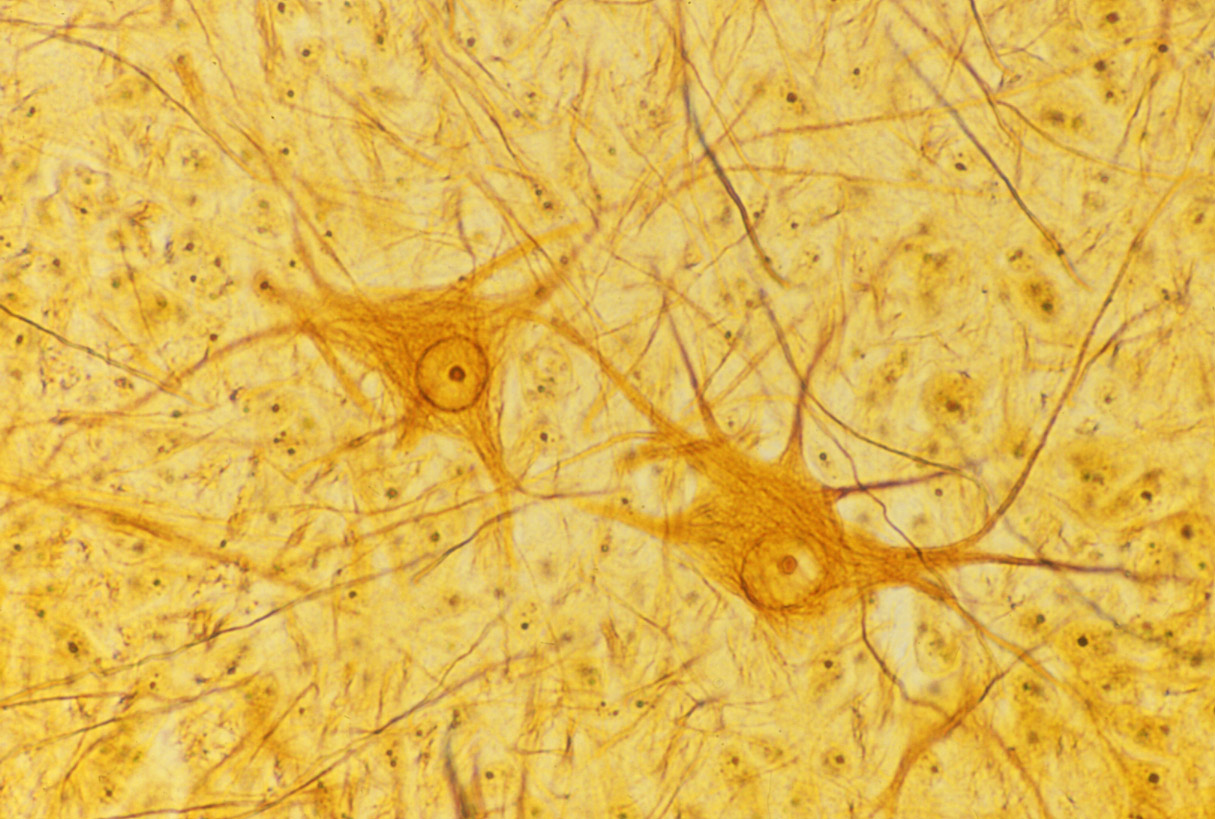}}
	\caption{\footnotesize{Microscope images of different types of biological tissues.}}
	\label{figura3}
\end{figure}

 {Next,} we adopt a standard adimensionalisation technique  and we consider $u$, $v$ 
and $w$ corresponding to rescaled concentrations of healthy tissue, tumour cells and excess 
of $H^+$ ions, respectively, which then satisfy the  {positivity requirement $u,v,w\geq 0$.} 
We introduce the following rescaled variables 
 {
\begin{equation*}
	u=\frac{U}{\kappa_1}, \qquad v=\frac{V}{\kappa_2}, \qquad w=\frac{\delta_3\,W}{\rho_3\kappa_2}\,,
\end{equation*}
and parameters
\begin{equation*}
	d=\frac{\delta_1\rho_3\kappa_2}{\delta_3\rho_1}, \qquad r=\frac{\rho_2}{\rho_1}, \qquad c=\frac{\delta_3}{\rho_1}.
\end{equation*}}
Because of biologically relevant considerations, 
 {the most important of those being to focalize on the effect of some specific biochemical 
signalling, \cite{McKiEtAl06} --specifically fluctuations of ionic concentration of pH--
we choose of  $D_3=D_3(x)$ as a space-dependent inhomogeneous diffusion function.}
Hence the adimensionalisation is performed with respect to the maximum intensity diffusion rate $\max D_3$, 
 {introducing the renormalised function $A=A(x)$ for the lactic acid, defined by
\begin{equation*}
  A(x):=\frac{D_3(x)}{\max D_3}\,.
\end{equation*}
Finally, introducing the rescaled independent variables
\begin{equation*}
	t \longmapsto \rho_1 t\,, \qquad x \longmapsto \sqrt{\frac{\rho_1}{\max D_3}}\,x\,,
\end{equation*}
we end up with
\begin{equation}\label{sistema GGbis}
	\begin{cases}
	u_t = u(1-u-dw)\\
	v_t = rv(1-v) + D\bigl[(1-u)v_x\bigr]_x\\
	w_t = c(v-w) + \bigl[A(x)w_x\bigr]_x
	\end{cases}
\end{equation}
where the additional (constant) parameter $D$ is defined by
\begin{equation*}
	D:=\frac{D_2}{\max D_3}>0\,.
\end{equation*}
From now on, for the sake of clarity, we will use the shortened version {\sf GG} model 
(or {\sf GG} system) to indicate \eqref{sistema GGbis}.}

Focusing on the biological interpretation of the {\sf GG} model \eqref{sistema GGbis},
the parameter $r$ denotes the reproduction rate of tumour tissue 
 {with respect to the reproduction rate of the normal one (explorations on uncontrolled tumor growth 
can be found in Subsection \ref{uncontrolled})} and $c$ 
 {the production induced by tumour tissue and spontaneous deactivation of $H^+$ ions.}
The variation of the tumour cells concentration described in the second equation is a consequence of its
intrinsic population dynamics and the diffusion mechanism also  {involves} the  {local} behaviour of the host tissue. 

As a matter of fact, there are two main properties characterizing the tumour dynamics, namely its aggressiveness
and its invasiveness, which are represented in system {\sf GG} \eqref{sistema GGbis} by the parameters $d$ and
$D$, respectively. 
For clinical application, the duality between highly aggressive but scarcely invasive and rapidly  {extending}
but less aggressive tumours constitutes a crucial issue \cite{bladder}.
The simulations illustrated in this work suggest that fast-growing tumours (i.e. large values of the growth rate $r$)
result also in enhanced invasiveness, hence these other properties are not independent. 
This fact is extremely interesting for applied scientists and clinical researchers since finding optimal strategies 
to control both these tumour parameters is still an open question \cite{Fasano 2}. 


\subsection{ {Propagation fronts for homogeneously diffusing acid}}

 {As stated before, one of the strengths of the {\sf GG} model \eqref{sistema GGbis} lies in its capability
to predict the presence of \emph{tumour-host hypocellular interstitial gap} between the invaded region and 
the invading one (Fig.\ref{figura4}), that is a separation layer practically depleted of cells between the healthy
and cancer populations and filled only with acid.
Such modelling prediction has been experimentally verified in both unfixed {\it in vitro} experiments 
and in flash-frozen tissues, by providing stronger evidence for this phenomenon \cite{GateGawl96}.
}
\begin{figure}[htb]\centering
	\includegraphics[scale=1.3]{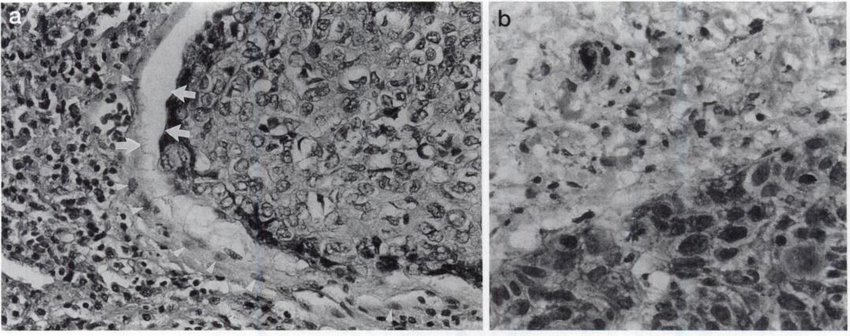}
	\caption{\footnotesize{Microscope images of tumour-host interstitial gap -- source \cite{GateGawl96}}} 
	\label{figura4}
\end{figure}

If $u \equiv 0$, namely in the absence of healthy tissue --or when the normal cells have 
been completely degraded by the effects of the tumour invasion-- then the second equation of system 
{\sf GG} \eqref{sistema GGbis}  {reduces to the classical \emph{Fisher--KPP equation} \cite{fisher, kpp37}}
\begin{equation}\label{kpp}
	v_t = rv(1-v) + Dv_{xx}\,,
\end{equation}
 {Equation \eqref{kpp}} exhibits  {infinitely many} front-type solutions connecting 
the equilibrium states given  {by $v=v_-:=1$ and $v=v_+:=0$.
The partial differential equation \eqref{kpp}} occurs in phase transition problems; 
in particular, Ronald Fisher proposed it in the context of population dynamics to describe
the spatial spread of an advantageous allele and explored its travelling wave solutions.
Indeed, it can be proven that for every wave speed $\theta \ge 2\sqrt{rD}$ there exists 
 {a \emph{propagation front} (or an \emph{invasion front})}
which switches from the  {stable} equilibrium state $v_- =1$ to the  {unstable} equilibrium $v_+ = 0$
 {with a monotone decreasing profile moving toward the right.} 
On the other hand, no
 {propagation front exists for $\theta<2\sqrt{rD}$ with a monotone profile, \cite{hadeler}.
Non-monotone profiles may exist but they are not biologically significant due to the loss of positivity,
a crucial property for the density function $v$.}

The travelling waves solutions to the {\sf GG} model \eqref{sistema GGbis} are special solutions 
of the form $\phi(x,t)=\Phi(x-\theta t)$\,, where $\theta$ is a (constant) \emph{wave speed} and the profile $\Phi$
propagates along the real line, with asymptotic states usually renormalised according to the problem context. 
 {By extension from Fisher--KPP equation \eqref{kpp},}
if $\Phi$ is required to satisfy boundary conditions of the form $\Phi(-\infty)=\phi_-$ and $\Phi(+\infty)=\phi_+$,
then the solution is called a 
 {propagation front (or, as before, an invasion front).}

 {Firstly, to avoid unnecessary analytical difficulties, 
we consider the case of a constant diffusion rate $A$ constant in the equation for the lactic acid.}
For analyzing the front-type solutions, we use the change of variables $z=x-\theta t$ 
and we make a slight abuse of notation by denoting
\begin{equation*}
	u(z)=u(x,t)\,, \qquad v(z)=v(x,t)\,, \qquad w(z)=w(x,t)\,.
\end{equation*}
By substitution into the {\sf GG} system \eqref{sistema GGbis}, we obtain
\begin{equation}\label{sistema GG TW1}
	\begin{cases}
	\theta u' + u(1-u-dw) = 0\\
	\theta v' + rv(1-v) + D\bigl[(1-u)v'\bigr]' = 0\\
	\theta w' + c(v-w) + Aw''= 0
	\end{cases}
\end{equation}
where the derivatives are referred to the variable $z$. 
In order to be consistent with the biological model, we consider $\theta$
and the following asymptotic conditions:
\begin{equation}\label{eq0}
	(u,v,w)(-\infty) = (0,1,1)\quad \text{and}\quad (u,v,w)(+\infty) = (1,0,0)\quad \text{for}\quad d\ge 1
\end{equation}
or
\begin{equation}\label{eq1}
	(u,v,w)(-\infty) = (1-d,1,1)\quad \text{and}\quad (u,v,w)(+\infty) = (1,0,0)\quad \text{for}\quad 0<d<1
\end{equation}
together with the asymptotic conditions $(u',v',w')(\pm \infty) = 0$ which correspond to the steady states 
of the  {underlying} dynamical system \eqref{sistema GG TW1}. 
In particular, if $d\ge 1$\,, due to the strong acidity effects, then the solutions \eqref{eq0} describe a process 
called \emph{homogeneous invasion} where the healthy tissue is completely destroyed behind the tumour front. 
Whilst if $0<d<1$\,, then the solutions \eqref{eq1} describe another situation called \emph{heterogeneous invasion}
with a residual concentration of healthy tissue given by the value $1-d$\,. 
Moreover, the profiles $u$,\, $v$ and $w$ are assumed to be monotonic, increasing in the case of $u$ 
and decreasing in the case of $v$ and $w$ (refer to Section \ref{homogeneous} for details).

\begin{remark}
 {At present day only partial rigorous results have been obtained about existence of propagating fronts 
for the ordinary differential system \eqref{sistema GG TW1} even in the case of constant $A$, \cite{DHMR22}.
To most of our knowledge, inspired by \cite{MascMoscSime24}, rigorous results for some reduction 
consisting of the assumption that (normalized) tumor cells and (normalized) acid are instantaneously 
at equilibrium have been proven in \cite{GallMasc22}.
Hence, for the complete {\sf GG} model \eqref{sistema GGbis}, we content ourselves with the 
numerical evidence of existence of propagating fronts taken from \cite{McGillen, Pierfi}.
Additionally, some formal expansions relative to the presence/absence of the interstitial gap
can be found in \cite{Fasano}.

Moreover, if for {\sf GG} model \eqref{sistema GGbis} a result analogous to the one proven 
for Fisher--KPP equation \eqref{kpp}, that is the existence of a minimal positive speed $\theta_\ast$
such that there exist a propagating front for any $\theta\geq \theta_\ast$, then it can be seen that 
the following implication holds
\begin{equation*}
	u\leq 1\qquad\Rightarrow\qquad 
	\theta_\ast^{GG}\leq \theta_\ast^{FKPP}=2\sqrt{rD}\,,
\end{equation*}
where $\theta_\ast^{FKPP}$ and $\theta_\ast^{GG}$ denote the minimal speed of Fisher--KPP equation
\eqref{kpp} and {\sf GG} system \eqref{sistema GGbis}, respectively.
Nevertheless, to our knowledge, no rigorous result of this type is known concerning {\sf GG} model.}
\end{remark}

 {The extension to the case of heterogeneous diffusive functions $A$, even under the simplifying assumption 
of a periodic profile $A$, is much more complicated because differentiation with respect to $z$ is not possible and 
additional dependency has to be taken into account, \cite{Derk15,Xin91}.
In particular, we observe the fact that in such a case invariance with respect to translation is lost.}
 

\section{ {Finite-volume method in one space dimension}}\label{numerical algorithm}

Version \eqref{sistema GGbis} of the {\sf GG} model allows to operate with fewer parameters $d$, $r$, $D$ 
and $c$, thus reducing their original range and coping with scaled functions $u$, $v$, $w$ and $A=A(x)$.

For the numerical investigation, the experimental domain is assumed to be the one-dimensional interval 
$[-L, L]$ with $L>0$ with homogeneous Neumann boundary conditions 
 {corresponding to the hypothesis isolated from its complementary set.
Moreover,}
we adopt a discretization strategy based on the \emph{finite volume method}, which guarantees
consistency and stability in terms of closeness to the physical properties of the model \cite{wesseling}. 
The numerical algorithm is semi-implicit in time and employs a non-uniform mesh for the  {spatial}
discretization.

Let $C_i=[x_{i-\frac{1}{2}},x_{i+\frac{1}{2}})$ be the
finite volume centred in $x_i=\frac{1}{2}(x_{i-\frac{1}{2}}+x_{i+\frac{1}{2}})$, 
for $i=1,2,\dots,N$, where $N$ is the number of vertices. 
Let $\Delta x_i=|x_{i+\frac{1}{2}}-x_{i-\frac{1}{2}}|$ be the variable cells' size, hence 
$|x_i-x_{i-1}|=\tfrac{1}{2}(\Delta x_{i-1}+\Delta x_i)$ is the length of an interfacial
interval (see Fig.\ref{figura2}). 

\begin{figure}[htb]
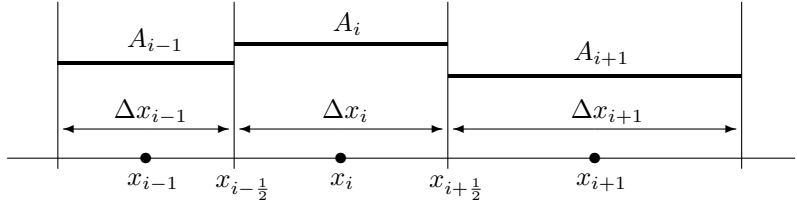

{\;
\put(82,18){\line(0,1){63}}
\put(63,22){\line(1,0){295}}
\put(115,22){\circle*{4}} 
\put(108,12){$x_{i-1}$}
\put(115,33){\vector(-1,0){31}}
\put(115,33){\vector(1,0){31}}
\put(103,37){$\Delta x_{i-1}$}
\put(141,11){$x_{i-\frac12}$}
\put(148,18){\line(0,1){63}}
\put(188,22){\circle*{4}} 
\put(184,12){$x_i$}
\put(188,33){\vector(-1,0){38}}
\put(188,33){\vector(1,0){38}}
\put(181,37){$\Delta x_i$}
\put(228,18){\line(0,1){63}}
\put(283,22){\circle*{4}}
\put(276,12){$x_{i+1}$}
\put(338,18){\line(0,1){63}}
\put(283,33){\vector(-1,0){53}}
\put(283,33){\vector(1,0){53}}
\put(274,37){$\Delta x_{i+1}$}
\put(221,11){$x_{i+\frac12}$}
\linethickness{0.4mm}
	\put(82,58){\line(1,0){66}}\put(108,64){$A_{i-1}$}
	\put(148,65){\line(1,0){80}}\put(185,70){$A_i$}
	\put(228,53){\line(1,0){110}}\put(275,59){$A_{i+1}$}
}
\caption{\footnotesize{Piecewise constant reconstruction on nonuniform mesh.}}
\label{figura2}
\end{figure}
We build a piecewise constant approximation of the diffusion function $A=A(x)$ by means
of its \emph{integral cell-averages}, namely
\begin{equation*}
	A_i = \frac1{\Delta x_i} \int_{C_i} A(x)\,dx \,,
\end{equation*}
and we perform the same projection also for the solution to the {\sf GG} system \eqref{sistema GGbis} by setting
\begin{equation*}
	u_i(t) = \frac1{\Delta x_i} \int_{C_i} u(x,t)\,dx\,,\quad
	v_i(t) = \frac1{\Delta x_i} \int_{C_i} v(x,t)\,dx\,,\quad
	w_i(t) = \frac1{\Delta x_i} \int_{C_i} w(x,t)\,dx\,.
\end{equation*}
Firstly, we consider the equation for the healthy cells concentration $u$
and its finite volume integral version given by
\begin{equation*}
	\frac{1}{\Delta x_i}\int_{C_i} u_t\,dx\, = \frac{1}{\Delta x_i} \int_{C_i} u(1-u-d\,w)\,dx\,,
\end{equation*}		
so that we obtain the following approximation
\begin{equation}\label{eq discr 1}
	u_i' = u_i(1-u_i- d\,w_i)\,,
\end{equation}
where the time dependence is dropped for shortness 
(with an abuse of notation for the symbols of the numerical variables).
Then, we consider the equation for the tumour cells concentration $v$ and its finite
volume integral version given by
\begin{equation*} 
	\frac{1}{\Delta x_i} \int_{C_i} v_t\,dx\, = \frac{r}{\Delta x_i} \int_{C_i} v(1-v)\,dx\, 
	+ \frac{D}{\Delta x_i} \int_{C_i}\, \big[ (1-u) v_x \big]_x\,dx \,.
\end{equation*}
Hence, we need a specific approach for the discretization of the nonlinear diffusion term. 
We proceed by evaluating the differential term at the mesh interfaces 
 {(as suggested by the Divergence Theorem which, the present case being one-dimensional, 
reduces to the Fundamental Theorem of Calculus)}
as follows
\begin{equation}\label{eq discr 2}
	\begin{aligned}
	\frac{D}{\Delta x_i} \bigg[ & \frac{(1-u_i)\Delta x_i+(1-u_{i+1})\Delta x_{i+1}}
		{\Delta x_i+\Delta x_{i+1}} \cdot \frac{v_{i+1}-v_i}{\frac{\Delta x_i}{2}+\frac{\Delta x_{i+1}}{2}}\\
	& -\, \frac{(1-u_{i-1})\Delta x_{i-1}+(1-u_i)\Delta x_i}{\Delta x_{i-1}+\Delta x_i}
		\cdot \frac{v_i-v_{i-1}}{\frac{\Delta x_{i-1}}{2}+\frac{\Delta x_i}{2}} \bigg]\,,
	\end{aligned}
\end{equation}
where we made use of interfacial differences to approximate the derivatives, and the other interfacial 
quantities are approximated by means of weighted averages, whose weights are the size of the adjacent 
finite volumes, so that $\Delta x_{i-1}/\Delta x_i$ and $\Delta x_i/\Delta x_{i+1}$ are employed 
at the interfaces $x_{i-\frac{1}{2}}$ and $x_{1+\frac{1}{2}}$, respectively. 
Therefore, we obtain
\begin{equation}\label{eq discrr}
	\begin{aligned}
	v'_i = r\,v_i(1-v_i) + \frac{D}{\Delta x_i} \bigg[ & \frac{(1-u_i)\Delta x_i+(1-u_{i+1})\Delta x_{i+1}}
		{\Delta x_i+\Delta x_{i+1}} \cdot \frac{v_{i+1}-v_i}{\frac{\Delta x_i}{2}+\frac{\Delta x_{i+1}}{2}}\\
	& -\,\frac{(1-u_{i-1})\Delta x_{i-1}+(1-u_i)\Delta x_i}{\Delta x_{i-1}+\Delta x_i}
		\cdot \frac{v_i-v_{i-1}}{\frac{\Delta x_{i-1}}{2}+\frac{\Delta x_i}{2}} \bigg]
	\end{aligned}
\end{equation}
Finally, we consider the equation for the lactic acid concentration $w$ and its finite
volume integral version given by
\begin{equation*}
	\frac{1}{\Delta x_i} \int_{C_i} w_t\,dx\, = \frac{c}{\Delta x_i} \int_{C_i} (v - w)\,dx
		+ \frac{1}{\Delta x_i} \int_{C_i} \bigl[ A(x)\,w_x \bigr]_x\,dx\,,
\end{equation*}
and, by proceeding as for \eqref{eq discr 2},   {we end up with}
\begin{equation}\label{eq discr 3}
	\begin{aligned}
	w'_i = c\,(v_i-w_i) + \frac{1}{\Delta x_i} \bigg( & \frac{A_i\Delta x_i+A_{i+1}\Delta x_{i+1}}
		{\Delta x_i+\Delta x_{i+1}} \cdot \frac{w_{i+1}-w_i}{\frac{\Delta x_i}2+\frac{\Delta x_{i+1}}2}\\
	& -\,\frac{A_{i-1}\Delta x_{i-1}+A_i\Delta x_i}{\Delta x_{i-1}+\Delta x_i} \cdot \frac{w_i-w_{i-1}}
		{\frac{\Delta x_{i-1}}2+\frac{\Delta x_i}2} \bigg)\,.
	\end{aligned}
\end{equation}

 {\subsection{The case of a uniform mesh}}

For the case of a uniform mesh with $\Delta x_i=\Delta x$ for all $i=1,2,\dots,N$, 
the semi-discrete system  {can be} simplified as follows
 {
\begin{equation}\label{eq discr 4bis}
	\left\{\begin{aligned}
	u'_i &= u_i(1-u_i - dw_i)\\
	v'_i &= r\,v_i(1-v_i) + \frac{D}{\Delta x}\bigg( \frac{2-u_i-u_{i+1}}{2} \cdot \frac{v_{i+1}-v_i}{\Delta x}
	\,- \,\frac{2-u_{i-1}-u_i}{2}.\frac{v_i-v_{i-1}}{\Delta x} \bigg)\\
	w'_i &= c\,(v_i-w_i) + \frac{1}{\Delta x}\bigg( \frac{A_i+A_{i+1}}{2} \cdot \frac{w_{i+1}-w_i}{\Delta x}
	\,- \,\frac{A_{i-1}+A_i}{2} \cdot \frac{w_i-w_{i-1}}{\Delta x} \bigg)
	\end{aligned}\right.
\end{equation}
}

\begin{remark}\label{remark2}
The use of the arithmetic average for dealing with diffusive equations is not standard in numerical analysis,
and the \emph{harmonic mean} is instead preferable for practical reasons. 
Hence, for example, we can consider the following alternative version of 
 {the third equation in \eqref{eq discr 4bis}} for the diffusion function $A$ of the lactic acid concentration
\begin{equation}\label{eq discr 6 harm}
	w'_i = c\,(v_i-w_i) + \frac{1}{\Delta x}\bigg( \frac{2A_iA_{i+1}}{A_i+A_{i+1}} 
	\cdot \frac{w_{i+1}-w_i}{\Delta x}\,- \,\frac{2A_{i-1}A_i}{A_{i-1}+A_i} \cdot \frac{w_i-w_{i-1}}{\Delta x} \bigg)\,.
\end{equation}
Nonetheless, it is important to observe that choosing the harmonic mean leads to numerical issues when
dealing with degenerate diffusion functions  so that the use of harmonic mean for diffusive
equations with singularities can be pursued with advanced techniques and artifacts, \cite{maddix}.
For instance, this is the case of 
 {the second equation in \eqref{eq discr 4bis}} 
for the tumour cells density, as soon as the (nonlinear) diffusion function $1-u$ becomes null somewhere. 
Then, the diffusion matrix which defines the numerical scheme becomes singular and, unfortunately, 
that is precisely the case of our simulations. 
Thus, we can not straightforwardly apply the harmonic average and we finally recurred to the arithmetic
mean, although this issue will be better explored in the future.
\end{remark}

For the time discretization of   {the ordinary differential system \eqref{eq discr 4bis}}
we employ a semi-implicit strategy by considering a time step $\Delta t_n=\abs{t^{n+1}-t^n}$
for $n=0,1,\dots$.
In particular, the reaction terms are treated explicitly, whilst the differential terms on the right-hand side can be 
approximated implicitly, so that we  {end up with the system}
\begin{equation}\label{eq discr 7}
	\left\{ \begin{aligned}
	u_i^{n+1} &= u_i^n + \Delta t \,u_i^n(1-u_i^n-d\,w_i^n)\\
	v_i^{n+1} &= v_i^n + r\,\Delta t\,v_i^n(1-v_i^n) + D \frac{\Delta t}{\Delta x}
		\bigg(\frac{2-u^{n+1}_i-u^{n+1}_{i+1}}{2}\cdot \frac{v^{n+1}_{i+1}-v^{n+1}_i}{\Delta x}\\
	& \hspace{135pt} - \frac{2-u^{n+1}_{i-1}-u^{n+1}_i}{2} \cdot \frac{v^{n+1}_i-v^{n+1}_{i-1}}{\Delta x} \bigg)\\
	w_i^{n+1} &= w_i^n + c\,\Delta t(v_i^n-w_i^n) + \frac{\Delta t}{\Delta x} \bigg( \frac{A_i+A_{i+1}}{2}
		\cdot \frac{w_{i+1}^{n+1}-w_i^{n+1}}{\Delta x}\\ &  \hspace{125pt} - \frac{A_{i-1}+A_i}{2}
		\cdot \frac{w_i^{n+1}-w_{i-1}^{n+1}}{\Delta x} \bigg)
	\end{aligned} \right.
\end{equation}
 {valid for the case $\Delta t_n=\Delta t$ with $\Delta t>0$ some fixed number.}
We finally impose Neumann-type boundary conditions $v_1^n=v_2^n$, $w_1^n=w_2^n$ and
$v_{N}^n=v_{N-1}^n$, $w_{N}^n=w_{N-1}^n$, for $n=1,2,\dots$
This implicit-explicit ({\sf IMEX}) mixed approach allows to make less expensive choices for the time step,
compared to fully explicit algorithms which would be heavily conditioned by stability restrictions and
consequently slower in computational time \cite{Quarteroni}.

The theoretical properties of the numerical algorithm are proven in \cite{ScanMascSime20} and, in particular, 
the consistency analysis (for regular solutions) which concludes that the algorithm is order 1 in time and order 
1 in space --despite the parabolic character of the dynamical system-- because of the inhomogeneity of the 
diffusion rates for the lactic acid equation. 
Nevertheless, a numerical study of the convergence error reveals an increased order of convergence, thus
suggesting to undertake further investigation.

\begin{remark}
In practical applications, we need to impose boundary conditions also for $u$, namely $u_1^n=u_2^n$ and
$u_{N-1}^n=u_N^n$, even if the mathematical model for the healthy tissue is indeed a purely ordinary
differential equation.
In order to be consistent with the biological context, a diffusion term should be considered also for the normal
cells population, but the diffusivity rate for $u$ would be extremely slower than those for $v$ and $w$, hence 
it can be omitted in the {\sf GG} model \eqref{sistema GGbis}.
Nevertheless, this simplification has to be taken into account when imposing numerical boundary conditions
for all three equations to efficiently solve the fully discrete system \eqref{eq discr 7}.
\end{remark}


 {\subsection{Discussion on the choice of finite volume algorithms}}\label{finitevolumes}

The diffusion term in the third equation of \eqref{sistema GGbis} can be rewritten as
\begin{equation}\label{splitting}
	\bigl[A(x)\,w_x\bigr]_x = A'(x)\,w_x + A(x)\,w_{xx}\,,
\end{equation}
thus revealing the typical characteristics of a transport mechanism induced by the non-homogeneity
of the diffusion rate function $A$. 
Moreover, the third equation in \eqref{eq discr 4bis} can also be rewritten as 
\begin{equation*}
	w'_i = c\,(v_i-w_i) + \frac{1}{2 \Delta x^2} \bigl[ A_{i+1}(w_{i+1}-w_i) 
		- A_{i-1}(w_i-w_{i-1}) + A_i(w_{i+1}-2w_i+w_{i-1}) \bigr]\,,
\end{equation*}
which corresponds to the analytical splitting of the derivative  {in \eqref{splitting}.}
Analogous considerations can be  {made} for the nonlinear diffusion term in 
 {the second equation in \eqref{eq discr 4bis}.}
However, this alternative is not useful for the numerical approach, as we will discuss 
 {in more detail} later on.

There are essentially two methods for the discretization of derivatives, namely finite differences
and finite volumes \cite{Quarteroni}. 
In the one dimensional setting, these strategies typically lead to the same results, but the procedures
through which they are obtained are significantly different. 
In order to explain this point, we analyze the two approaches applied to the {\sf GG} system \eqref{sistema GGbis} 
and, in particular, to the equation for $w$ by focusing on the differential term $\bigl[A(x)\,w_x\bigr]_x$ with 
$A$ replaced by its piecewise constant projection on the spatial mesh.

Firstly, we rewrite this term by splitting the derivatives as done in \eqref{splitting} and 
the finite volume integral formulation of the right-hand side is
\begin{equation}\label{eq finitevolume1}
	\frac{1}{\Delta x} \int_{x_{i-\frac{1}{2}}}^{x_{i+\frac{1}{2}}}\, A'(x) w_x \, dx 
	+\frac{A_i}{\Delta x} \int_{x_{i-\frac{1}{2}}}^{x_{i+\frac{1}{2}}} \, w_{xx} \,dx\,,
\end{equation}
where we have considered a uniform spatial mesh. 
We notice that we took out $A_i$ from the second integral because it is constant inside each finite volume 
 {$\bigl[x_{i-\frac{1}{2}},x_{i+\frac{1}{2}}\bigr]$.} 
We discretize the second integral as follows
\begin{equation}\label{eq discret}
	\frac{A_i}{\Delta x}\bigg(\frac{w_{i+1}-w_i}{\Delta x}-\frac{w_i-w_{i-1}}{\Delta x} \bigg)
	=\frac{A_i}{\Delta x^2}\big(w_{i+1}-2w_i+w_{i-1}\big)\,,
\end{equation}
since the derivatives are approximated at the interfaces $i-\frac{1}{2}$ and $i+\frac{1}{2}$, respectively.\\
For the first integral in \eqref{eq finitevolume1} we can take, for instance,
\begin{equation*}
	\frac{1}{2} \frac{A_i-A_{i-1}}{\Delta x} \cdot \frac{w_i-w_{i-1}}{\Delta x} 
	+ \frac{1}{2} \frac{A_{i+1}-A_i}{\Delta x} \cdot \frac{w_{i+1}-w_i}{\Delta x}\,,
 \end{equation*}
where the product of derivatives is approximated with the average values at the interfaces
$i-\frac{1}{2}$ and $i+\frac{1}{2}$, respectively, and $\Delta x$ in the denominator is simplified 
with the one resulting from the length of the  {entire interval}.
Finally, we can come back to recover the original discretisation of  
 {$\bigl[A(x)\,w_x\bigr]_x$ as done in the third equation in system \eqref{eq discr 7}.}

Now let us apply the finite differences for discretizing $A\,w_{xx}$, so that we obtain the same result
as in the right-hand side of \eqref{eq discret}. 
The difference occurs for the product $A'(x)w_x$, which is typically discretized by centered finite differences as follows
\begin{equation*}
	\frac{A_{i+1}-A_{i-1}}{2 \Delta x} \cdot \frac{w_{i+1}-w_{i-1}}{2\Delta x}\,,
\end{equation*}
and, therefore, it seems not to be possible to recover from the sum of these two parts any (consistent) finite difference 
discretization of the  {differential operator in divergence form \eqref{splitting}}, 
differently from the finite volumes approach.\\


\subsection{Space-averaged estimate of the propagation speed}

Once the numerical framework has been established, the natural step to proceed with experiments 
and simulations consists  {of} defining a wave speed estimation for the numerical solution. 
 {Let us stress that the speed of propagation of a front is a crucial and macroscopically measurable
ingredient not only from a theoretical and numerical point of view, but also for furnishing clinically
reliable predictions \cite{DosSantosEtAl24}.}

 {To provide} a numerical approximation for the wave speed at time $t^n$, we employ the 
 {\emph{LeVeque--Yee formula}, a}
space-averaged estimate, originally proposed in \cite{leveque}, and successfully applied to 
the case of a reactive version of the Goldstein--Kac model for correlated reaction-diffusion systems 
in \cite{lattanzio, lattanzio 2}. 

We briefly derive the main analytical concepts behind its formulation:
given an increment $h$, let $\phi$ be a differentiable function having constant asymptotic 
states $\phi_\pm$ (with $\phi_-\neq\phi_+$). 
Then we can write
\begin{equation*}
	\begin{aligned}
	\int_{\mathbb{R}}\,\big[\phi(\xi+h)-\phi(\xi) \big]\,d\xi
		=h\int_{\mathbb{R}}\int_{0}^{1}\,\phi'(\xi+\eta h)\,d\eta d\xi
		=h\int_{0}^{1}\!\!\int_{\mathbb{R}}\,\phi'(\mu)\,d\mu \, d\eta
		=h\big(\phi_+-\phi_- \big)\,.
	\end{aligned}
\end{equation*}
Setting  {$h=-\theta\Delta t$}, we deduce the following integral equation for the wave speed
 {\begin{equation}\label{wspeed}
	\theta=\frac{1}{\Delta t}\,\int_{\mathbb{R}} \frac{\phi(\xi)-\phi(\xi-\theta\Delta t)}{\phi_--\phi_+}\,d\xi.
\end{equation}}
We recall that invasion fronts are computed by imposing the change of variable  {$z=x-\theta t$}
in the {\sf GG} system \eqref{sistema GGbis}. 
For instance, if we focus on the invasion front with speed $s>0$ of the tumour cells density, hence 
 {$\phi=\phi(z)$} 
corresponds to $v(x,t)$ in the original variables and the  {discretized version of the space averaged 
wave speed estimate} over a uniform spatial mesh at time $t^n$ is given by the LeVeque--Yee formula, 
 {that is
\begin{equation}\label{levequeyee}
	\theta^n=\frac{\Delta x}{\Delta t}\cdot\frac{1}{v_--v_+}\,\sum_{i=1}^{N}\big(v_i^{n+1}-v_i^n\big),
\end{equation}
with $N$ is the number of vertices and $v_\pm$ the stationary states of $v$ 
in the {\sf GG} system \eqref{sistema GGbis}. 
Let us stress that $n\mapsto \theta^n$ defines a numerical sequence of real numbers;
therefore the choice of the value of the natural number $n$ is crucial and deserves some caution.
}

It is important to underline the strength of estimate \eqref{levequeyee}  {relying}
on its independence from the dynamics of the solution provided by the {\sf GG} model \eqref{sistema GGbis}, 
thus being always numerically computable  {once chosen the numerical scheme.} 


 {\section{Experimental validation of the numerical scheme}}\label{homogeneous}

 {The numerical computations proposed in this section provide support for strategies to efficiently struggle 
against cancer spread in both experimental and clinical applications. 
In principle, one may interfere with the tissue inhomogeneity to slow down the tumour front or rather act
pharmacologically only on regions where the cancer cells spread faster. 
Moreover, the fact that some specific acid diffusion profiles lead to inefficient solutions could explain why
some tissue seems to be naturally tumour-free. 
For instance, soft tissue sarcomas of muscles, nerves and blood vessels are very rare, maybe because
of the fibrous tissues structure which obstructs the destructive acid infiltration \cite{hudzik}.}

 {In what follows --after a short discussion on the homogeneous diffusion case for the acid equation--
we focus on two main prototypical cases: piecewise constant diffusion function with a single jump and 
periodic diffusion function with a single frequency.}

To validate the numerical algorithm, an extensive set of simulations has been performed in 
\cite{Pierfi, ScanMascSime20} using the scheme \eqref{eq discr 7}, with a uniform spatial mesh, 
and the numerical results are consistent with those presented in the literature \cite{McGillen} and
also with the analytical approach illustrated in \cite{Fasano}.

The experiments are carried out with the parameters available in \cite{Fasano}, as listed in Table \ref{tavola1}. 
Of course, the choice of parameters is crucial and different choices could be considered 
 {(as done in \cite{McGillen}).}
We assume $T$ as the final time and $L$  {as} the semi-length of the space interval $[-L,L]$, 
while the spatio-temporal mesh is built by fixing $\Delta x=0.005$ and $\Delta t=0.01$.

\begin{table}[htp]\centering
\begin{tabular}{c|c|c|c|c|c|c}
	$d$                       & $r$   &  $D$            &  $c$ & $A$ & $L$  & $T$      \\	\hline
	$0.5,1.5,2.5,12.5$    & 1   &  $4\cdot10^{-5}$  &  70  & 1 & 1  & 20 
\end{tabular}
\caption{\footnotesize{Numerical values for the simulation parameters -- source \cite{Fasano}}}
\label{tavola1}
\end{table}

\begin{figure}[hbt]\centering
	\subfloat[][initial profiles]{\includegraphics[width=.35\textwidth]{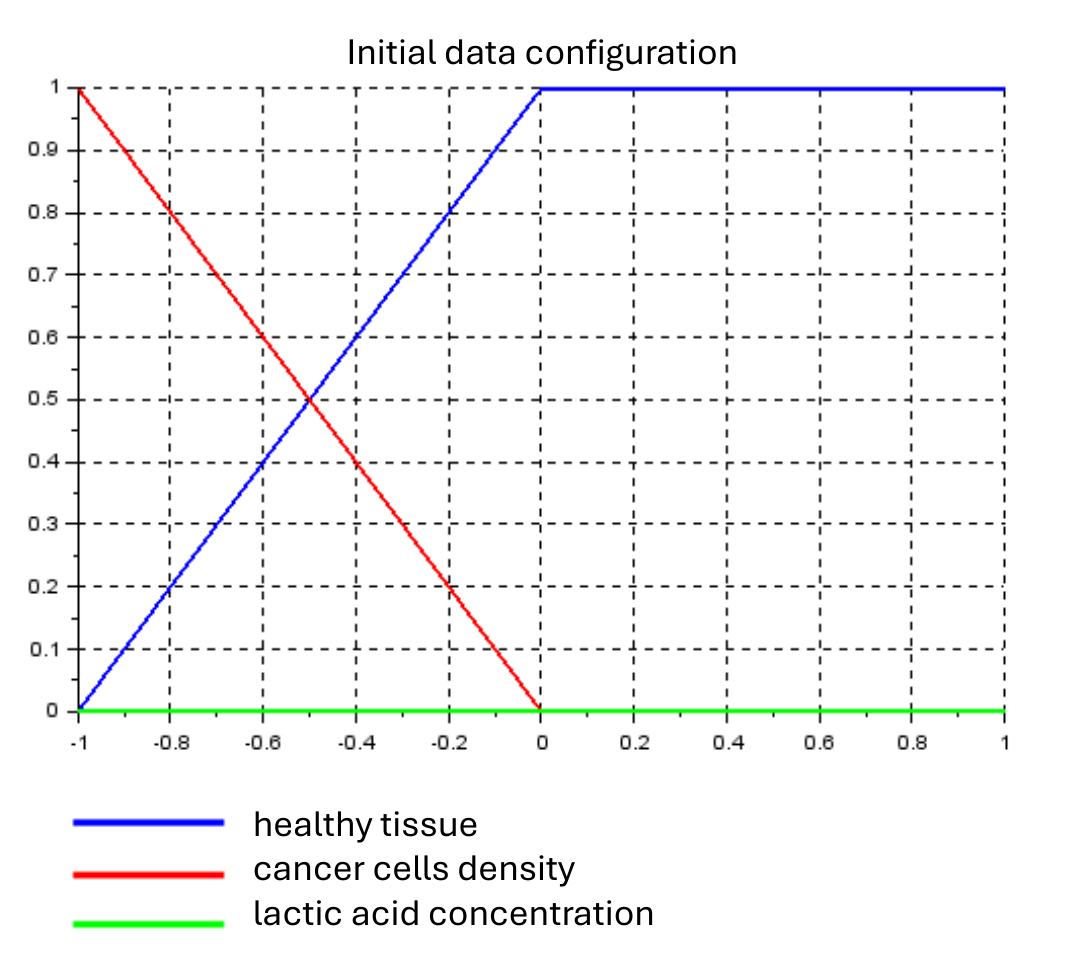}}\qquad
	\subfloat[][heterogeneous invasion, $d=0.5$]{\includegraphics[width=.35\textwidth]{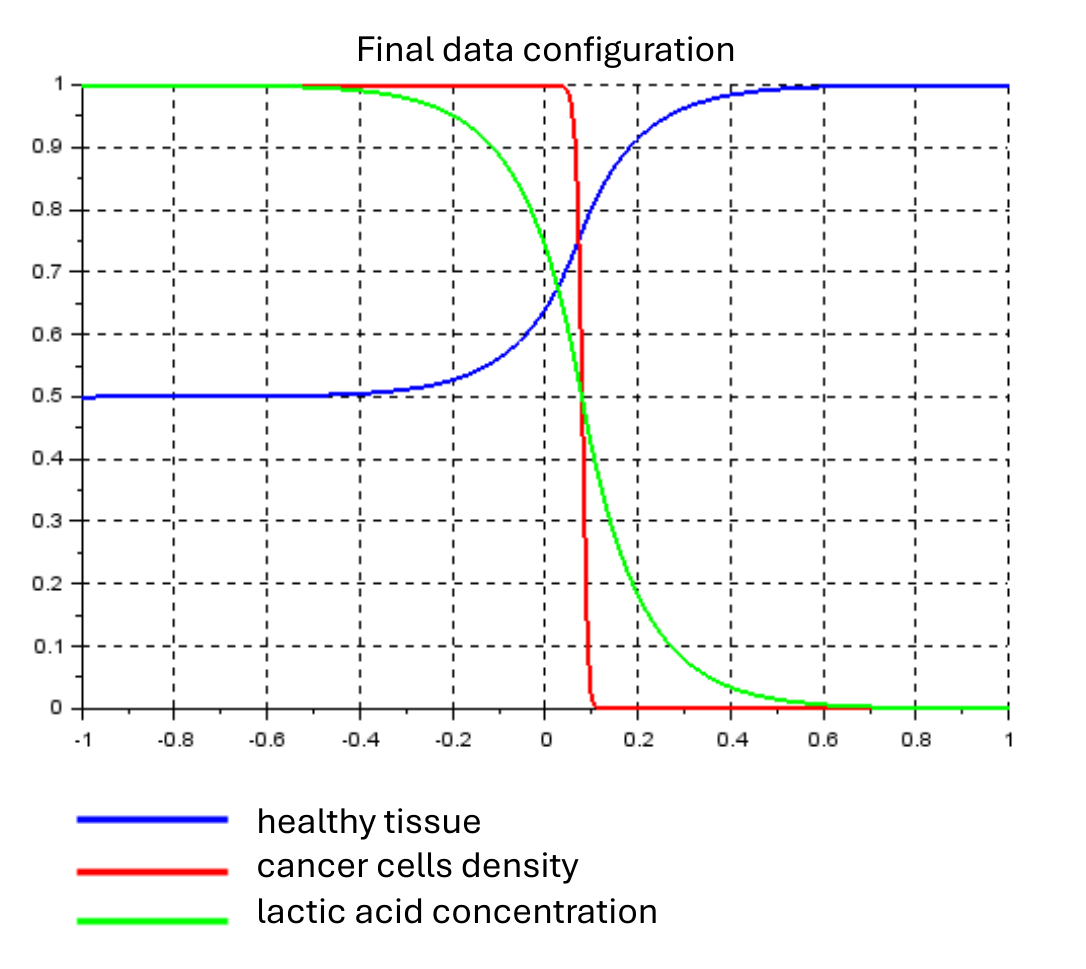}}\\
	\subfloat[][hybrid configuration, $d=2.5$]{\includegraphics[width=.35\textwidth]{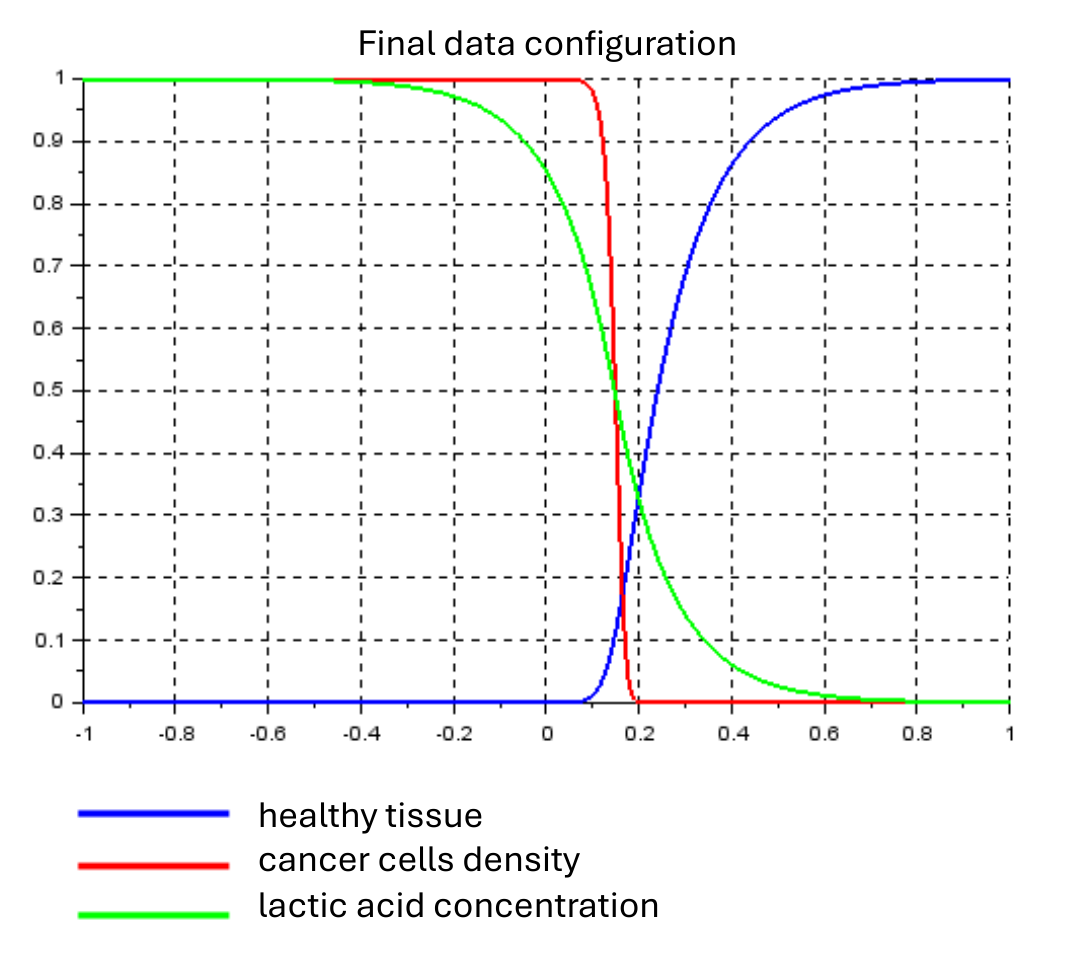}}\qquad
	\subfloat[][homogeneous invasion, $d=12.5$]{\includegraphics[width=.35\textwidth]{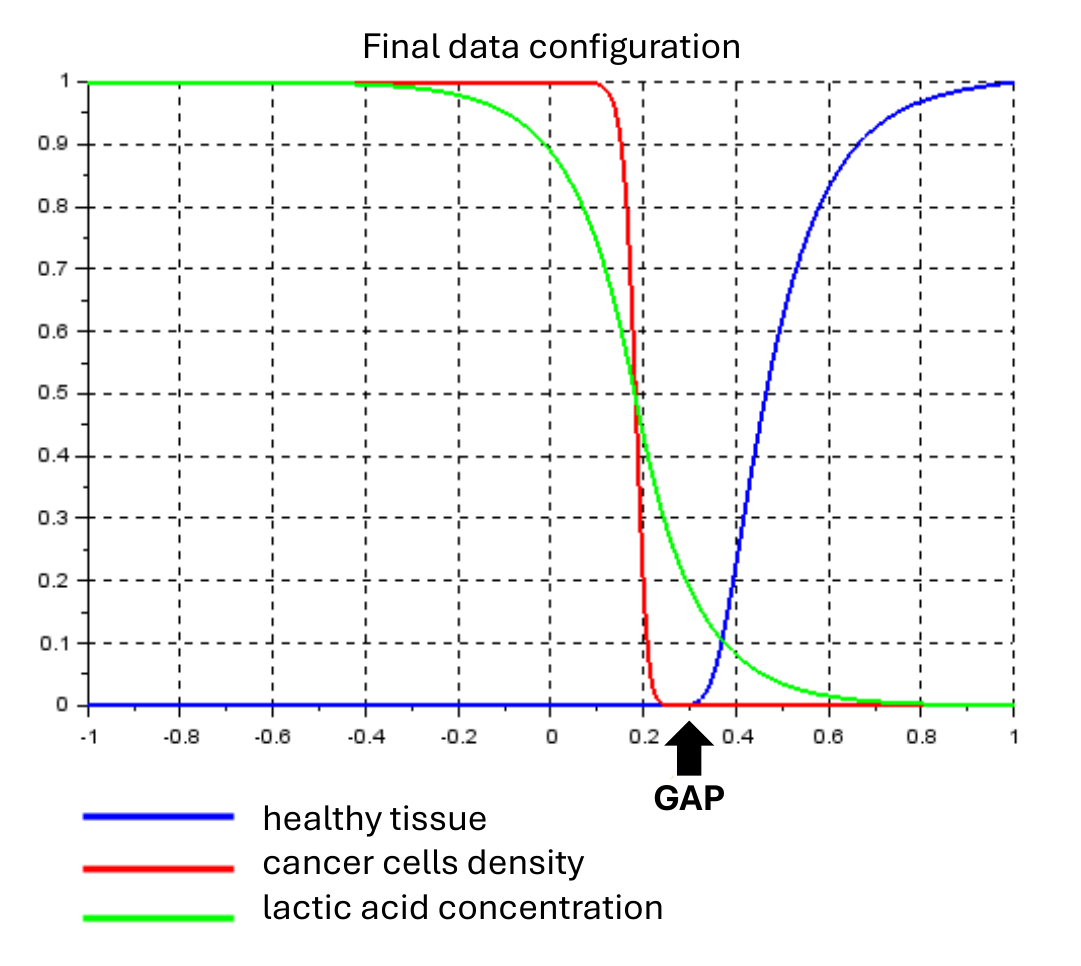}}
	\caption{\footnotesize{Different configurations of the numerical solution: comparison between heterogeneous
	evolution (b) and existence of the spatial interstitial gap within the homogeneous invasion (d).}}
	\label{finalconfig1}
\end{figure}

\begin{remark}
We point out that the choice of the spatial mesh indeed depends on the tissues structure and geometry 
since the space step $\Delta x$ should be small enough to capture the inhomogeneities of the lactic acid
diffusion function (refer to Section \ref{numerical algorithm}). 
This numerical issue becomes even more crucial in multidimensional simulations, for which unstructured
nonuniform meshes are usually adopted to account for different geometrical characteristics
without penalising the computational time.
On the other hand, the choice of the time step $\Delta t$ depends not only on biological and physical
considerations, but also on the numerical stability which is established by the CFL condition \cite{Quarteroni}.
\end{remark}

For the choice of the initial profiles, a piecewise linear decreasing density is taken into account for the cancer 
cells extending out from its core, where $v=1$, and getting towards zero; for the healthy cells density, 
the starting profile is simply obtained through a reflection, by imposing a complementary behaviour with 
respect to the cancer cells density; finally, the extracellular lactic acid concentration is initially equal to zero. 

The results reported in Fig.\ref{finalconfig1}(b-c-d)  {exhibit three different types of behaviours tuned}
by the parameter $d$  {which measures the destructive efficacy of the microenvironmental}
acidity on the healthy tissue,  {to be regarded} as an indicator of the tumour aggressiveness. 
From a qualitative point of view, all solutions evolve as forward propagation fronts moving from left 
to right with positive wave speed. 
The plot displayed in Fig.\ref{finalconfig1}(b) corresponds to a phenomenological regime known 
as \emph{heterogeneous invasion}, which turns out to happen when $d<1$. 
It is characterized by the coexistence of tumour and healthy tissue behind the wavefront, 
because a group of normal cells survives to the low aggressiveness of the tumour. 
On the other hand, when $d\gg1$, a different evolution shape takes place, the so-called
\emph{homogeneous invasion} shown in Fig.\ref{finalconfig1}(d), which is the most aggressive configuration.
Indeed, the healthy tissue is being completely destroyed behind the advancing tumour cells
wavefront because of the intensity of acidity induced into the environment. 
A narrow overlapping zone actually persists for increasing values of $d>1$, which produces 
 {\emph{(heterogeneous) hybrid configurations},} as shown in Fig.\ref{finalconfig1}(c).

From a mathematical point of view, the strong dissimilarity in terms of  {the} steepness of the wave profiles 
for $u$ and $v$ observed in Fig.\ref{finalconfig1} is justified by the fact that somehow $u$ inherits the parabolic 
regularity of the lactic acid concentration $w$ through the reaction term (see the first equation of the {\sf GG}
system \eqref{sistema GGbis}, whereas the diffusion constant $D$ of the neoplastic tissue is typically very small 
(see Table \ref{tavola1}). 
In fact, when passing from the  {Gatenby--Gawlinky model \eqref{sistema GG}} 
to its non-dimensionalized version \eqref{sistema GGbis}, that parameter is deduced as $D=D_2/D_3$
and it is physically relevant to assume the value of $D_3$ larger than $D_2$.


\subsection{Piecewise constant  {acid diffusion} with single jump}\label{subsect:single_jump}

 {Typical biological tissues are heterogeneous at macroscopic, mesoscopic and microscopic levels.
Intrinsic heterogeneity plays a crucial role also at intra-tumoral scale, \cite{MaruAlmePoly12}.
Here, we focus on the effect of space dependency on diffusivity of the lactic acid due to local fluctuations 
of the micro-environment.}

In this subsection, we provide numerical simulations for the {\sf GG} system \eqref{sistema GGbis} 
by applying  the scheme \eqref{eq discr 7} in presence of a heterogeneous diffusion coefficient $A$. 
We simplify the framework by taking a piecewise constant function.
There are biological and physical reasons behind this particular choice: indeed, we suppose that $A$ 
represents the diffusion of the lactic passing from a specific tissue or organ to another  {adjacent one.
As a consequence, the diffusivity experiences a sudden jump corresponding to} 
the different nature and geometry of the  {micro-environment.}
We aim  {to investigate} the effects of this transition on the motility of lactic acid, the consequences for
the consumption of healthy tissue and, finally, the modification occurring to the tumour propagation fronts.

We assume the parameters $d$, $r$, $D$ and $c$ as described above  {(see Table \ref{tavola1}),} 
$T$ as the final time and $L$ the length of the spatial interval $[0,L]$, with the spatio-temporal mesh built 
by fixing $\Delta x=0.005$ and $\Delta t=0.01$.  
As already discussed, the value of $\Delta x$ depends on the particular organ or tissue considered, 
whilst $\Delta t$ has to be chosen according to the CFL condition \cite{Quarteroni}.

\begin{remark} 
For the parameters of {\sf GG} system \eqref{sistema GGbis}, we assume $0 < D \ll 1$, 
because its form derives from the ratio  {$D:={D_2}/{D_3}$} between the diffusivity of malignant cells
and that of the $H^+$ ions, therefore it is reasonable to take $D_3$ much larger than $D_2$. 
Finally,  {the parameter} $d$ measures the destructiveness of pH against the healthy tissue 
 {to be interpreted as an indirect action of the tumour aggressiveness by means of the underlying
acidity micro-environment.}
\end{remark}

We deal with different initial profiles, either piecewise constant  {(i.e., with a jump)} or piecewise linear 
 {(i.e., continuous),} and different values of $A$, so that we distinguish two main cases with relative sub-cases. 
The first type of initial data is composed  {of} a decreasing piecewise constant density for the cancer cells
extending out from its core, where $v=1$, and having a jump at $x=\frac{L}{4}$; for the healthy cells density,
the starting graph is simply obtained through a reflection, by imposing a complementary behaviour with 
respect to $v$.
Finally, the extracellular lactic acid concentration is initially equal to zero. 
The corresponding graph  {is} shown in Fig.\ref{initialdata3}(b), referring to the spatial interval  $[0,1]$
whilst the second type is similar to that used in the previous Subsection as shown in Fig.\ref{initialdata3}(a).
\begin{figure}[htb]\centering
	\subfloat[][piecewise constant]{\includegraphics[width=.35\textwidth]{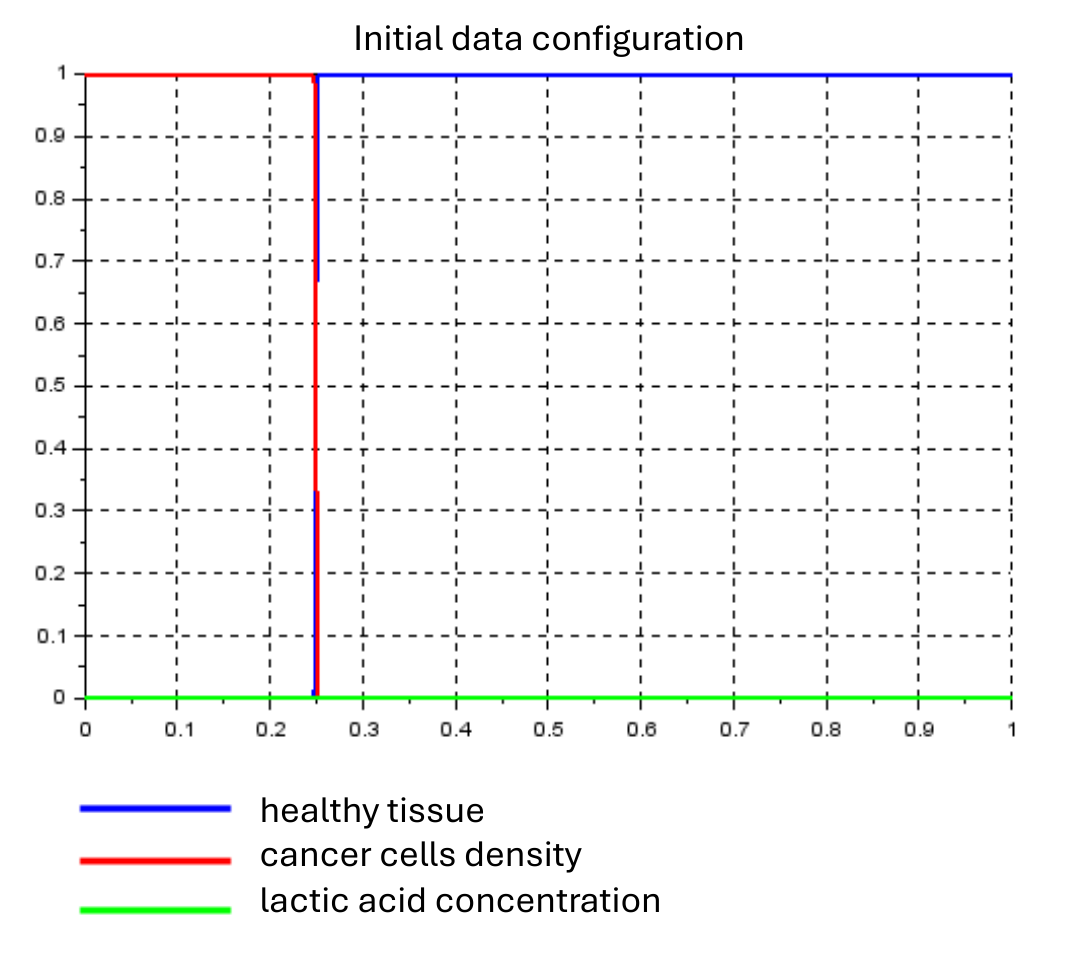}}\qquad
	\subfloat[][piecewise linear]{\includegraphics[width=.35\textwidth]{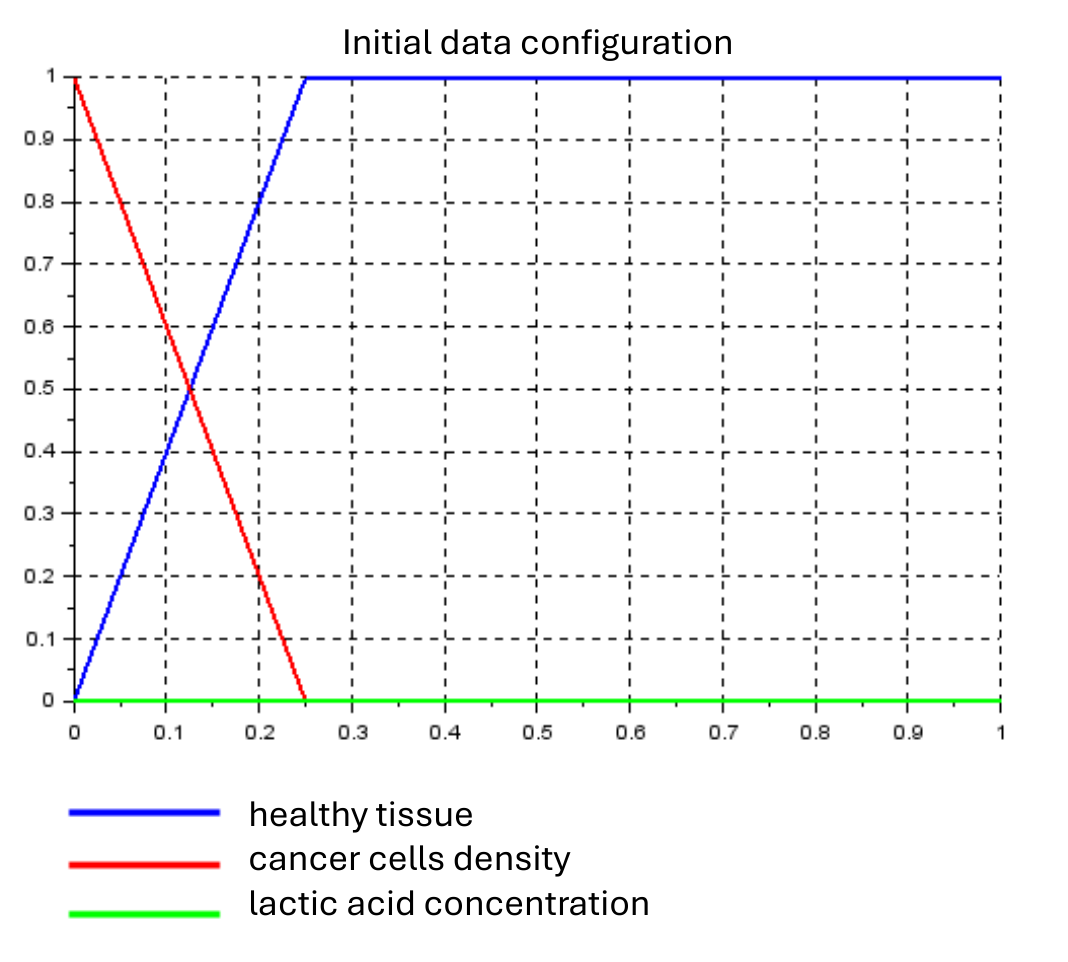}}
	\caption{\footnotesize{Initial profiles for the numerical simulations.}}
	\label{initialdata3}
\end{figure}
 
\begin{remark}
There is a specific experimental purpose for the choice of two different initial profiles. 
Indeed, the Riemann-type initial data in Fig.\ref{initialdata3}(a)  {is representative of}
many \emph{in vitro} experiments where a group of healthy cells is removed to inoculate a small colony
of cancer ones, hence the initial profile has a strong discontinuity which depends on sudden modifications
of the  {local} environment. 
On the other hand, the initial data in Fig.\ref{initialdata3}(b) refers to \emph{in vivo} situations, where the 
earlier development of cancer cells within a healthy tissue happens gradually  {(i.e., linearly, in the 
present case)} without immediately destroying the host environment.
\end{remark}

\begin{remark}
 {Mathematically speaking, both initial data represented in Fig.\ref{initialdata3} define a solution,
at least for short times and for large values of $x$, with a component $v$ that take values in a region
where the diffusion degenerates because $u$ is initially equal to 1.
Nevertheless, this is only a transient phenomenon since, for $t>0$, $u$ is strictly less than $1$.
Indeed, tumor tissue $v$ is strictly positive at values close to the left side of the interval, i.e., $x=0$,
generating acid concentration $w$. 
On its turn, $w$ diffuses with infinite speed of propagation, determining $u$ to be strictly decreasing
starting from $1$ because of the form of the time derivative of $u$.}
\end{remark}

 {To start with, we consider the following diffusion function}
\begin{equation}\label{single-jump}
	A(x):=\left\{\begin{aligned}
		a_1 &\quad& \text{if $x \in \bigl[0,\tfrac{5}{8}L\bigr]$\,,} \\
		a_2 &\quad& \text{if $x \in \bigl[\tfrac{5}{8}L,L\bigr]$\,.}
	\end{aligned}\right.
\end{equation}
We initially refer to the experimental data reported in Table \ref{tavola1}  {except for} 
$A$  which is not anymore homogeneous.
 {Moreover,} we select a piecewise constant initial profile as shown in Fig.\ref{initialdata3}(a).
The results obtained with $A$ strictly increasing $a_1<a_2$ are displayed in Fig.\ref{finalconfig4}, 
specifically for $a_1=0.1$ and $a_2=1$. 

\begin{figure}[htb]\centering
	\subfloat[][heterogeneous invasion, $d=0.5$]{\includegraphics[width=.4\textwidth]{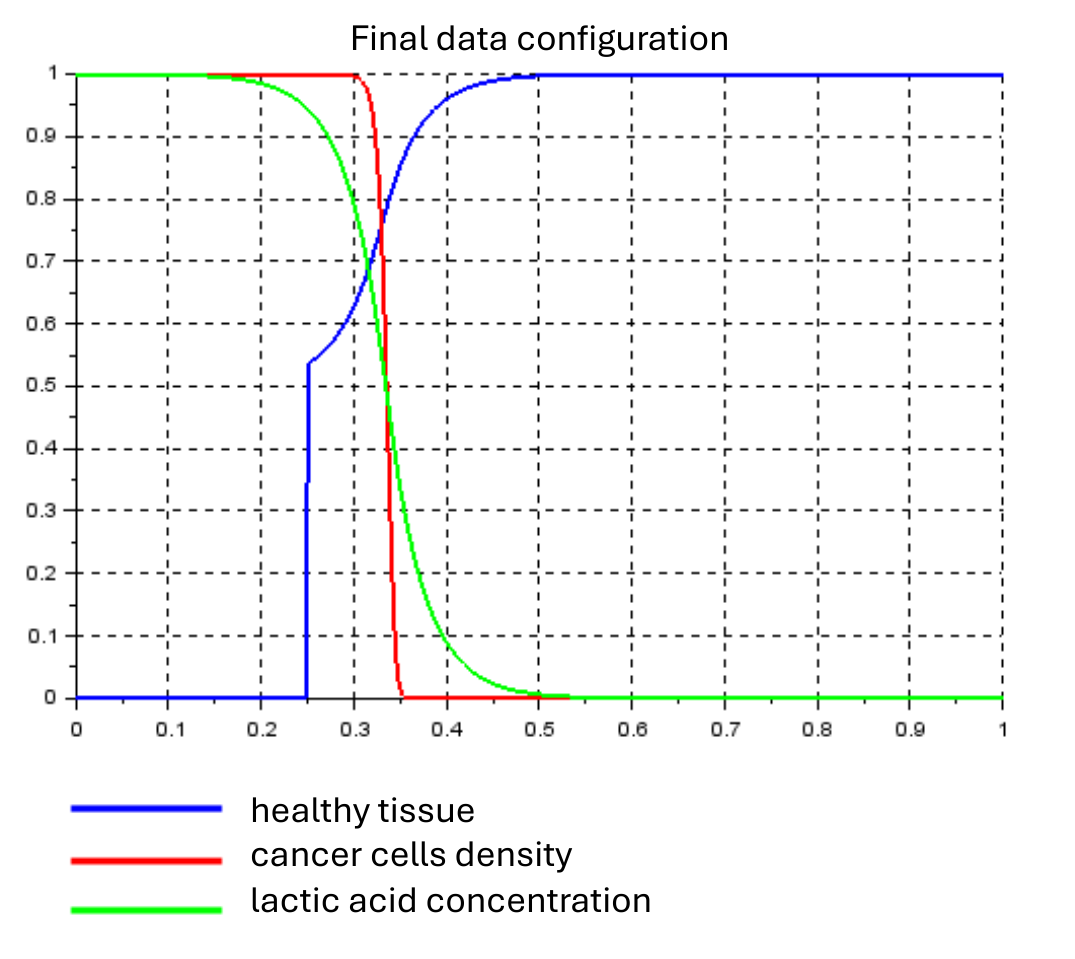}}\qquad
	\subfloat[][hybrid configuration, $d=1.5$]{\includegraphics[width=.4\textwidth]{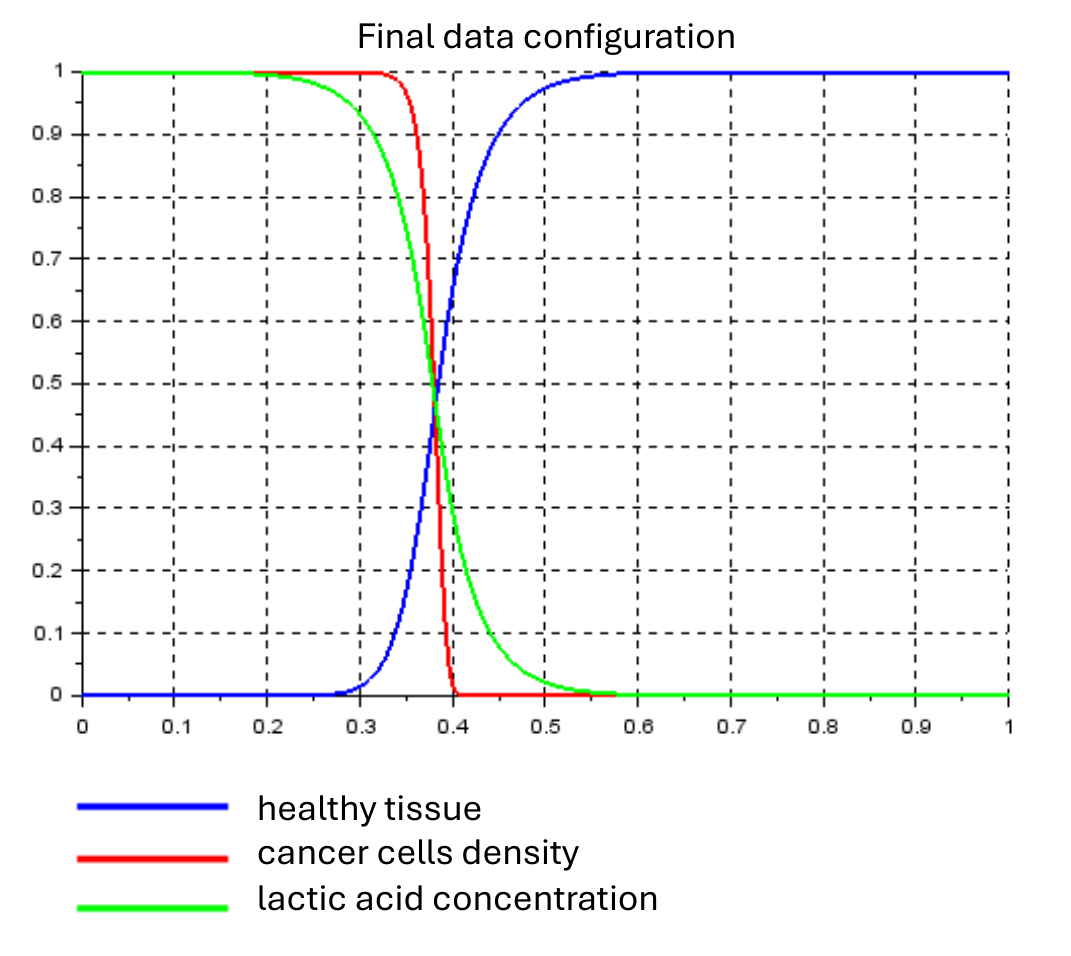}} \\
	\subfloat[][hybrid configuration, $d=12.5$]{\includegraphics[width=.4\textwidth]{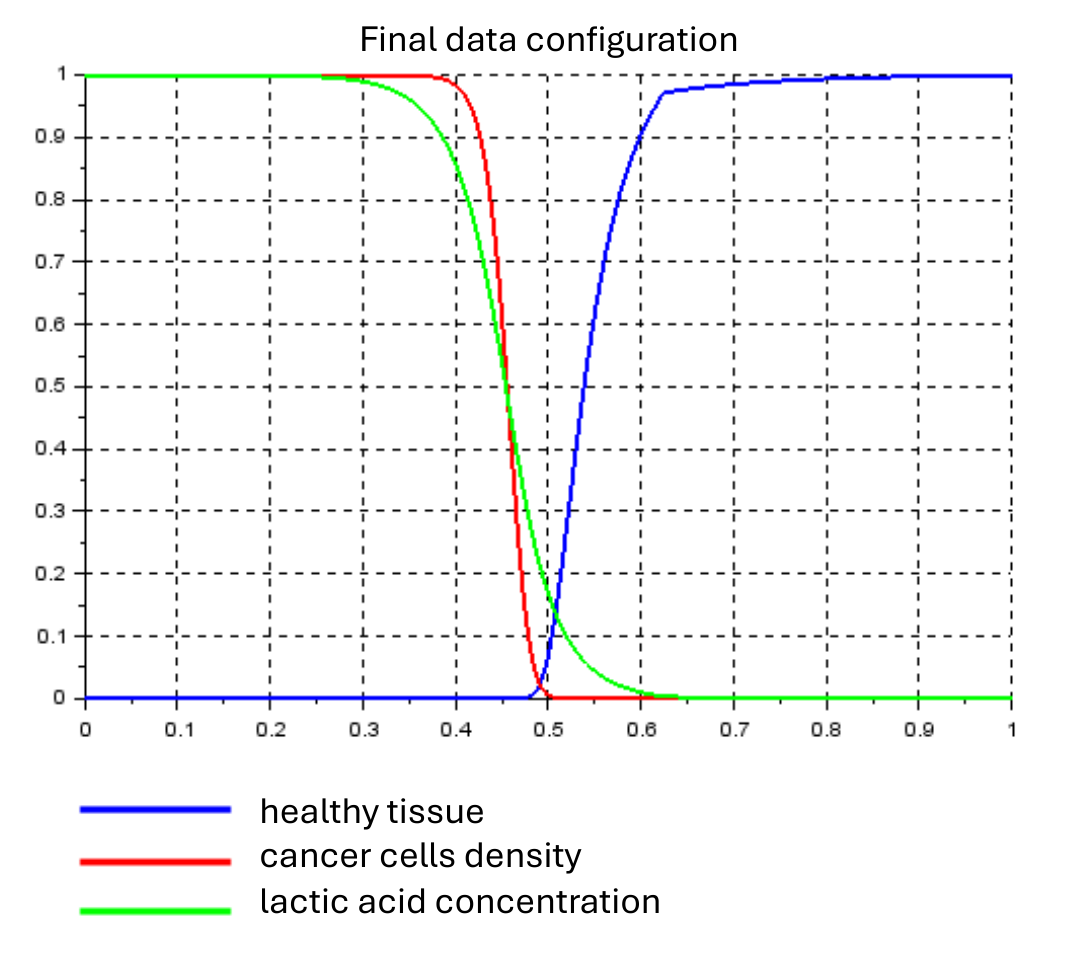}}\qquad
	\subfloat[][homogeneous invasion, $d=35$]{\includegraphics[width=.4\textwidth]{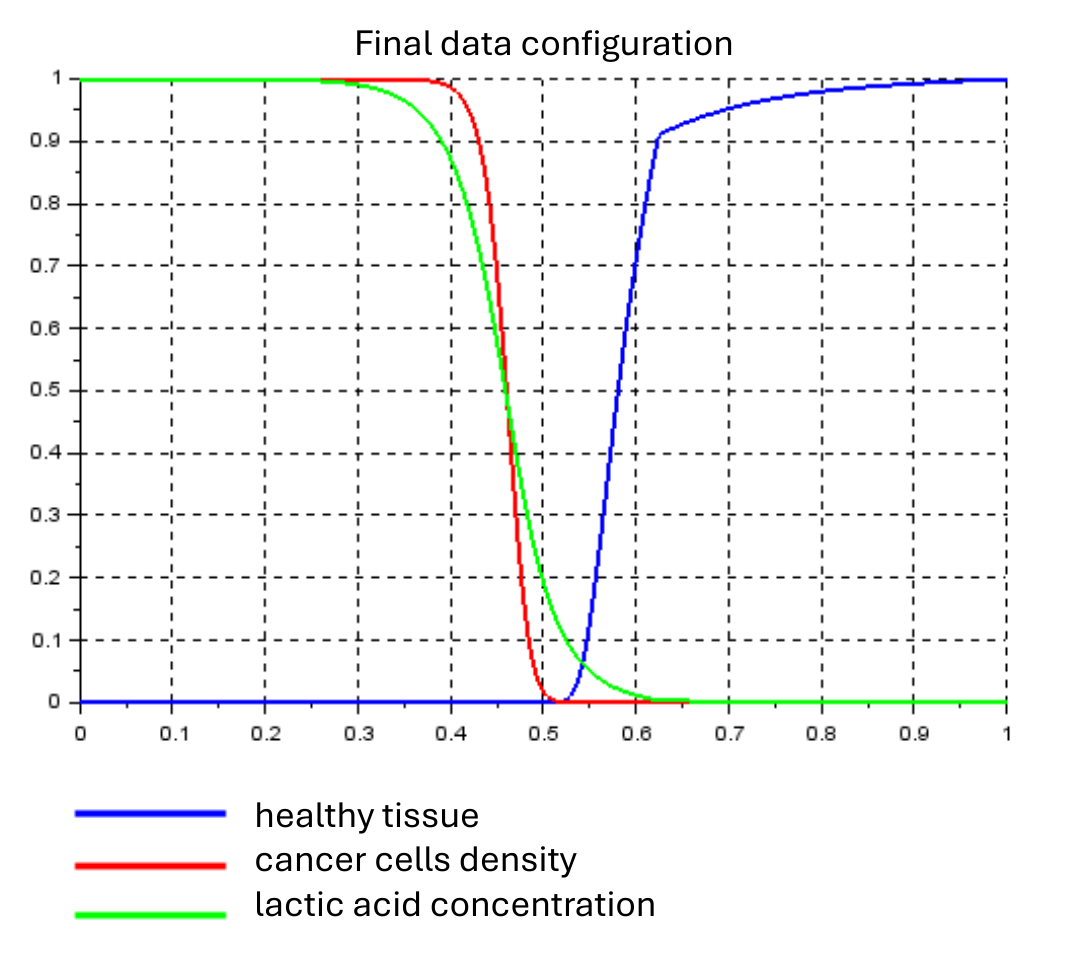}}
	\caption{\footnotesize{Different configurations of the numerical solution in presence of heterogeneous
	diffusion ($a_1<a_2$): comparison between heterogeneous evolution (a) and existence of the interstitial 
	gap within the homogeneous invasion (d).}}
	\label{finalconfig4}
\end{figure}

We notice the formation of the spatial interstitial gap when $d\gg 1$, and in particular for a larger coefficient, 
 {here} $d=35$, because of the slowdown effect induced by the smaller values of the diffusion rates
for the lactic acid.
If we decrease the value of $|a_1-a_2|$, this modification of the diffusion inhomogeneity significantly
affects the formation of the gap, which occurs for smaller values of $d$.

\begin{remark} 
At this stage, we point out another advantage of using the finite volume method, together
with those already discussed in Subsection \ref{finitevolumes}. 
Indeed, when considering a heterogeneous function $A$ with discontinuities, the mathematical problem
associated with {\sf GG} system \eqref{sistema GGbis} loses its regularity, since the characteristics of 
parabolicity are weaker, and this phenomenon is clear from Fig.\ref{finalconfig4}. 
From an analytical point of view, this means that the functional $A_x w_x$ may not be well defined 
at $x=\frac{5}{8}L$, which is the point of irregularity. 
The integral formulation of the finite volume approach allows to correctly deal with this issue.
\end{remark}

The results obtained with $A$ strictly decreasing ($a_1>a_2$) are displayed in Fig.\ref{finalconfig6}, 
specifically for $a_1=1$ and $a_2=0.1$. 
We notice again the formation of the spatial interstitial gap when $d\gg 1$, but in this case for $d=12.5$ 
the gap is largely present.

\begin{figure}[]\centering
	\subfloat[][heterogeneous invasion, $d=0.5$]{\includegraphics[width=.35\textwidth]{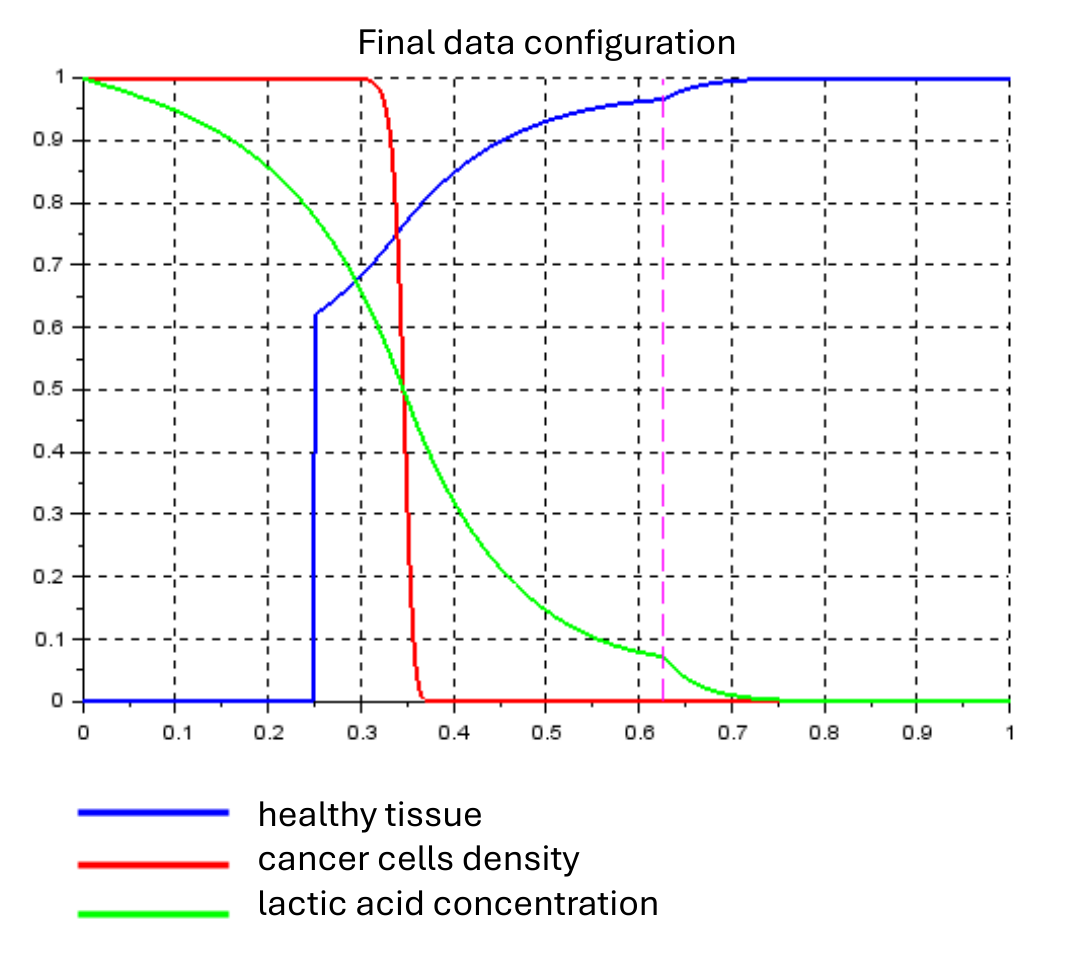}}\qquad
	\subfloat[][homogeneous invasion, $d=12.5$]{\includegraphics[width=.35\textwidth]{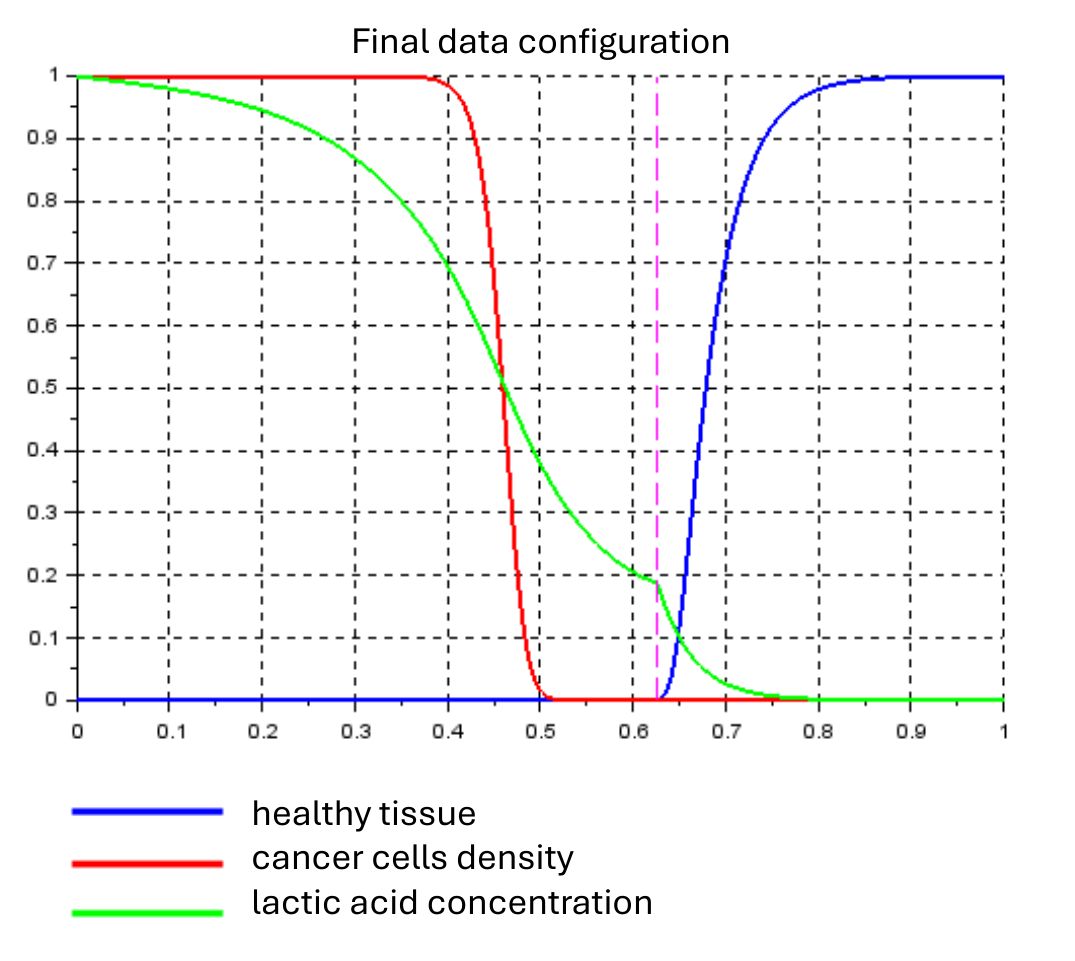}}
	\caption{\footnotesize{Configurations of the numerical solution in presence of heterogeneous diffusion ($a_1>a_2$): 
	comparison between heterogeneous evolution for $d$ small (a) and existence of the spatial interstitial 
	gap within the homogeneous invasion for $d$ large (b).}}
	\label{finalconfig6}
\end{figure}

It is also interesting to remark that the final profile on the healthy tissue density for hybrid and heterogeneous
regimes clearly manifests a discontinuity in its derivative, which is inherited from the lactic acid configuration 
through the nonlinear reaction term in the first equation of {\sf GG} system \eqref{sistema GGbis}. 
The magenta vertical (dashed) line on the graphs of Fig.\ref{finalconfig6} indicates the abscissa of discontinuity
for the diffusion function $A$, which passes exactly through the point where the lactic acid and the healthy tissue
concentrations have a discontinuity in their derivative,
 {hence the solution in some of its components exhibits some corner points.}

Moreover, if we decrease the value of $|a_1-a_2|$, this does not affect the formation and the width of the gap, 
but the final profile of the lactic acid concentration is smoother and the corner point is not anymore visible 
(see Fig.\ref{finalconfig7}).

\begin{figure}[htb]\centering
	\subfloat[][diffusion rates $a_1=1$ and $a_2=0.8$]{\includegraphics[width=.35\textwidth]{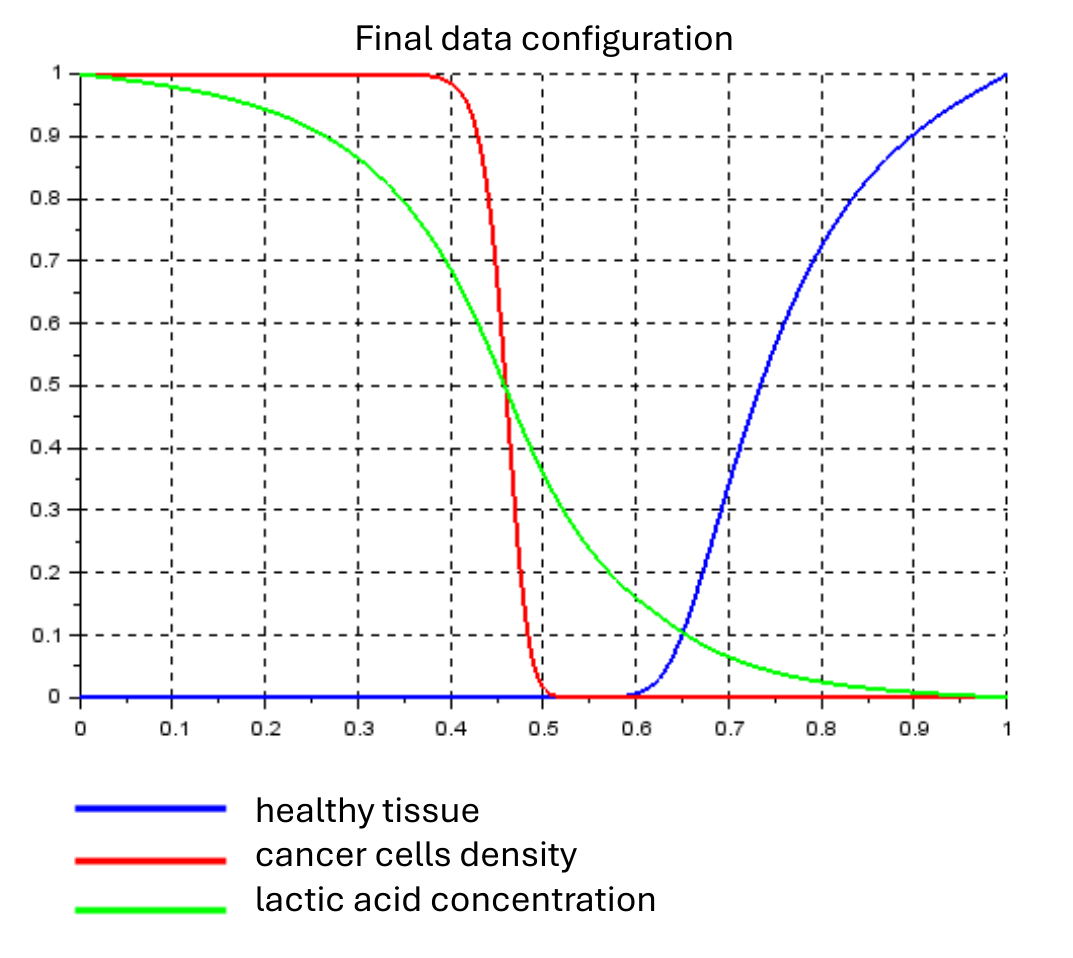}}\qquad
	\subfloat[][diffusion rates $a_1=0.3$ and $a_2=0.1$]{\includegraphics[width=.35\textwidth]{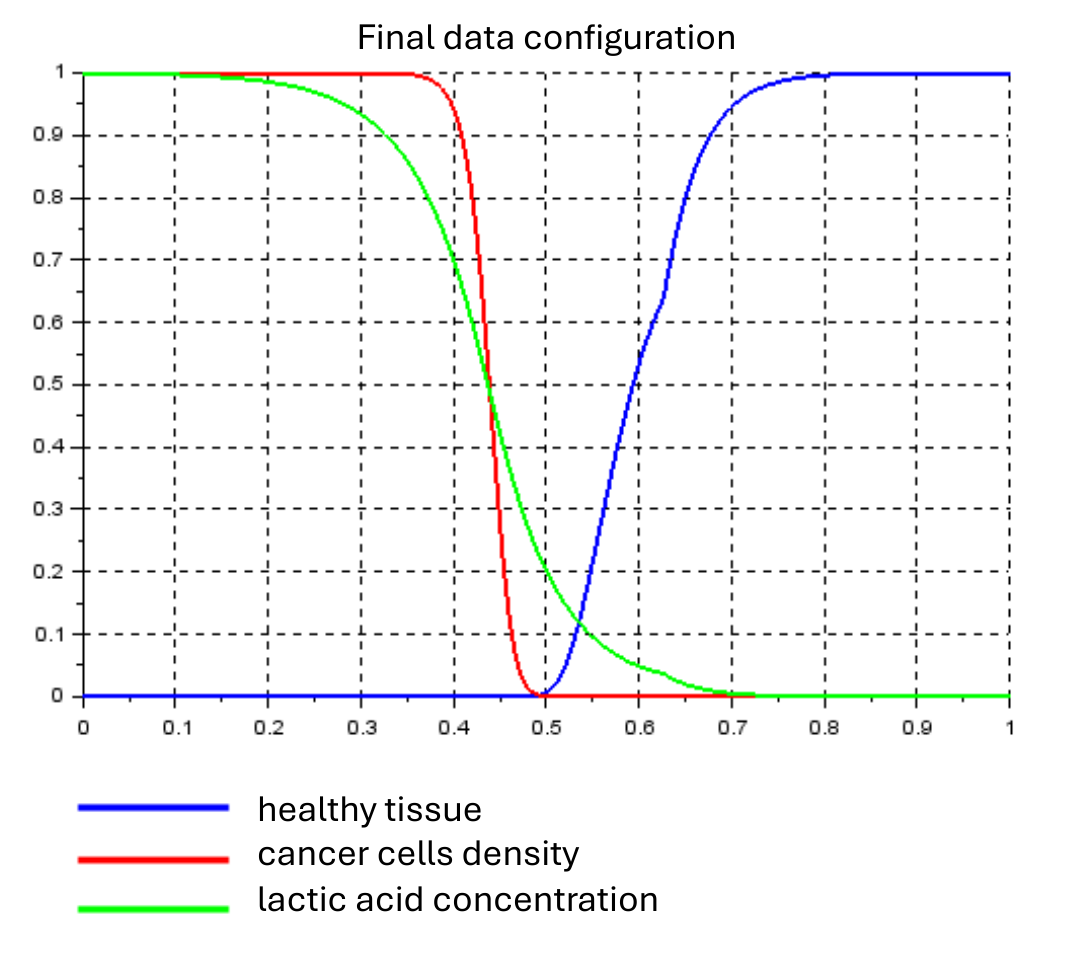}} 
	\caption{\footnotesize{Homogeneous invasion from  initial data in Fig.\ref{initialdata3}(b)
	in presence of heterogeneous piecewise constant diffusion function $A$ 
	with same value of the jump $|a_1-a_2|$ and $d=12.5$.}}
	\label{finalconfig7}
\end{figure}

The presence of corner points in the solution profiles for $u$ and $w$, which is observed 
when $A$ becomes discontinuous, obviously disappears in the case of homogeneous diffusion 
(see Fig.\ref{finalconfig1}).
Furthermore, the profile of solution $u$ becomes smoother as $d$ becomes larger enough for the
interstitial gap to emerge, as shown in Fig.\ref{finalconfig6}, because then the discontinuity point of 
$A$ falls into the region where the healthy tissue is null. 

Then, we change the initial profile and we use the one displayed in Fig.\ref{initialdata3}(b). 
Performing numerical simulations with the same data used above, we can see substantial
differences only for small values of $d$. 
\begin{figure}[bht]\centering
	\subfloat[][diffusion rates $a_1=0.1$ and $a_2=1$]{\includegraphics[width=.35\textwidth]{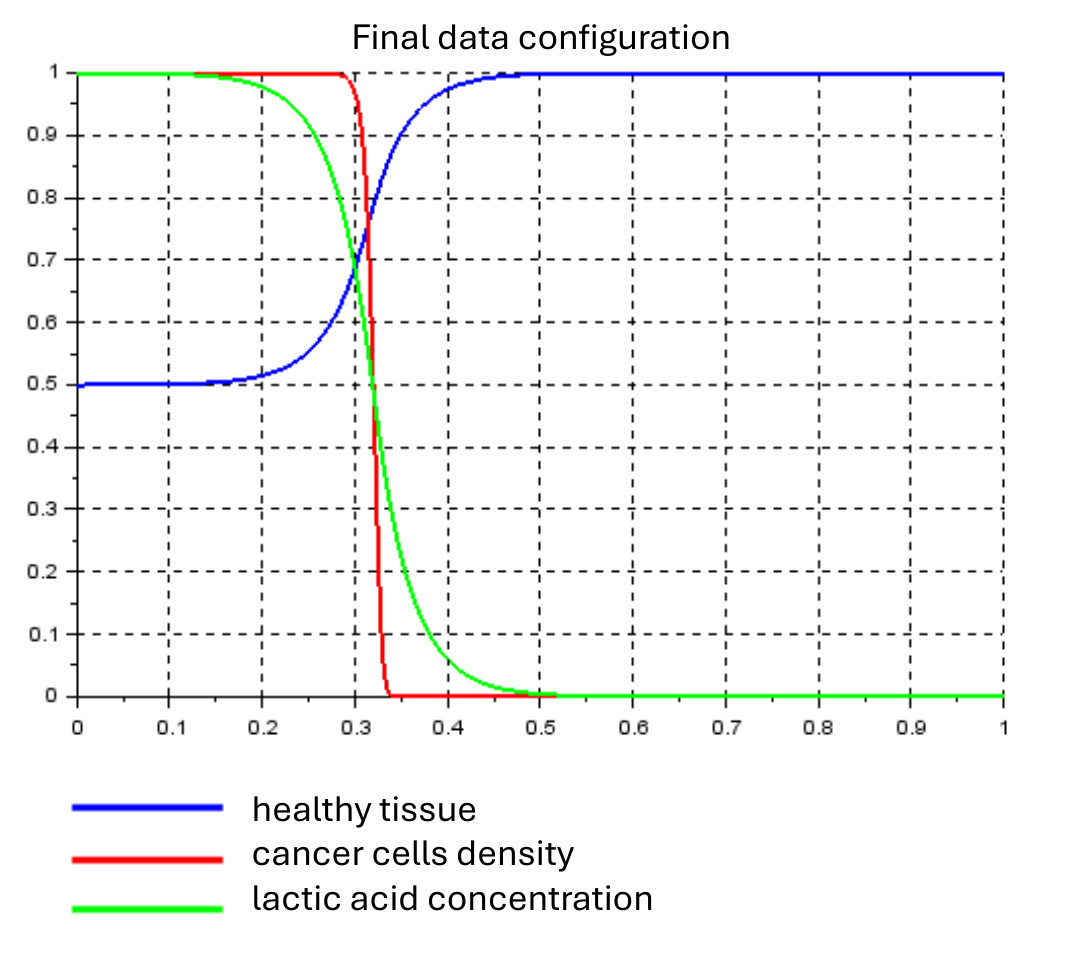}}\qquad
	\subfloat[][diffusion rates $a_1=1$ and $a_2=0.1$]{\includegraphics[width=.35\textwidth]{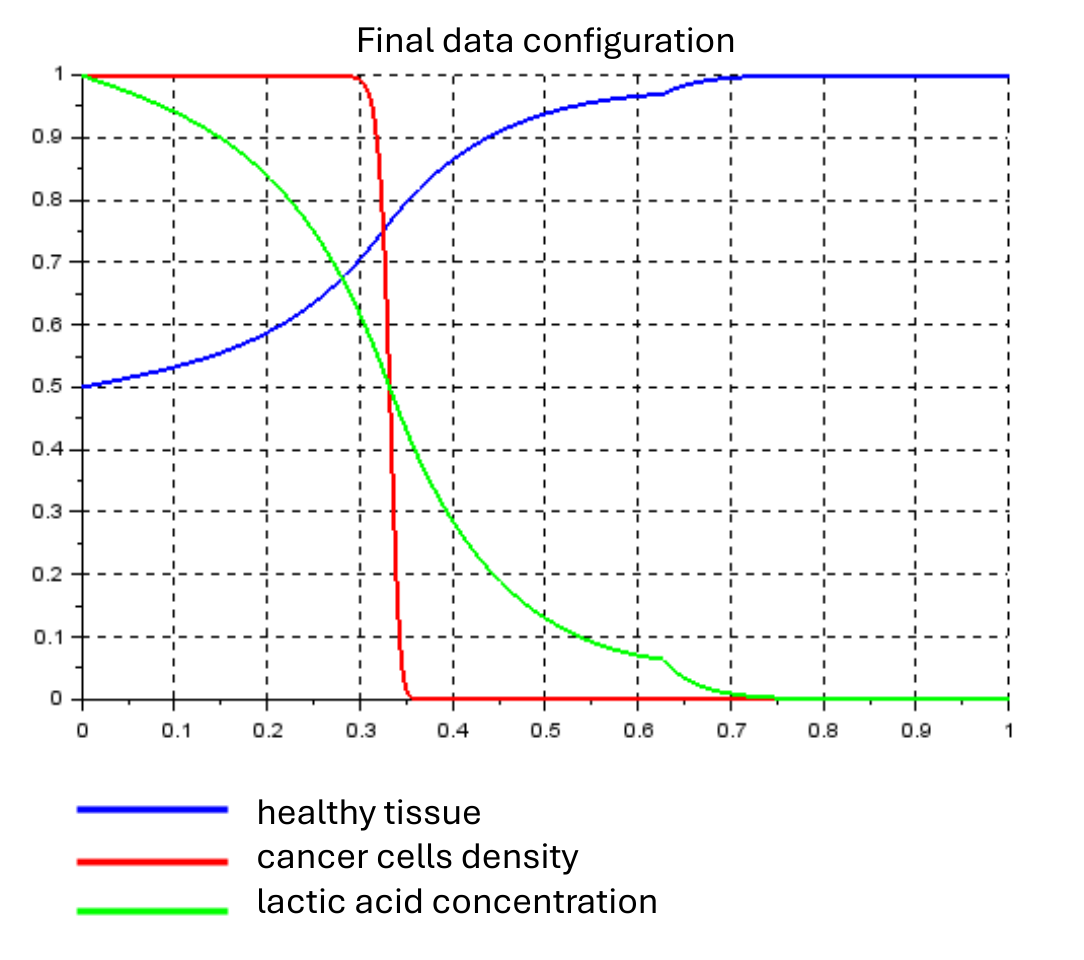}}
	\caption{\footnotesize{Heterogeneous invasion from initial data in Fig.\ref{initialdata3}(b)
	in presence of heterogeneous piecewise constant diffusion function $A$ for $d=0.5$\,.}}
	\label{finalconfig11}
\end{figure}
The results obtained for $d=0.5$ and the other parameters of Table \ref{tavola1} are reported in Fig.\ref{finalconfig11} 
with a comparison between the case of $A$ strictly increasing ($a_1<a_2$) and strictly decreasing ($a_1>a_2$).
These results suggest that if the tumour is not very aggressive (i.e. $d$ is small), then the healthy tissue is not 
completely destroyed.


\subsection{Periodic  {acid} diffusion function}\label{subsect:periodic}

In this section, we provide numerical simulations of {\sf GG} system \eqref{sistema GGbis} using 
the scheme \eqref{eq discr 7} in presence of a heterogeneous function $A$, which is supposed to be periodic. 
This choice is related to the structure of most parts of tissues, which are porous and heterogeneously 
permeable and hence responsible  {for} a non-uniform invasion by the lactic acid. 
Such a function $A$ may thus represent the effect of increased diffusivity for the $H^+$ ions due
to easy passage through tissues for half its period, and the effect of reduction due to obstacles 
along the way \cite{chaplain,ChapLola06}.

\begin{remark}
It is worthwhile to underline that the choice of $\Delta x$ has to be made carefully, in order to avoid trivial 
interpolations of the periodic diffusion function $A$ if the space-step is too close to (a multiple of) its period.
\end{remark}

We choose the initial profiles in Fig.\ref{initialdata3}(b) on the spatial interval $[0,1]$ 
and we refer to Table \ref{tavola1} for the numerical parameters.

\begin{figure}[tbh]\centering
	\subfloat[][\emph{heterogeneous invasion, $d=0.5$}]{\includegraphics[width=.35\textwidth]{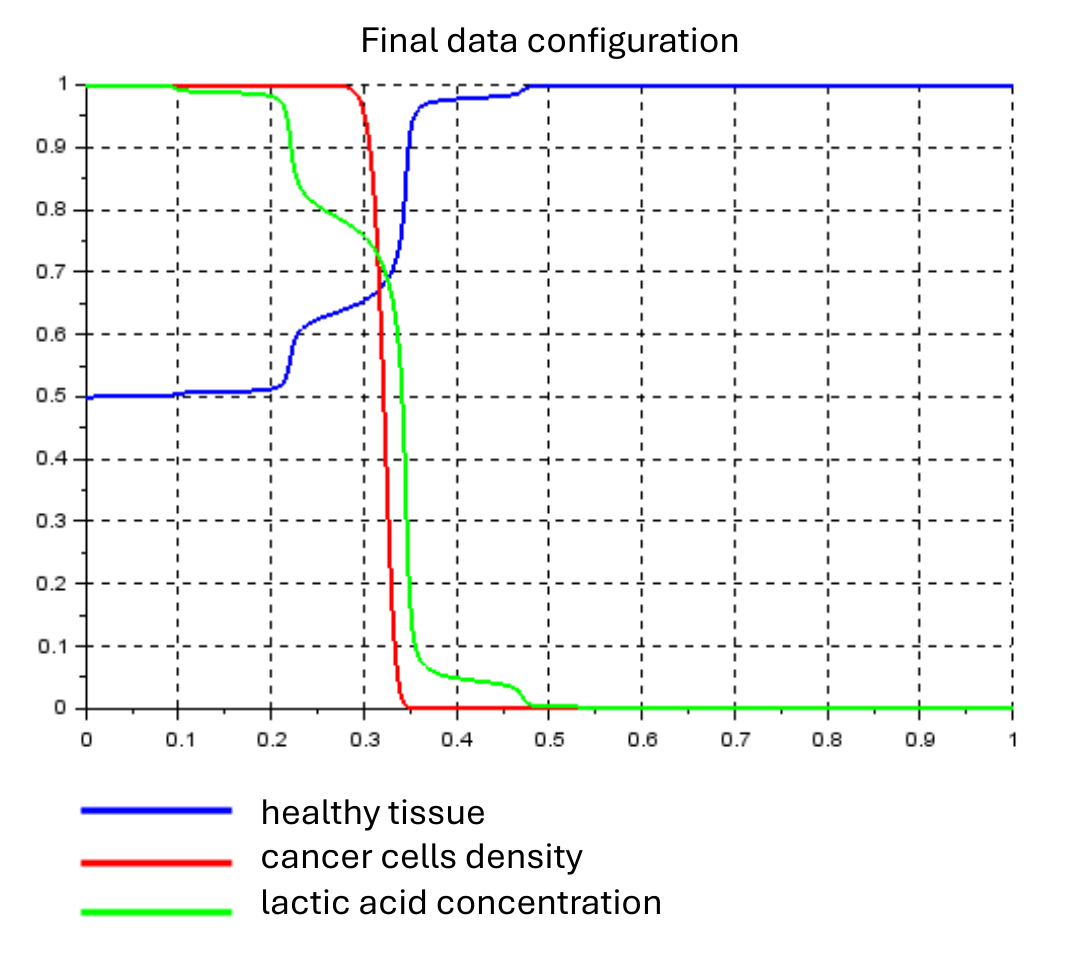}}\qquad
	\subfloat[][\emph{hybrid configuration, $d=1.5$}]{\includegraphics[width=.35\textwidth]{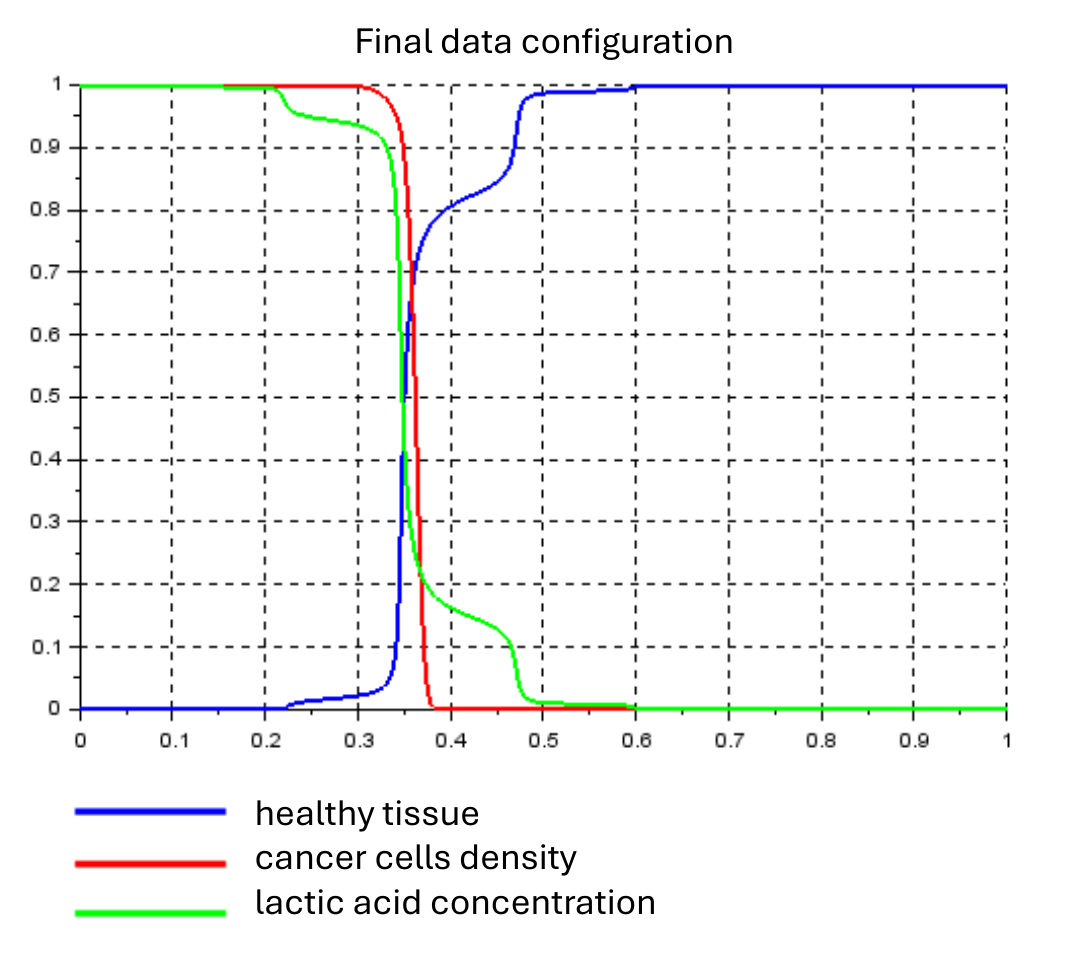}} \\
	\subfloat[][\emph{hybrid configuration, $d=30$}]{\includegraphics[width=.35\textwidth]{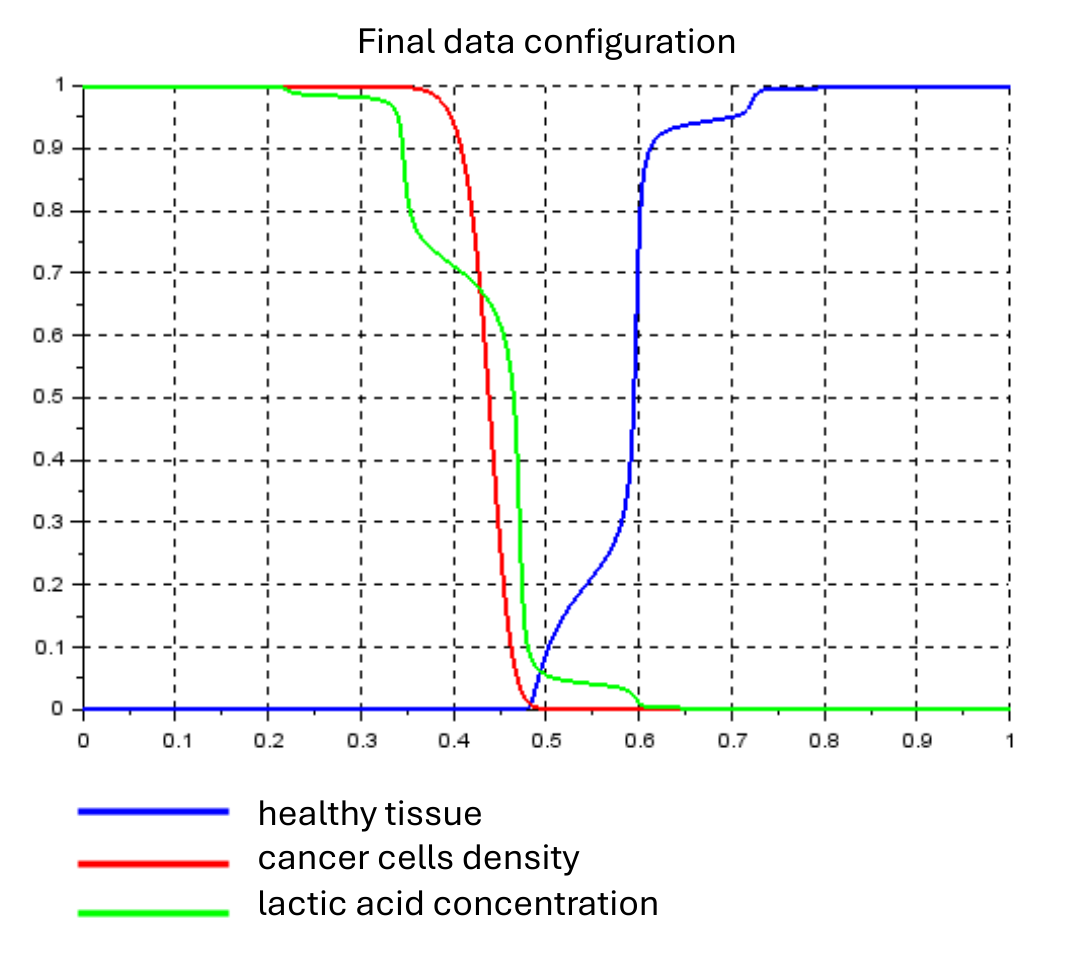}}\qquad
	\subfloat[][\emph{homogeneous invasion, $d=60$}]{\includegraphics[width=.35\textwidth]{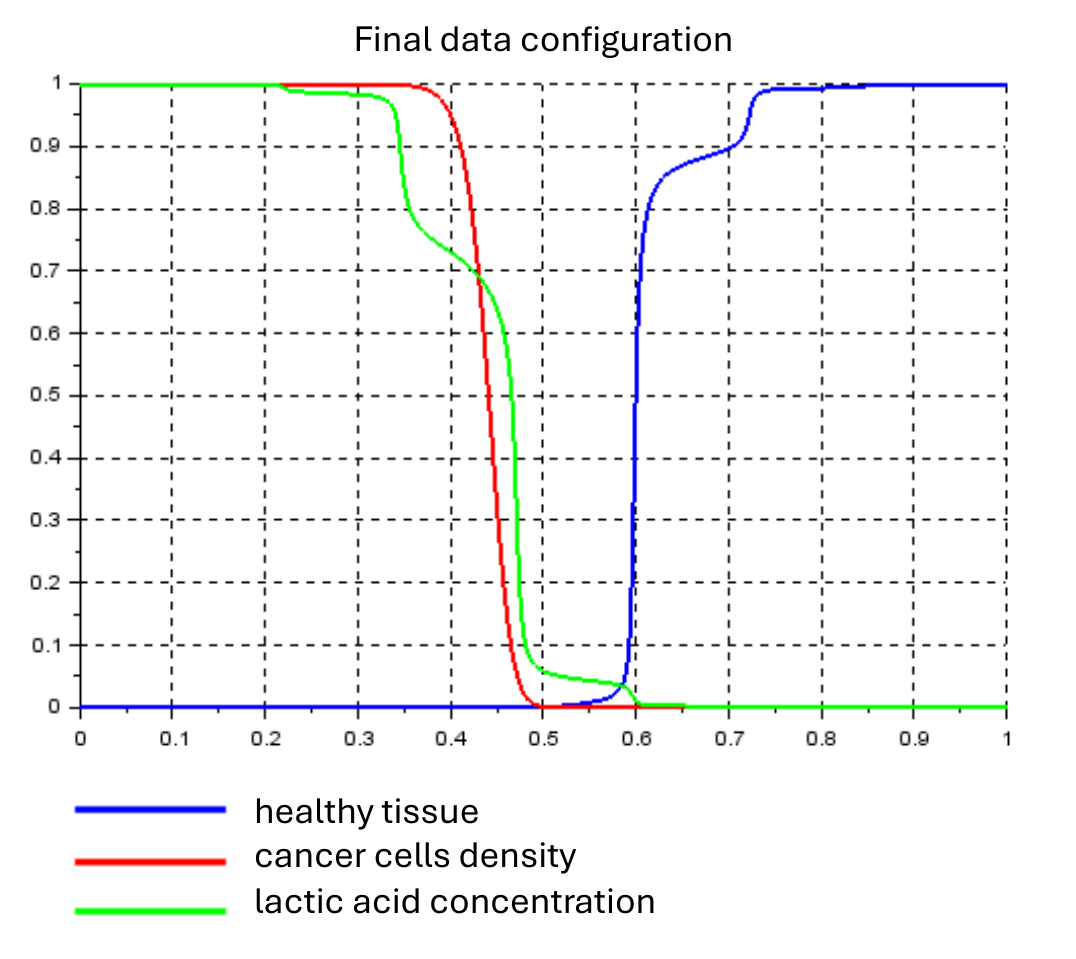}}
	\caption{\footnotesize{Different configurations of the numerical solution in presence of periodic diffusion $A$ with
	frequency $\omega=50$ and existence of interstitial gap within the homogeneous invasion (d).}}
	\label{perconfig3}
\end{figure}

We analyse different configurations for the diffusion function $A$ by changing the amplitude and frequency
of its oscillations, according to the following expression 
\begin{equation}\label{periodico1}
	A(x) = \tfrac12\bigl\{\alpha_0 \bigl[1-\sin(\omega x)\bigr]+ \alpha_1 \bigl[1+\sin(\omega x)\bigr]\bigr\}
	= \tfrac12 (\alpha_1+\alpha_0) + \tfrac12 (\alpha_1-\alpha_0) \sin(\omega x)\,,
\end{equation}
where $\alpha$ is the fixed amplitude of each oscillation and $\omega$ is the frequency 
(so that ${1}/{\omega}$ is the period).
In order to be physically consistent with the {\sf GG} model \eqref{sistema GGbis}, it is assumed that 
$0<\alpha_0\le A(x)\le \alpha_1$ for all $x$ in the domain, thus $\alpha_0$ defining the  {uniform positive}
constant of  {diffusivity function $A=A(x)$.}

A modification of the parameter $\omega$ in \eqref{periodico1} implies a change in the frequency of oscillations: 
we suppose  $\omega\gg1$ and we  initially choose $\alpha_0=0.1$ and $\alpha_1=1$. 
Numerical simulations with $\omega=50$ on the spatial interval $[0,1]$ are displayed in Fig.\ref{perconfig3}. 

We propose three relevant considerations. 
First of all, despite the perturbations observed in Fig.\ref{perconfig3} for the profiles $u$ and $w$ 
(with those for $u$ inherited from $w$ through the reaction term $-duw$), the solutions still exhibit a front-type
behaviour which is then very robust (such stability is actually preserved for smaller values of the frequency). 
Secondly, the value of $\omega$ affects the formation of the interstitial gap, which becomes visible for $d=60$.
Indeed, analogous simulations are reported in Fig.\ref{perconfig4} (see  {Appendix \ref{appendix}}), where we 
increase the frequency value to $\omega=100$ and the appearance of the interstitial gap is again for $d=50$.

Moreover, the perturbations are more intense for higher values of $\omega$, whilst they are obviously damped
for larger values of $d$ since the healthy tissue becomes null.
Finally, looking at the wave speed approximation given by the LeVeque-Yee formula \eqref{levequeyee} for different
frequencies of the diffusion function $A$ (see Fig.\ref{wavespeed3}), we can infer that the value of $\omega$
does not affect the propagation velocity (the solutions for $\omega=50$ and $\omega=100$ converge to the
same asymptotic value 
 {$\theta^*\approx 0.012$).}
\begin{figure}[bht]\centering
	\subfloat[][\emph{$\omega=50$}]{\includegraphics[width=.35\textwidth]{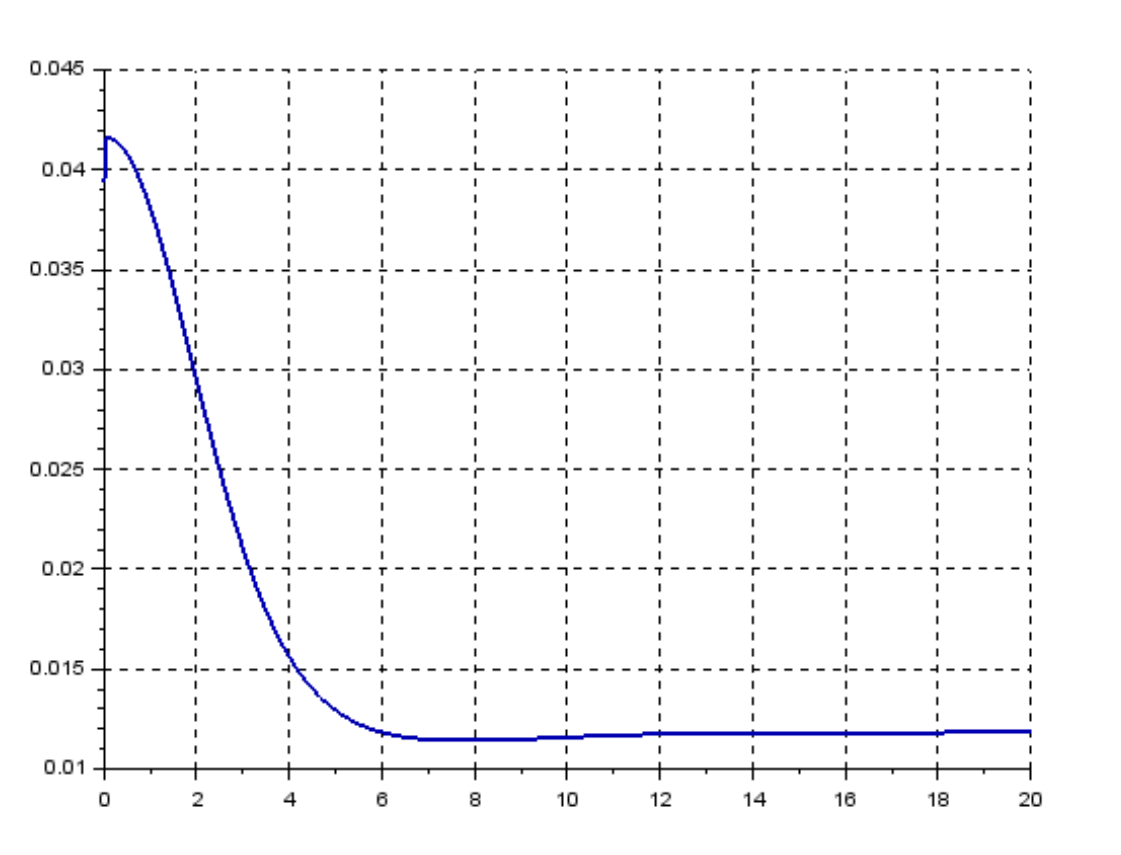}}\qquad
	\subfloat[][\emph{$\omega=100$}]{\includegraphics[width=.35\textwidth]{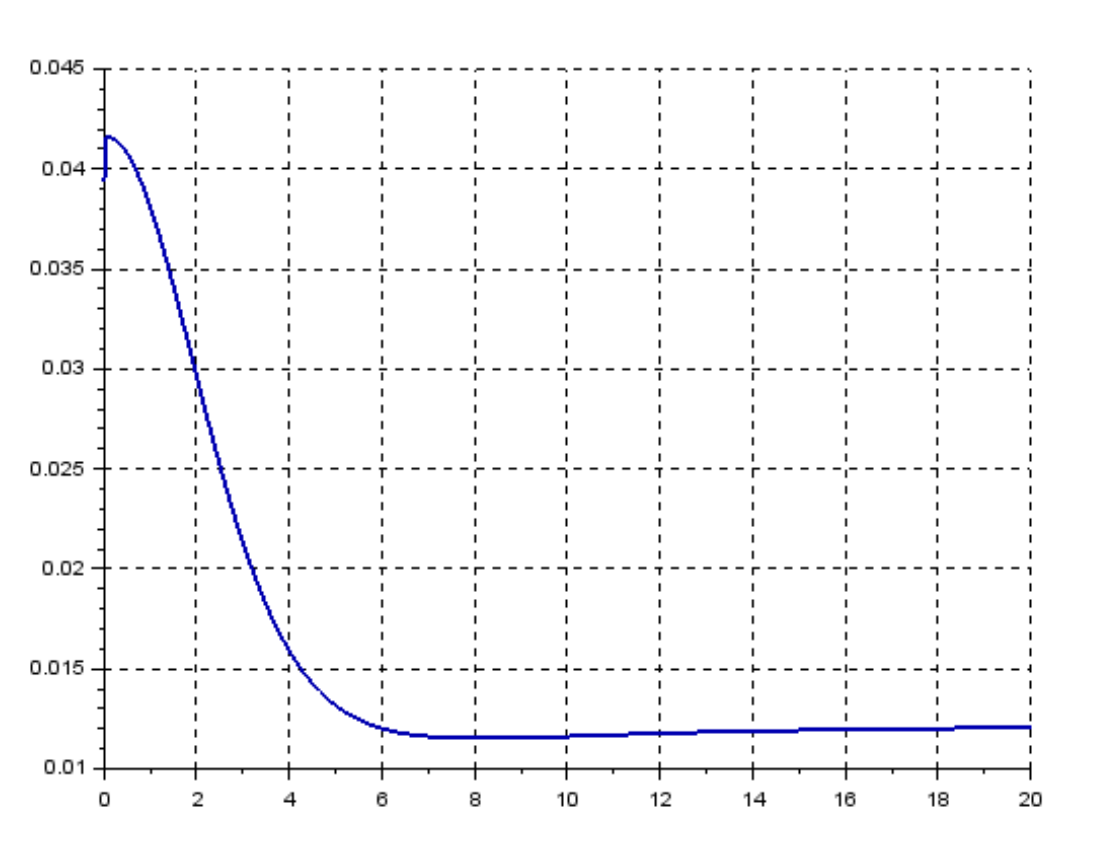}}
	\caption{\footnotesize{Wave speed approximation for different values of the frequency $\omega$ and $d=20$.}}
	\label{wavespeed3}
\end{figure}


 {\subsubsection{Modifying amplitude of the periodic acid diffusion}}\label{amplitude-intensity}

The amplitude and intensity of the diffusion function $A$ also  {modify} the behaviour of the numerical solution, 
especially when the value of $\omega$ is very large. 
For the following simulations, we use the formula \eqref{periodico1} with different values of $(\alpha_0,\alpha_1)$.

Firstly, we choose $(\alpha_0,\alpha_1)=(0.4,0.6)$, so that the amplitude is $|\alpha_1-\alpha_0|=0.2$
while the other parameters are those in Table \ref{tavola1} with $\omega=50$. 
\begin{figure}[bht]\centering
	\subfloat[][$d=0.5$]{\includegraphics[width=.4\textwidth]{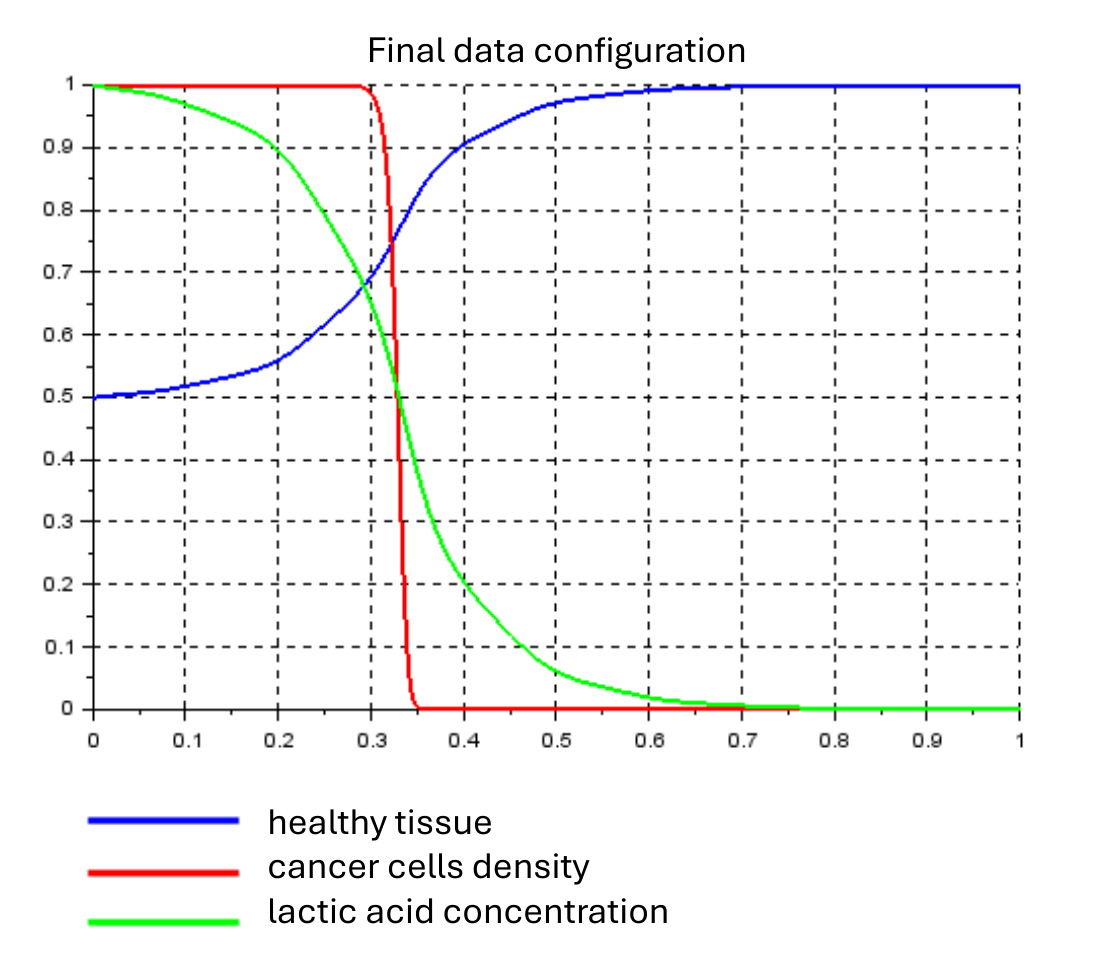}}\qquad
	\subfloat[][$d=20$]{\includegraphics[width=.4\textwidth]{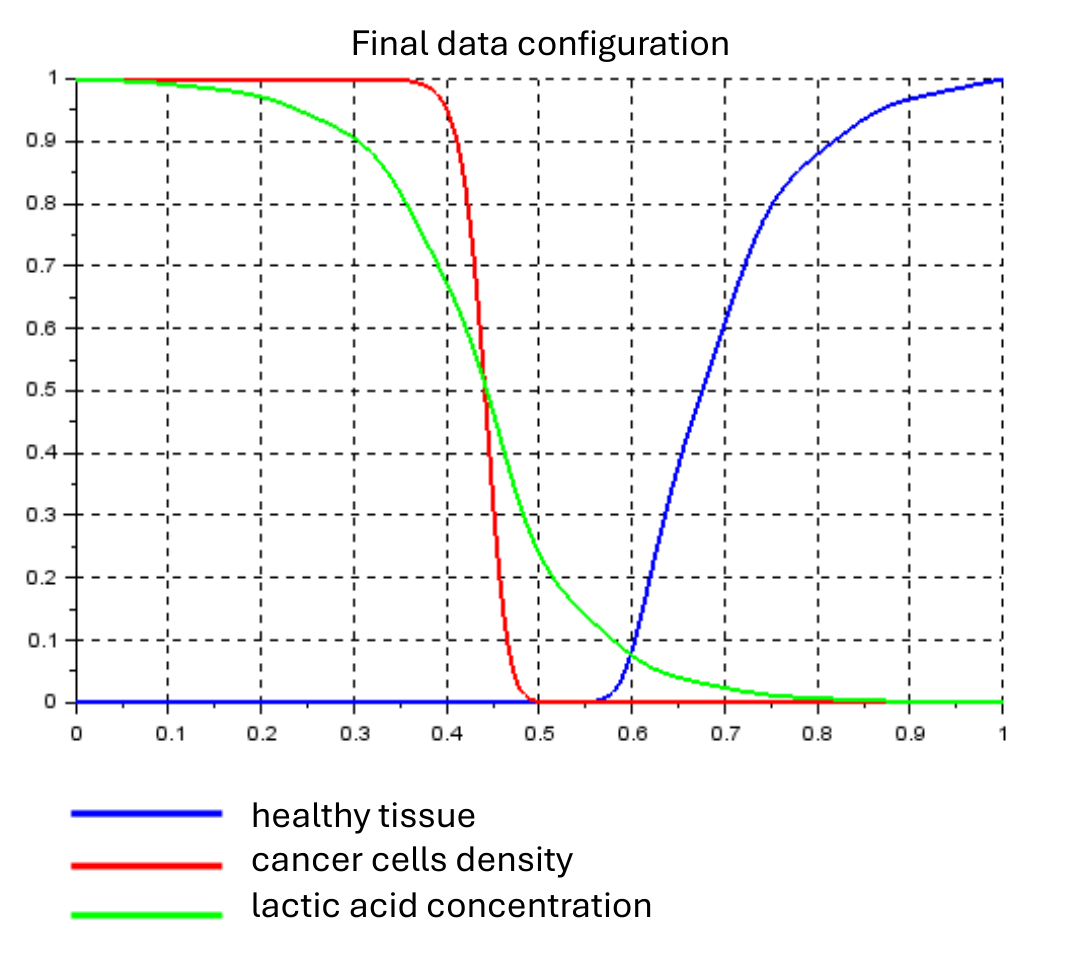}}
	\caption{\footnotesize{Heterogeneous (left) and homogeneous invasions (right)
	with $\omega=50$ and $(\alpha_0,\alpha_1)=(0.4,0.6)$.}}
	\label{perconfig7}
\end{figure}
The numerical solutions displayed in Fig.\ref{perconfig7} are smoother with respect to Fig.\ref{perconfig3}
and  {these} characteristics do not depend on the values of $\alpha_0$ and $\alpha_1$, 
but on their difference in absolute value, as shown in Fig.\ref{perconfig8} and Fig.\ref{perconfig11} where 
 {$(\alpha_0,\alpha_1)=(0.8,1)$ and $(\alpha_0,\alpha_1)=(0.1,0.3)$, respectively (see  {Appendix \ref{appendix}}).}
In Fig.\ref{perconfig12} and Fig.\ref{perconfig13}, we choose a very small intensity for $A$, 
which is $|\alpha_1-\alpha_0|=0.05$  and we compare with the case 
 {$(\alpha_0,\alpha_1)=(0.95,1)$ with $(\alpha_0,\alpha_1)=(0.01,0.06)$, respectively.} 
It is observable that the more the acid concentration is weak, the more the interstitial gap appears 
for large value of $d$ (see Fig.\ref{perconfig13}, where the gap is not yet present for $d=200$).

We stress the fact  {that the front profile of numerical solution is clearly distinguishable}
even if we increase the value of $\omega$, see \cite{ScanMascSime20}.\\


 {\subsection{Uncontrolled growth: the role of tumor reproduction}\label{uncontrolled}}

 {Next, we focus on the role played by the tumor reproduction parameter $r$.}
In order to investigate the variations of the wave speed propagation, we compare the approximation 
provided by the LeVeque--Yee formula \eqref{levequeyee} for the case in Fig.\ref{perconfig3} 
($\omega=50$ and $r=1$) with that in Fig.\ref{perconfig10}, which is constructed using the same 
parameters, frequency and amplitude but with a larger value $r=10$.  

\begin{figure}[bht]\centering
	\subfloat[][$d=0.5$]{\includegraphics[width=.35\textwidth]{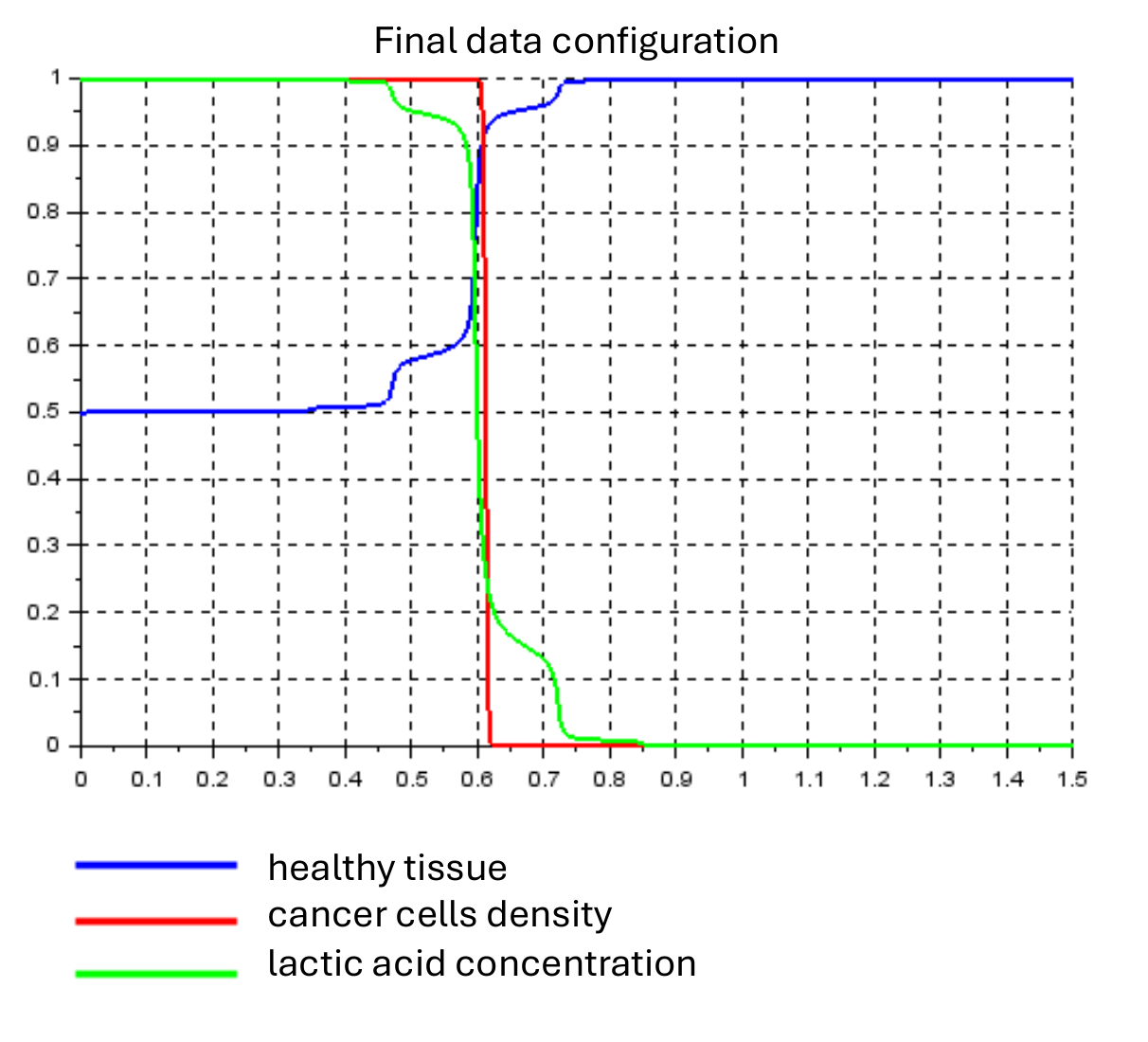}}\qquad
	\subfloat[][$d=1.5$]{\includegraphics[width=.35\textwidth]{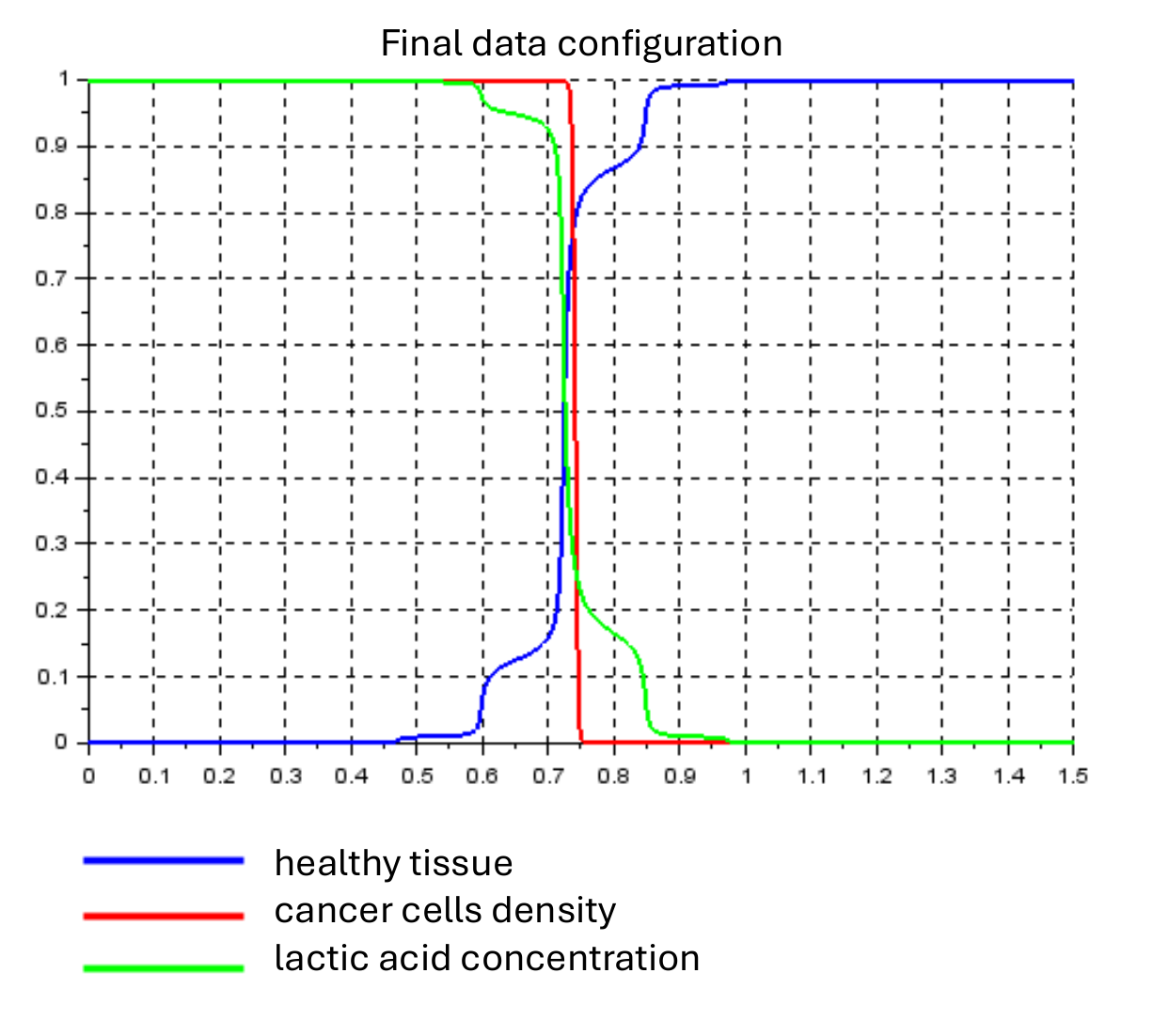}} \\
	\subfloat[][$d=30$]{\includegraphics[width=.35\textwidth]{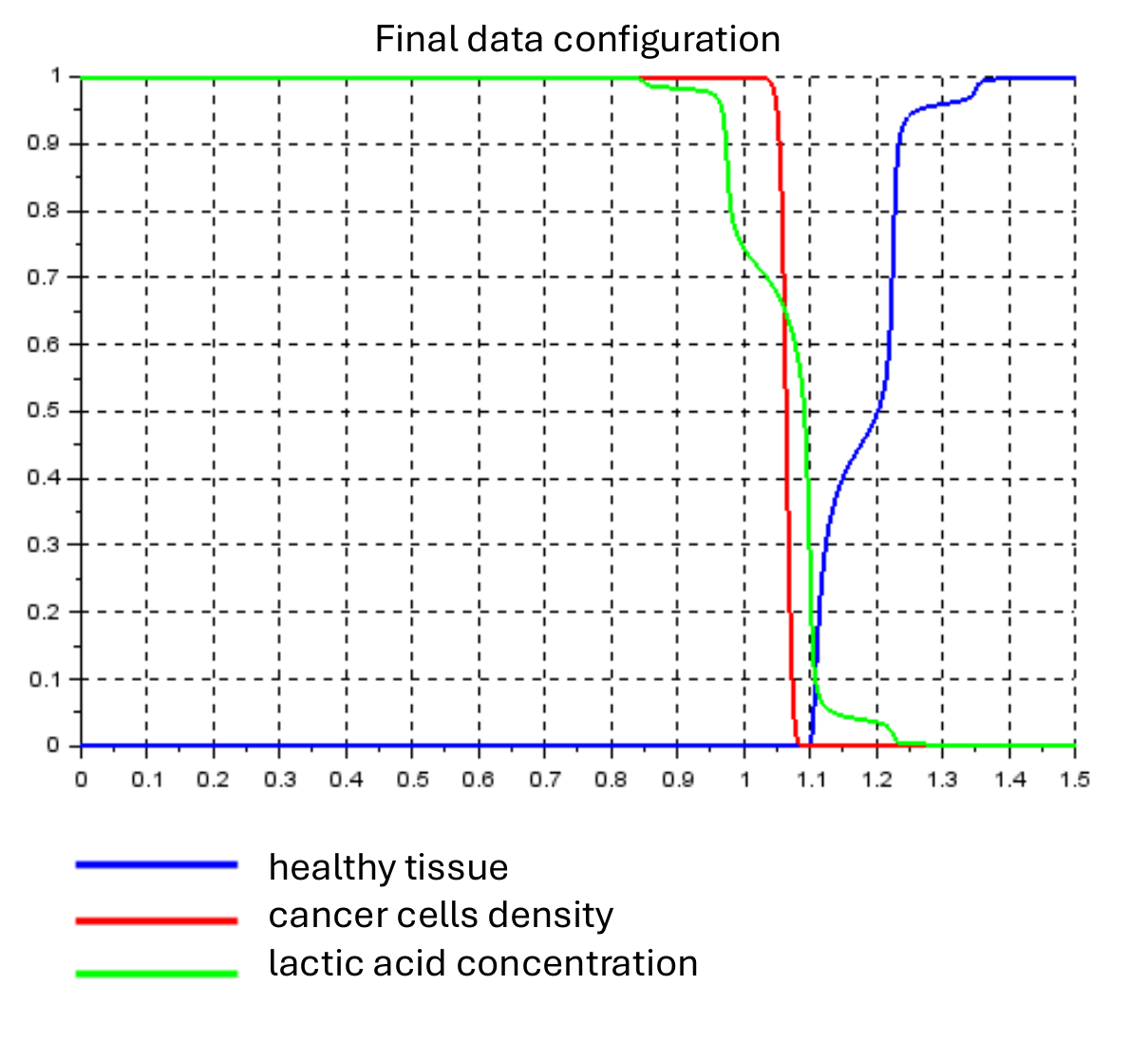}}\qquad
	\subfloat[][$d=60$]{\includegraphics[width=.35\textwidth]{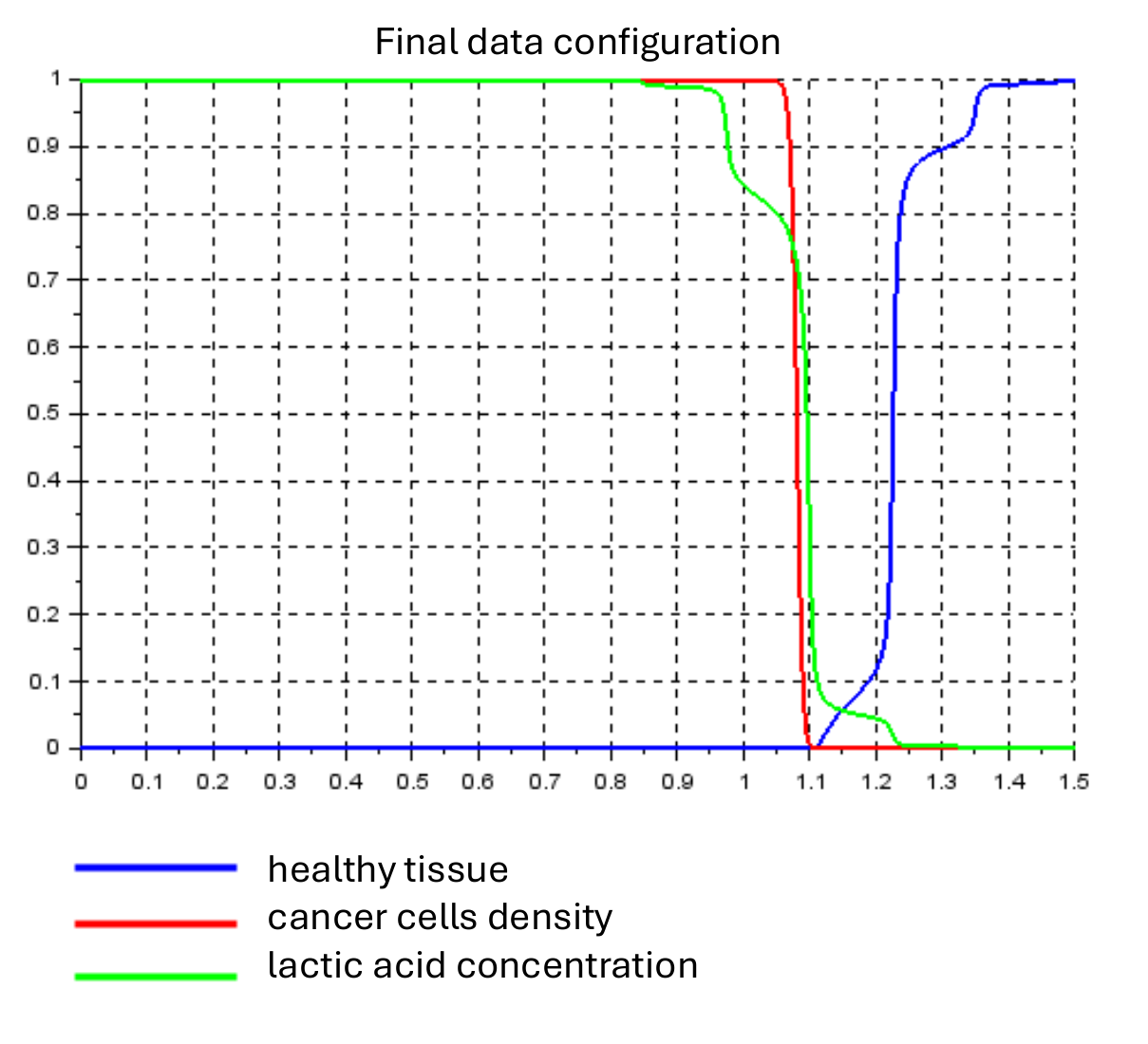}} 
	\caption{\footnotesize{Different configurations of the numerical solution in presence of periodic
	diffusion $A$ with frequency $\omega=50$ and uncontrolled growth $r=10$.}}
	\label{perconfig10}
\end{figure}

 {The choice of value of $r$ is arbitrary and it has no specific biological interpretation,
but the idea is to control the emerging patterns in such a different regime.
However, such parameter $r$ is truly significant because it is directly related to the duplication process:
hence it should be regarded as a measurable quantity in experiments both {\it in vitro} and {\it in vivo}. 
}

The results of the wave speed approximation are displayed in Fig.\ref{wavespeed5}  {(see Appendix \ref{appendix}),}
where we have used the final time $T=40$ to better appreciate the asymptotic behaviour.
 {Also,} we have  {widened} to the right the spatial interval in order to see the swipe of fronts.

Two  {major} considerations can be made. 
Firstly, it is evident that both $d$ and $r$ affect the propagation speed, 
and in particular the fronts propagate faster for larger values of $d$ and $r$. 
Secondly, we can observe that if $d$ is small (for instance, $d=0.5$ or $d=1.5$), namely if the final
configuration is heterogeneous or hybrid, then the wave speed does not approach an asymptotic threshold, 
but it seems almost periodic; moreover, the closer the final configuration is to a homogeneous regime, 
the more this phenomenon attenuates. 
This behaviour could be explained by the fact that, for the same frequency and amplitude, 
if the diffusion has greater intensity, then the acid particles are so fast  {that} their velocity is not affected 
by the tissue inhomogeneity; on the contrary, if the acid diffusion is weak, then particles are possibly 
entrapped in pores and their speed trend is oscillatory. \\

	
\subsection{Is homogenization possible in the long run?}\label{homogenisation}

In this concluding  {subsection,} we make an attempt to determine a  {constant}
diffusion rate $A$ which may approximate a periodic diffusive trend in the long run. 
We take a solution to {\sf GG} system \eqref{sistema GGbis} with periodic function  {$A=A(x)$}
and we look for a   {positive constant $A_{\textrm{\tiny eff}}>0$,} called \emph{effective} $A$, 
such that the new solution computed for this homogeneous diffusion possibly exhibits the same 
asymptotic wave speed as the original one. 
	
Since we are dealing with periodic diffusion functions, the harmonic mean seems a reasonable choice
to predict an approximate value --with respect to other types of averages such as the geometric or
weighted mean-- and for a positive function $A$ defined as follows
 {\begin{equation}\label{harmonic_mean}
	m_{\textrm{\tiny h}}(A):=\left(\frac{1}{p}\int_{0}^{p}\,\frac{dy}{A(y)}\right)^{-1}\,,
\end{equation}
where $p$ is the period of the function $A$ 
(e.g., for the sinusoidal function \eqref{periodico1}, there holds $p=2\pi/\omega$).}
Then, we compare the wave speed trends for different periodic diffusions $A$ with the corresponding
ones obtained  {by setting $A_{\textrm{\tiny eff}}:=m_{\textrm{\tiny h}}(A)$ assuming as correct 
the LeVeque--Yee formula \eqref{levequeyee}.}
	
First, we analyse the case of piecewise constant periodic $A$. 
We take the function $b: \mathbb{R}\longmapsto\mathbb{R}$ defined as
\begin{equation}\label{b}
	b(y)=\left\{\begin{aligned}
	\alpha_1 	& \quad & \text{if $y\in(0,\beta)$}\\
	\alpha_0 	& \quad & \text{if $y\in(\beta,1)$}
	\end{aligned}\right.
\end{equation}
where $0<\beta<1$ is the discontinuity point. 
Then, we extend the function \eqref{b} by periodicity and rescale it to the interval $[0,1]$ 
by using the following formula
\begin{equation}\label{pc}
	A_{\textrm{\tiny pc}}(x):=b\left({x}/{\epsilon}\right)\,,
\end{equation}
where $\epsilon$ is the number of oscillations inside the interval. 
In this case, we can compute explicitly the harmonic mean and we have
\begin{equation*}
	m_{\textrm{\tiny h}}(A_{\textrm{\tiny pc}})=\bigg(\int_{0}^{1}\,\frac{dy}{b(y)}\bigg)^{-1}
	=\bigg(\int_{0}^{\beta}\,\frac{dy}{\alpha_1}+\int_{\beta}^{1}\,\frac{dy}{\alpha_0}\bigg)^{-1}
	=\frac{\alpha_0\alpha_1}{\alpha_0\beta+(1-\beta)\alpha_1}\,.
\end{equation*}
We report the numerical wave speed approximations, obtained by using the LeVeque--Yee formula \eqref{levequeyee},
with different values of $\alpha$, $\beta$ and $\omega=2\epsilon$ in order to verify if 
 {$A_{\textrm{\tiny eff}}=m_{\textrm{\tiny h}}(A_{\textrm{\tiny pc}})$} 
provides an effective equivalent for the periodic diffusion function $A_{\textrm{\tiny pc}}$. 
	
In Fig.\ref{hom01}, we show the results obtained with $(\alpha_0,\alpha_1)=(0.01,1)$ and $\omega=100$,
together with the profile of $A$. 
\begin{figure}[htb]\centering
	\subfloat[][\emph{$d=0.5$}]{\includegraphics[width=.35\textwidth]{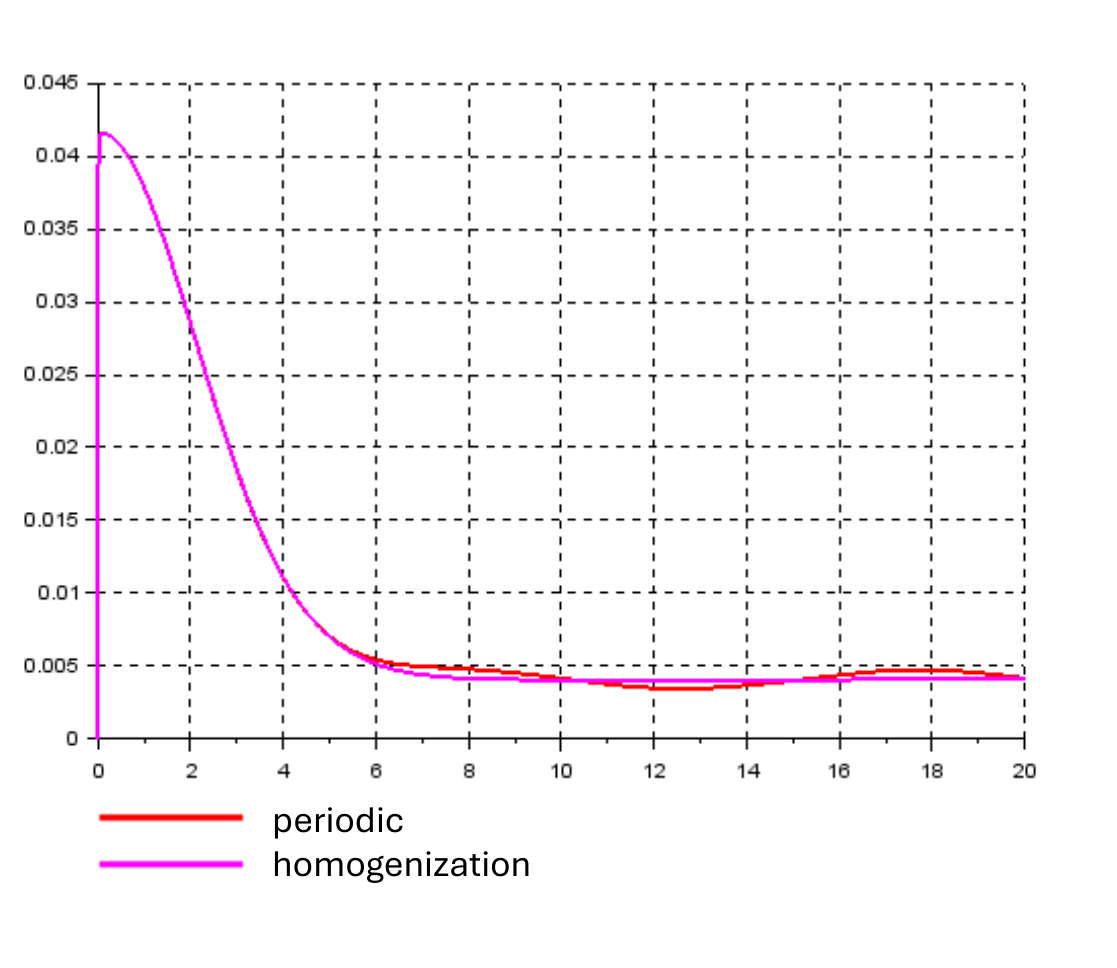}}\qquad
	\subfloat[][\emph{$d=1.5$}]{\includegraphics[width=.35\textwidth]{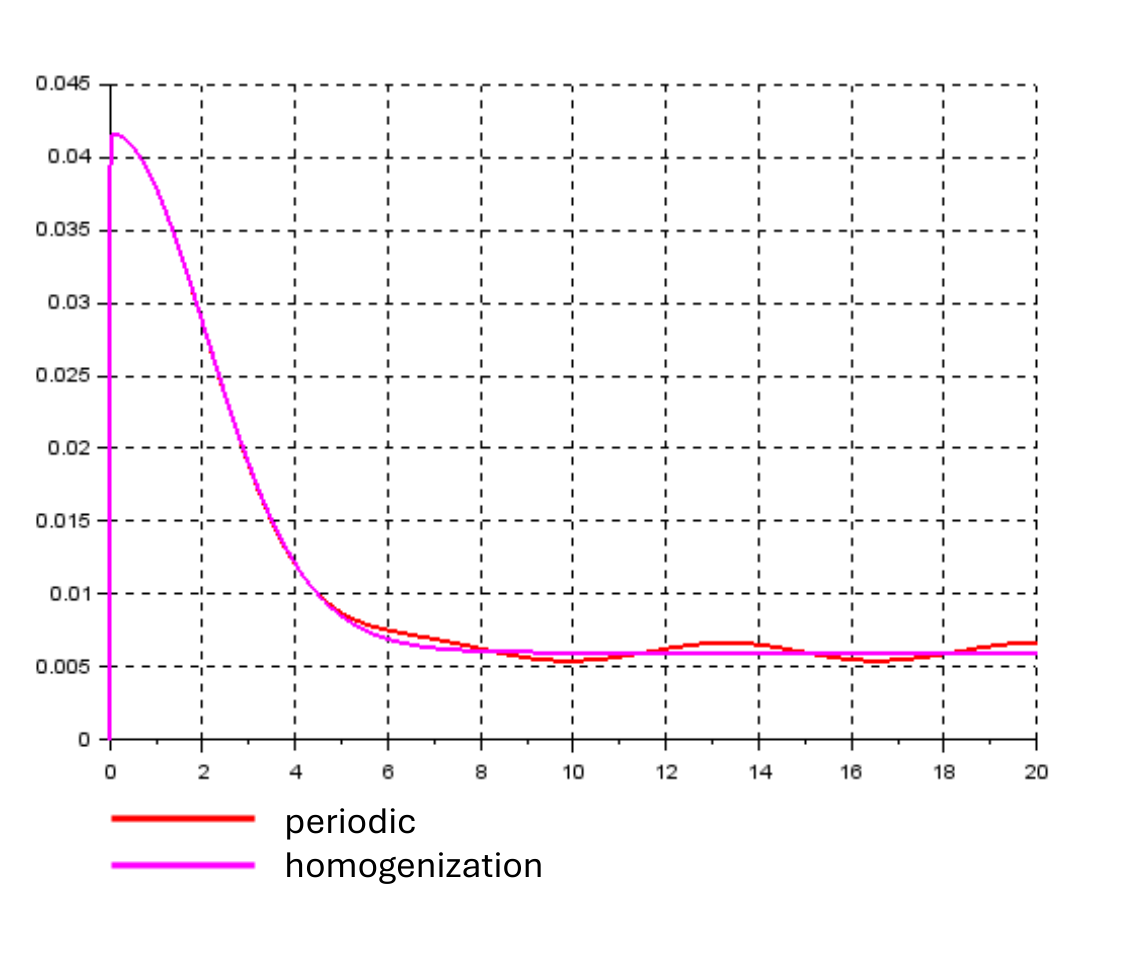}} \\
	\subfloat[][\emph{$d=30$}]{\includegraphics[width=.35\textwidth]{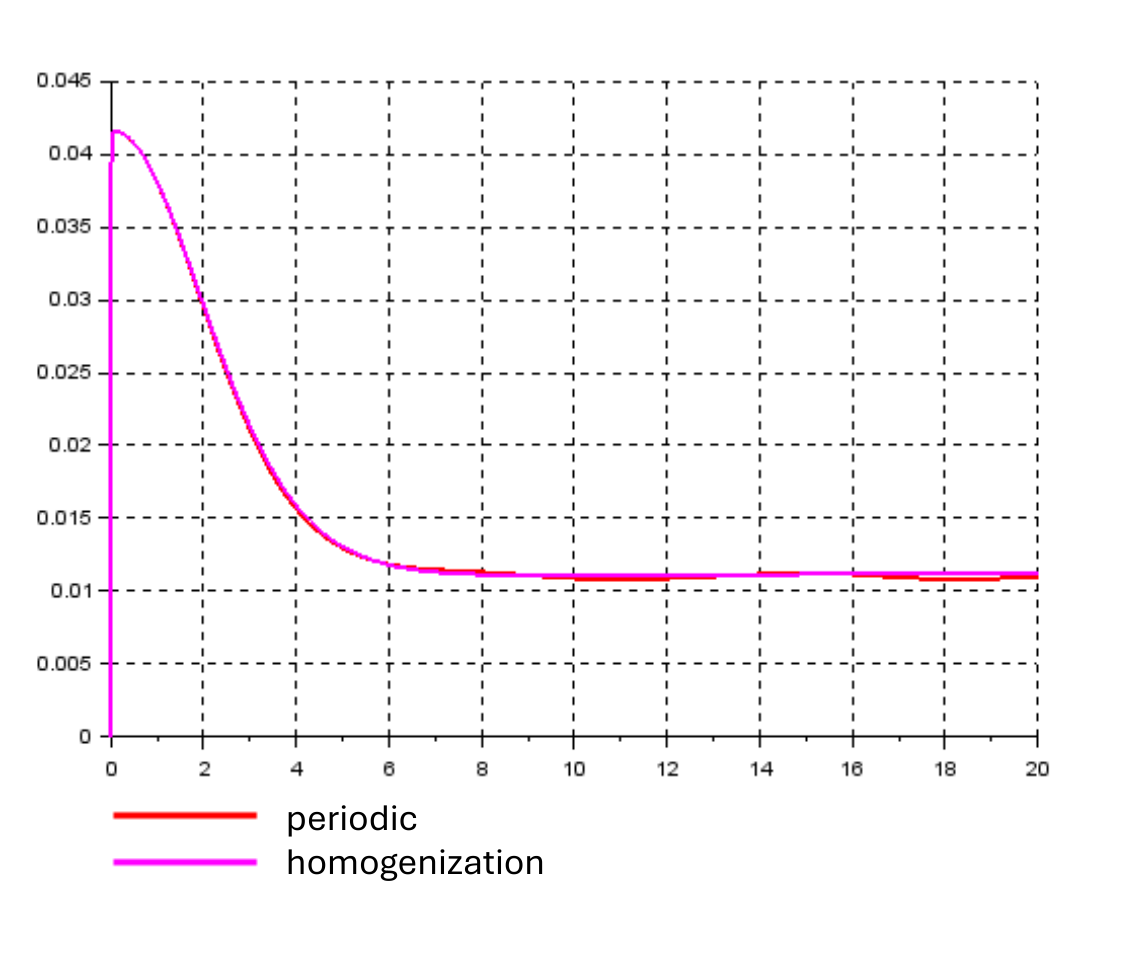}}\qquad
	\subfloat[][\emph{$d=60$}]{\includegraphics[width=.35\textwidth]{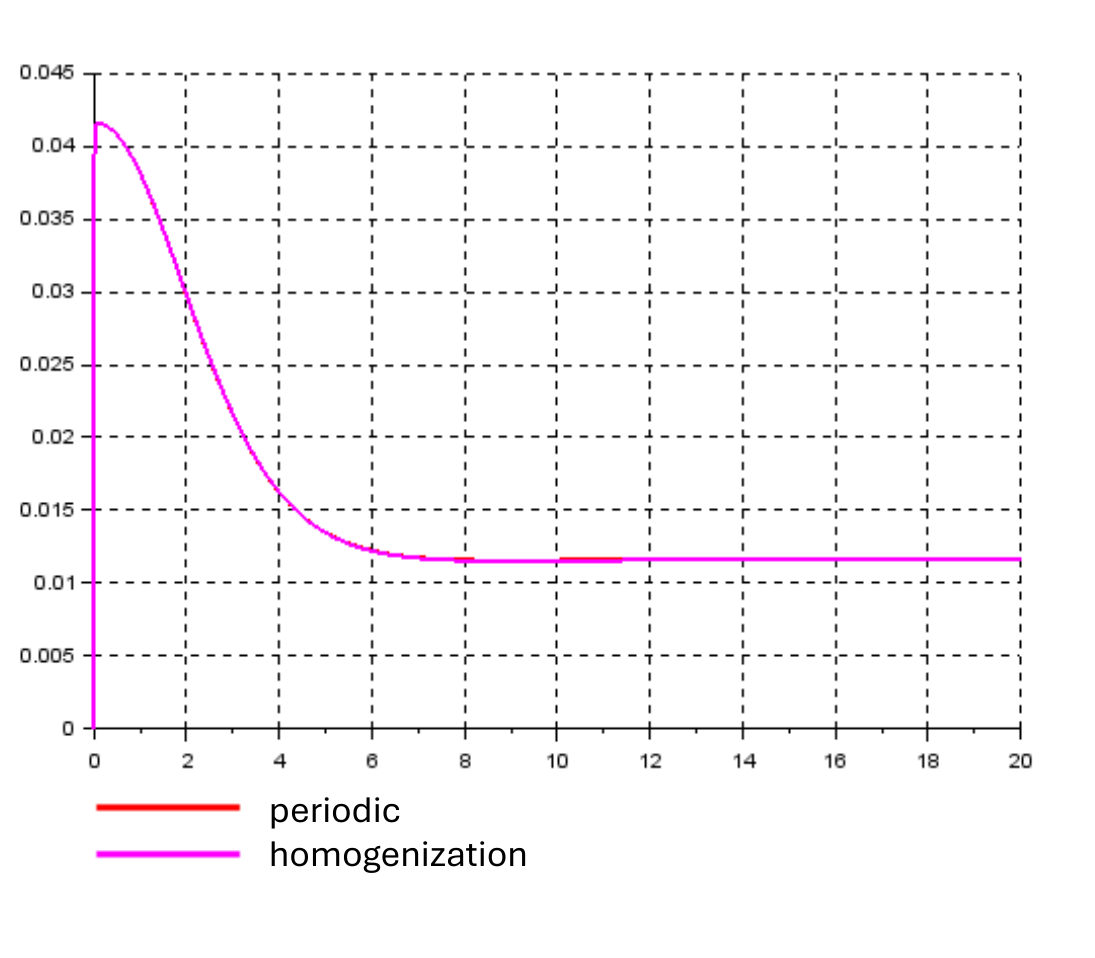}}
	\caption{\footnotesize{Comparison between wave speed trends for piecewise constant periodic diffusion 
	 {$A_{\textrm{\tiny pc}}$, see \eqref{pc}, with $(\alpha_0,\alpha_1)=(0.01,1)$, $\omega=100$ (red line) 
	and for the effective diffusion $A_{\textrm{\tiny eff}}=m_{\textrm{\tiny h}}(A_{\textrm{\tiny pc}})$ (magenta line).} }}
	\label{hom01}
\end{figure}
We notice that for the cases $d=0.5$ and $d=1.5$ the wave speed does not converge asymptotically
to any value, but rather oscillates between a minimum and a maximum, while $A_{\textrm{\tiny eff}}$
seems to tend to the mean value between these two.
 {In some cases (e.g., for large $d$) only the magenta line is readable.
This is because the two graphs overlap and the red line is drawn before the magenta one.}

Interesting considerations can be made  when we reduce the amplitude of oscillations $\omega=50$
--taking, for example,  {$(\alpha_0,\alpha_1)=(0.95,1)$ or rather $(\alpha_0,\alpha_1)=(0.4,0.6)$--}
because homogenisation occurs for every value of $d$ considered.
The results are omitted and we refer to \cite{ScanMascSime20} for an extensive set of simulations.

Now we consider the sinusoidal diffusion function $A$  {in \eqref{periodico1}.}
We point out that, in general, it is difficult to calculate analytically the formula for the harmonic mean
$m_{\textrm{\tiny h}}(A)$, thus we construct a numerical approximation of the harmonic mean of $A$
by using a quadrature formula \cite{Quarteroni}.
For the results illustrated in Fig.\ref{hom1}, we consider the periodic diffusion function $A$ in
\eqref{periodico1} with the same parameters  {as} Fig.\ref{hom01}, and we notice that, as for 
 {the} relative case with piecewise constant periodic $A$, homogenisation occurs only for $d\gg1$.
 {Similar conclusions can be drawn} when reducing the amplitude of oscillations, 
since  {simulations provide numerical evidence that} homogenisation occurs for all $d$, 
 {see \cite{ScanMascSime20}.}
\begin{figure}[bht]\centering
	\subfloat[][\emph{$d=0.5$}]{\includegraphics[width=.35\textwidth]{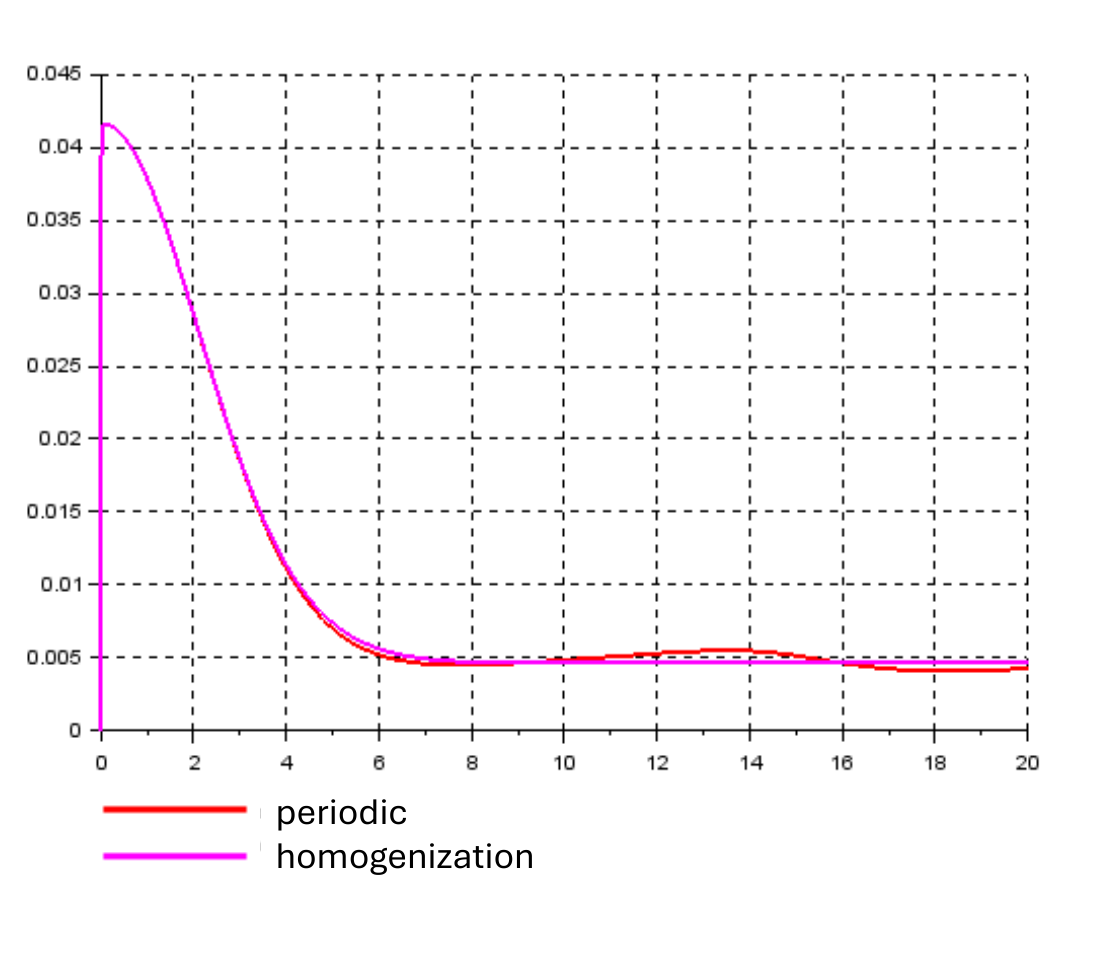}}\qquad
	\subfloat[][\emph{$d=1.5$}]{\includegraphics[width=.35\textwidth]{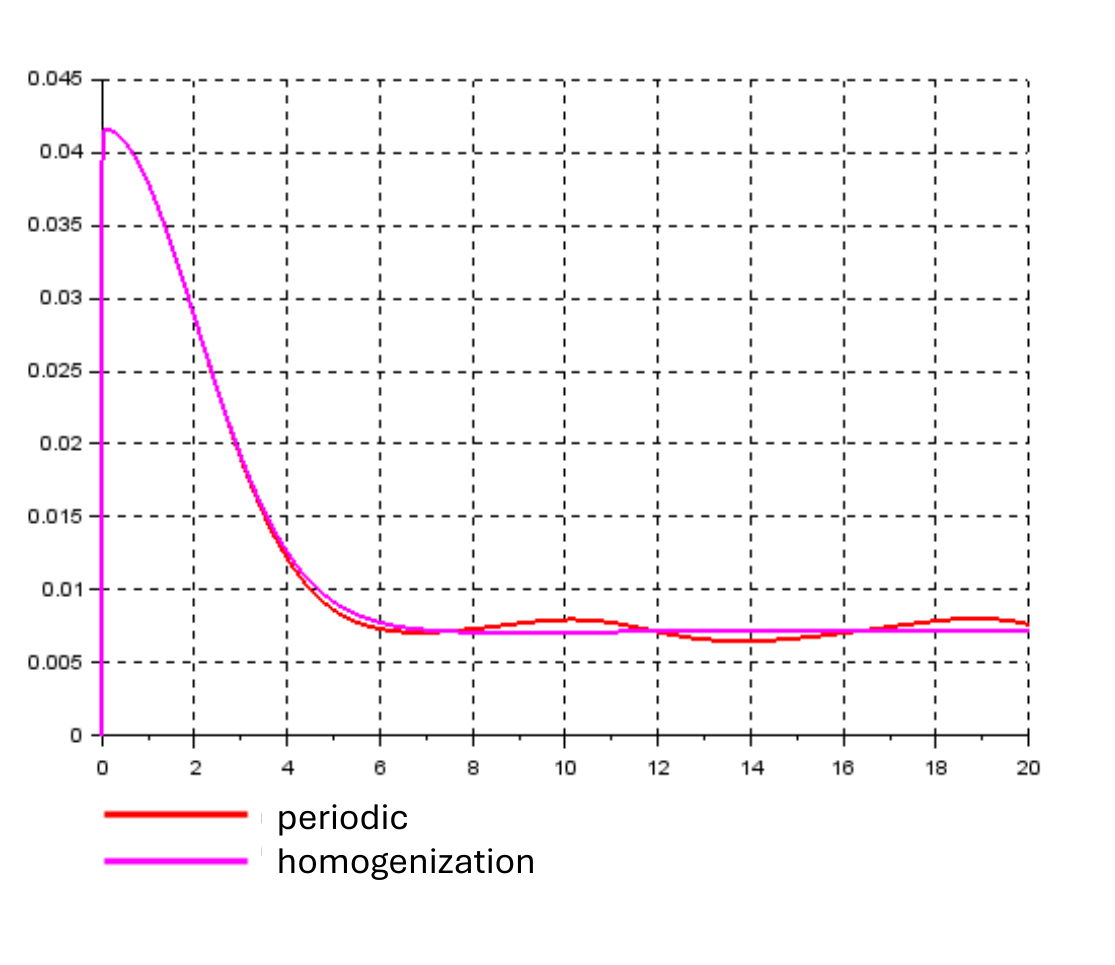}} \\
	\subfloat[][\emph{$d=30$}]{\includegraphics[width=.35\textwidth]{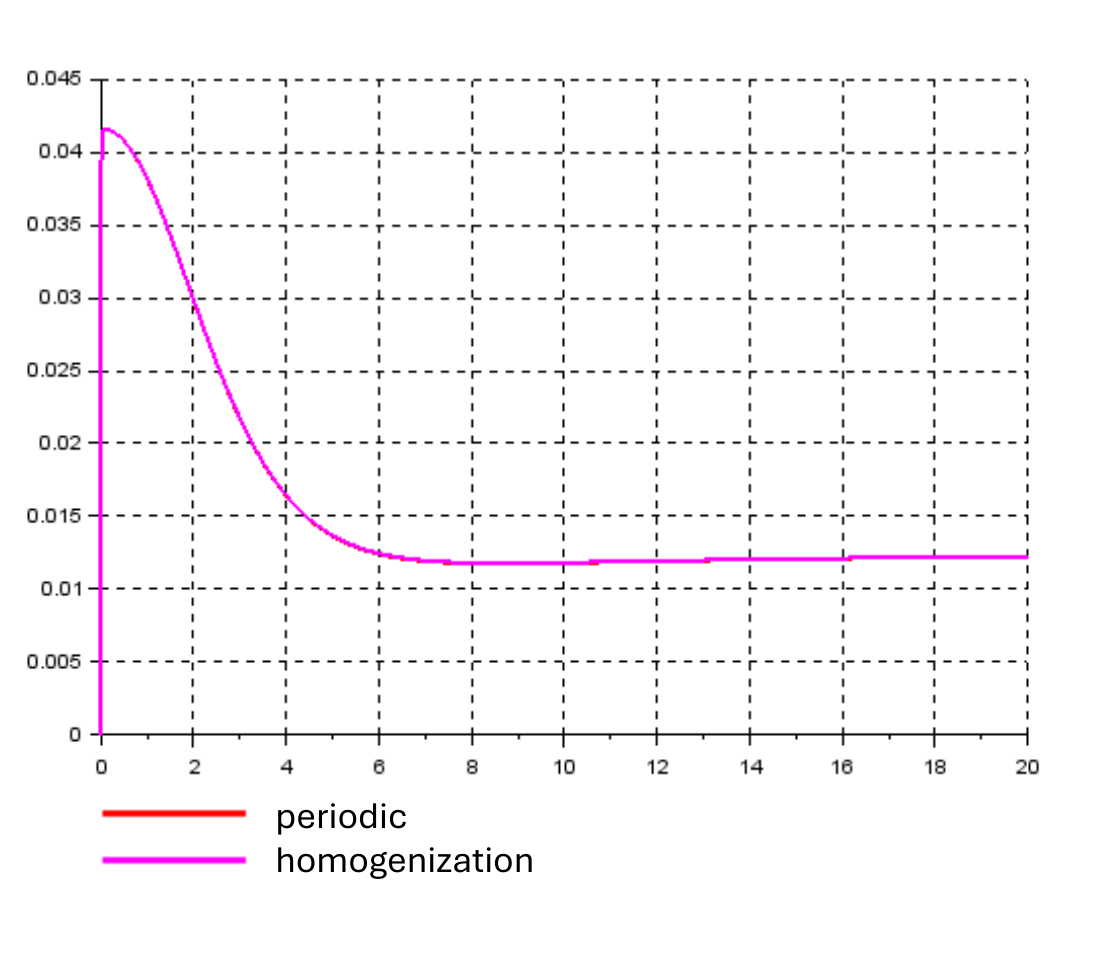}}\qquad
	\subfloat[][\emph{$d=60$}]{\includegraphics[width=.35\textwidth]{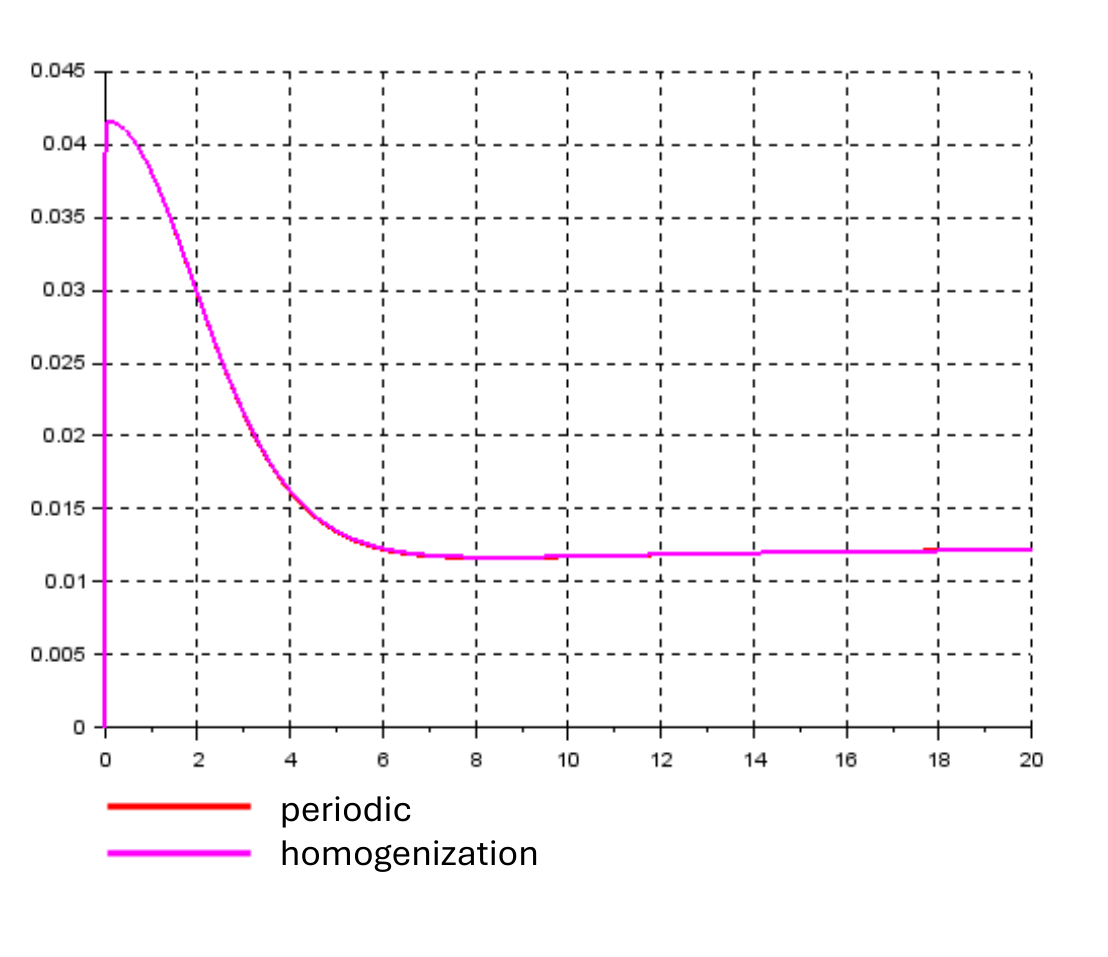}}
	\caption{\footnotesize{Comparison between wave speed trends for sinusoidal diffusion $A$, 
	 {see formula \eqref{periodico1}, with $(\alpha_0,\alpha_1)=(0.01,1)$, $\omega=100$ (red line)}
	and for the  {effective diffusion $A_{\textrm{\tiny eff}}$ (magenta line) given}
	by the harmonic mean $m_{\textrm{\tiny h}}(A)$.}}
	\label{hom1}
\end{figure}

Looking at the previous simulations, it is evident that a good agreement with $A_{\textrm{\tiny eff}}$ equal to the
harmonic mean of the periodic function $A$ always occurs for homogeneous configurations (i.e. when $d$ 
is such that a spatial interstitial gap appears). 
This fact can be justified by observing that the interstitial gap is actually an area where both $u$ and $v$ become
null, therefore in this case the equation for $v$ in {\sf GG} system \eqref{sistema GGbis} is a Fisher--KPP equation. 
It can be proven that homogenization for Fisher--KPP equations is possible \cite{elsmaily}, thus we can
expect  a similar phenomenon to occur also for the model {\sf GG} \eqref{sistema GGbis} when we 
consider a homogeneous configuration.
 {At the present day, such an issue has not yet been explored.}

For the sake of reproducibility of the numerical results illustrated and discussed in this section, in Table \ref{tavola3}, 
we  {collect} all parameters with the relative outcome in terms of homogenisation,
to be used as  {a} benchmark for comparative investigations.

\begin{table}[htp]
	\centering
	\begin{tabular}{c|c|c|c|c|c}
	$d$   	&  $\omega$ 	& $\alpha_0$ 	& $\alpha_1$  	& p.wise const. & sinusoidal \\ \hline
	$0.5$   	&  100  		&  0.01  		& 1 			& NO  		& NO  \\
	$1.5$   	&  100  		&  0.01  		& 1 			& NO 		& NO  \\
	$30$   	&  100  		&  0.01  		& 1			& NO 		& HOM  \\
	$60$    	&  100  		&  0.01  		& 1 			& HOM  		& HOM  \\ \hline
	$0.5$    	&  50  		&  0.95  		& 1 			& HOM  		& HOM  \\
	$1.5$    	&  50  		&  0.95  		& 1			& HOM		& HOM  \\
	$30$    	&  50 		&  0.95  		& 1			& HOM		& HOM  \\
	$60$     	&  50			&  0.95		& 1			& HOM		& HOM  \\ \hline
	$0.5$	&  50			&  0.4		& 0.6			& HOM		& HOM  \\
	$1.5$	&  50			&  0.4  		& 0.6			& HOM		& HOM  \\
	$30$    	&  50			&  0.4		& 0.6			& HOM		& HOM  \\
	$60$     	&  50  		&  0.4  		& 0.6 		& HOM  		& HOM  \\
	\end{tabular}
	\caption{\footnotesize{Numerical parameters for the simulations 
	 {referring to $A$ given by \eqref{periodico1} and corresponding numerical experiments.}
	The values for $D=4\cdot10^{-5}$, $r=1$, $c=70$ are fixed and the spatial interval is $[0,1]$.}}
	\label{tavola3}
\end{table}

It is evident that there are two special cases for which the two profiles for $A$ have different behaviours
referring to homogenization: these are the ones  {for $(\alpha_0,\alpha_1)=(0.01,1)$,}
$\omega=100$ and $r=1$, see Fig.\ref{hom01} and Fig.\ref{hom1}(c).
The reasons  {seem} to be found in the final configurations: while in the piecewise constant case there is not
a purely homogeneous invasion, in the sinusoidal case the formation of the interstitial gap has already started.
\begin{figure}[]\centering
	\subfloat[][piecewise constant]{\includegraphics[width=.35\textwidth]{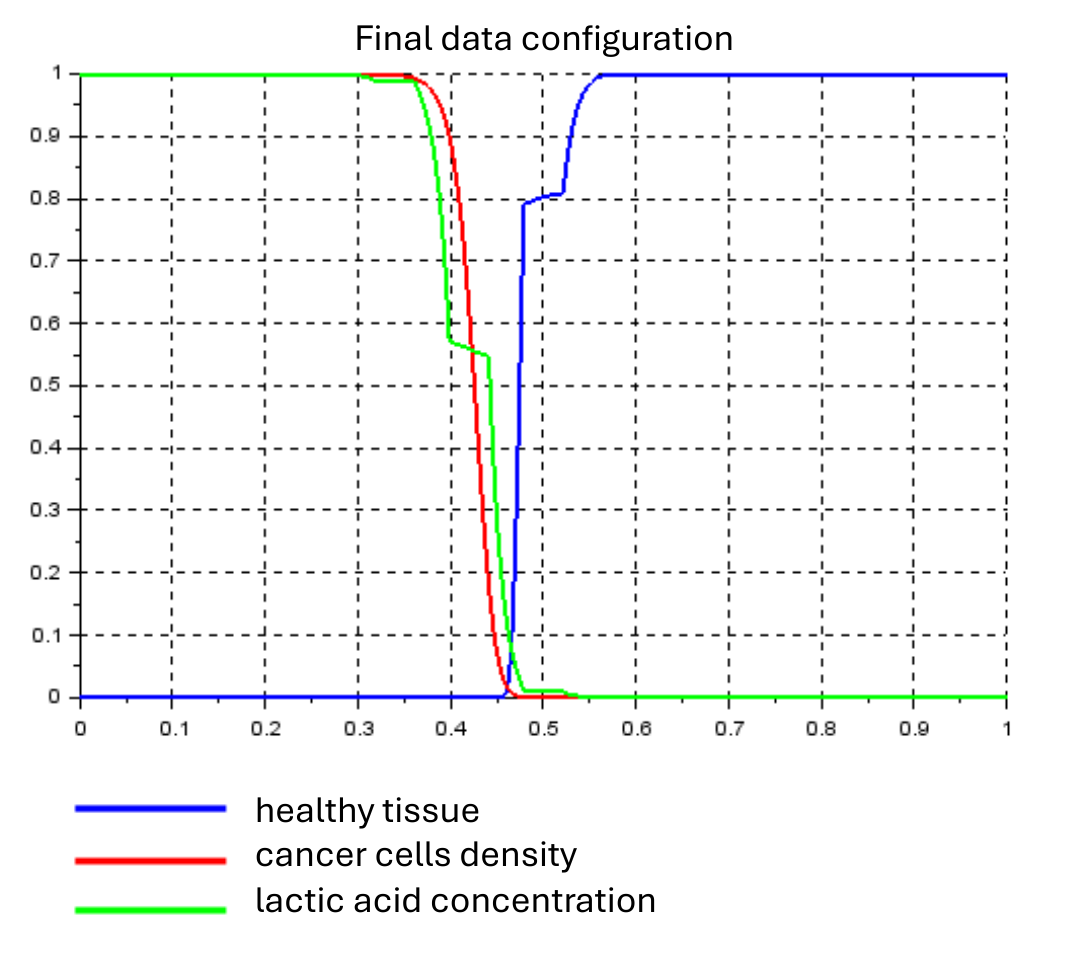}}\qquad
	\subfloat[][sinusoidal]{\includegraphics[width=.35\textwidth]{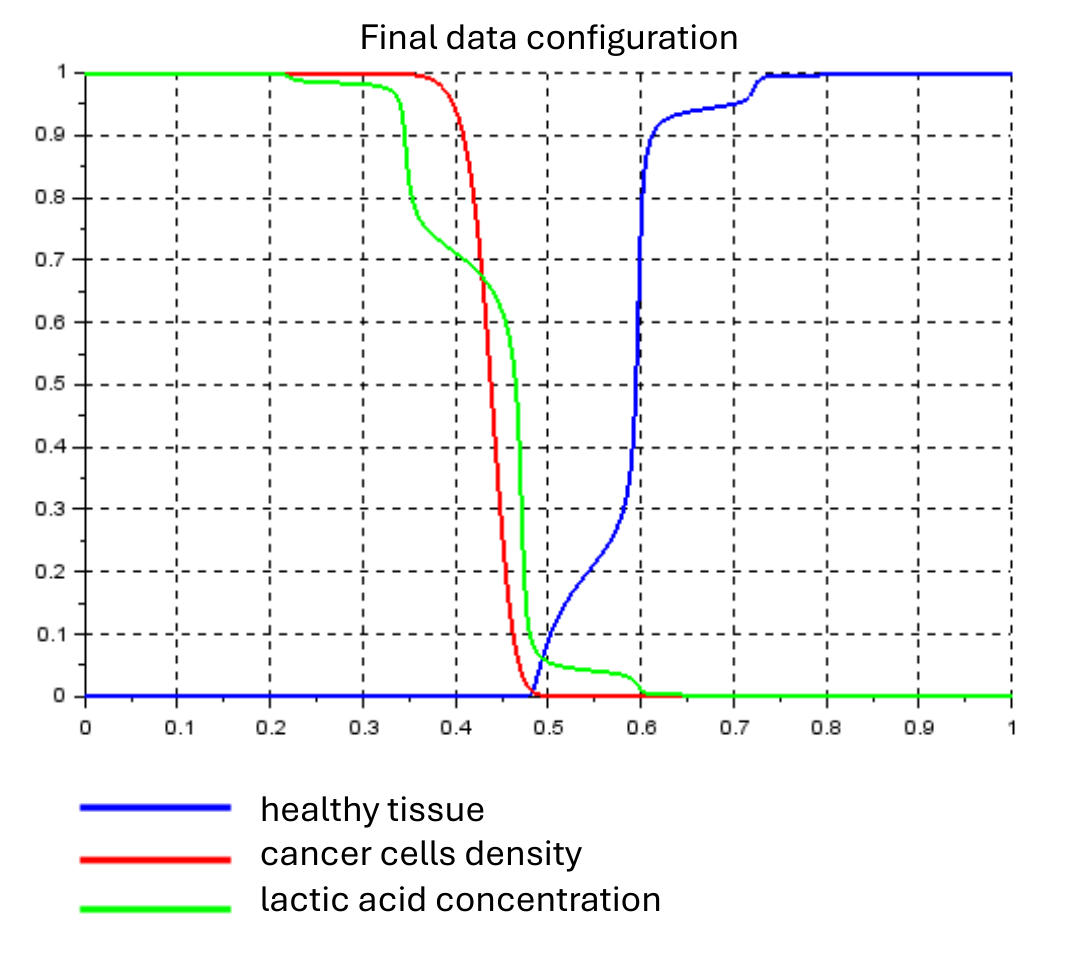}}
	\caption{\footnotesize{Case $d=30$, $r=1$: comparison between final configurations for piecewise constant,
	 {as in \eqref{pc},} and sinusoidal periodic diffusion $A$,  {see \eqref{periodico1}}, 
	in the case of different behaviour  {relative} to homogenisation.}}
	\label{hom}
\end{figure}

 {\section{Final comments and perspectives}}\label{sect:conclusions}

 {The presence of intra-tumour heterogeneity is increasingly recognized as one of the primary 
drivers of heterogeneity in patient outcomes.
Appropriate clinical procedures could profit by such elementary modeling identifying alternative targets 
which could complement the standard therapies and approaches.

In this work, we have explored the cancer invasion dynamics using a one-dimensional variation of the
Gatenby--Gawlinski model modified by means of a heterogeneous diffusive term for lactic acid, 
describing some kind of local ionic pH signalling. 
To explore the consequences of such scenario, we performed some numerical simulations using 
a finite volume method and considered for the acid diffusion a single jump piecewise constant function
\eqref{single-jump} and then a single-frequency periodic one \eqref{periodico1}.
From the computations with \eqref{single-jump}, see Subsection \ref{subsect:single_jump},
we learn that the formation of the spatial interstitial gap is affected by the amplitude of the jump.
In particular, numerical experiments show that a decrease in the inhomogeneity affects formation 
of the gap, according to a smaller cancer aggressiveness parameter $d$.

From the simulations where for the acid diffusion we focussed a periodic function, see Subsection \ref{subsect:periodic},
we found that the solution still have a front type behaviour where the frequency of the considered acid 
function does not produce effects on the wave speed, but influences the appearance of the interstitial gap. 
Also, we limited ourselves to the case of periodic fluctuations in the form of sinusoidal function \eqref{periodico1} 
and piecewise constant \eqref{pc}, our intention being to focus on the case of {\it random media} in some future 
works, \cite{Xin09}.}

 {These results are very preliminary and deserve further investigation.
The issue of homogenisation must be analysed in terms of analytical convergence, especially taking 
into consideration the large sensitivity of the model with respect to the parameter $d$. 
Furthermore, we emphasize that the use of a non-uniform mesh could better model the non-homogeneity
of human tissues, thus our numerical algorithm will be modified in order to be more consistent with the 
biological context.}

 {Let us stress that the results obtained in this work are basic to study the cancer diffusion mechanisms
also in multiple dimensions. 
In particular, the computational costs related to numerical simulations in two and three dimensions,
with real input data, can lead to long computational times.
For these reasons, HPC (High Performance Computing) techniques could be essential to produce
numerical results comparable with real clinical data. 
In this regard, a two dimensional code for the modified Gatenby--Gawlinski model proposed in this
work based on GPU (Graphical Processing Unit) architecture has been already developed by the 
authors and will be used in a joint work.}
\\

\section*{Acknowledgements}

 {This research work has been partially supported by the French National Research Agency (ANR) 
under the UCA JEDI Investments for the Future project --reference no. ANR-15-IDEX-01-- and by 
the  {European Union - NextGenerationEU Project} under the Italian Ministry of University and 
Research (MUR) through the National Centre for HPC, Big Data and Quantum Computing 
CN$\_$00000013-CUP:E13C22001000006.}

 {Chiara Simeoni is thankful to the Department of Mathematics G. Castelnuovo, Sapienza, Universit\`a 
di Roma (Italy) for the invitation during which this work has been conceived.
Corrado Mascia is thankful to the Laboratoire J.A. Dieudonn\'e,  Universit\'e C\^ote d'Azur, Nice (France)
for the invitation during which this work has been completed.}


\appendix

 {
\section{Appendix. More numerical simulations}\label{appendix}

Here, some numerical simulations for different choices of parameters $(\alpha_0,\alpha_1)$,
$\omega$ and $r$ are considered. 
Detailed commentaries are contained in the previous Sections.

\begin{figure}[bht]\centering
	\subfloat[][$d=0.5$]{\includegraphics[width=.32\textwidth]{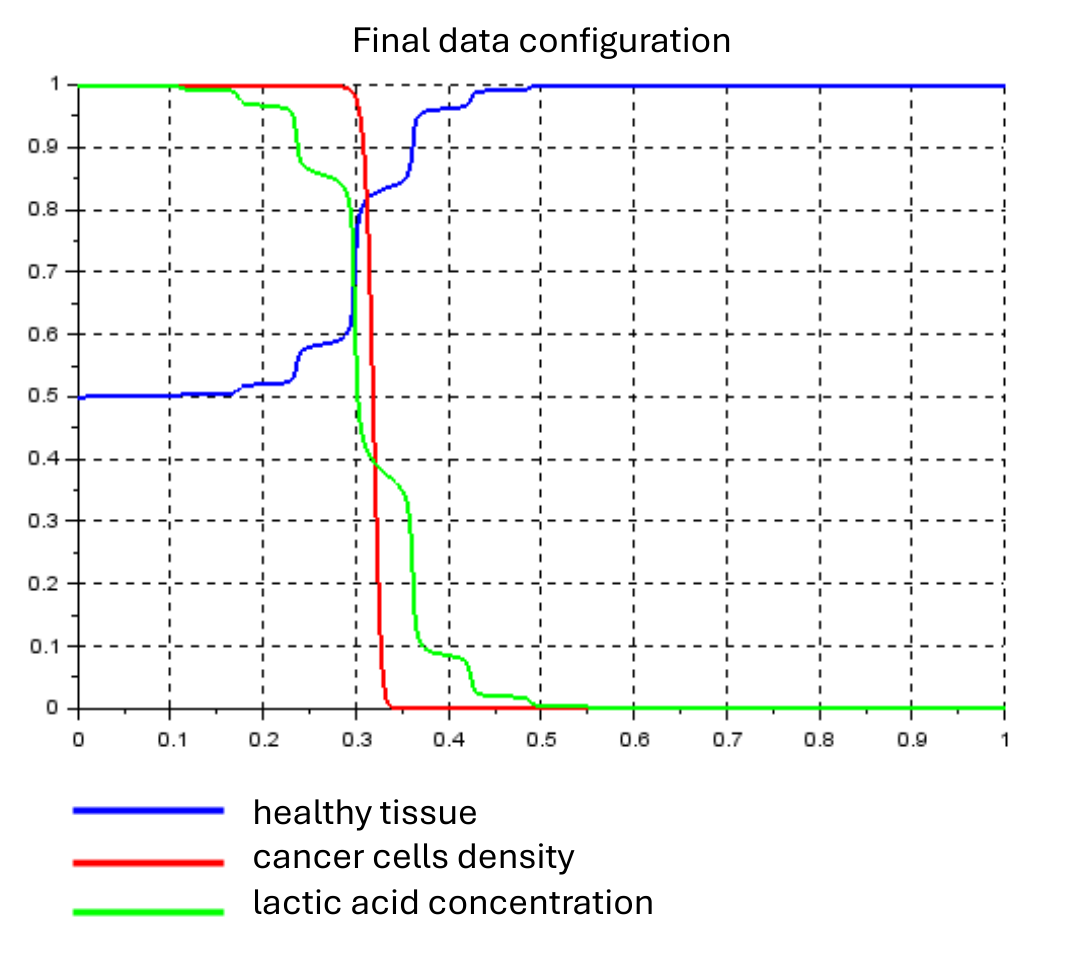}}\qquad
	\subfloat[][$d=1.5$]{\includegraphics[width=.32\textwidth]{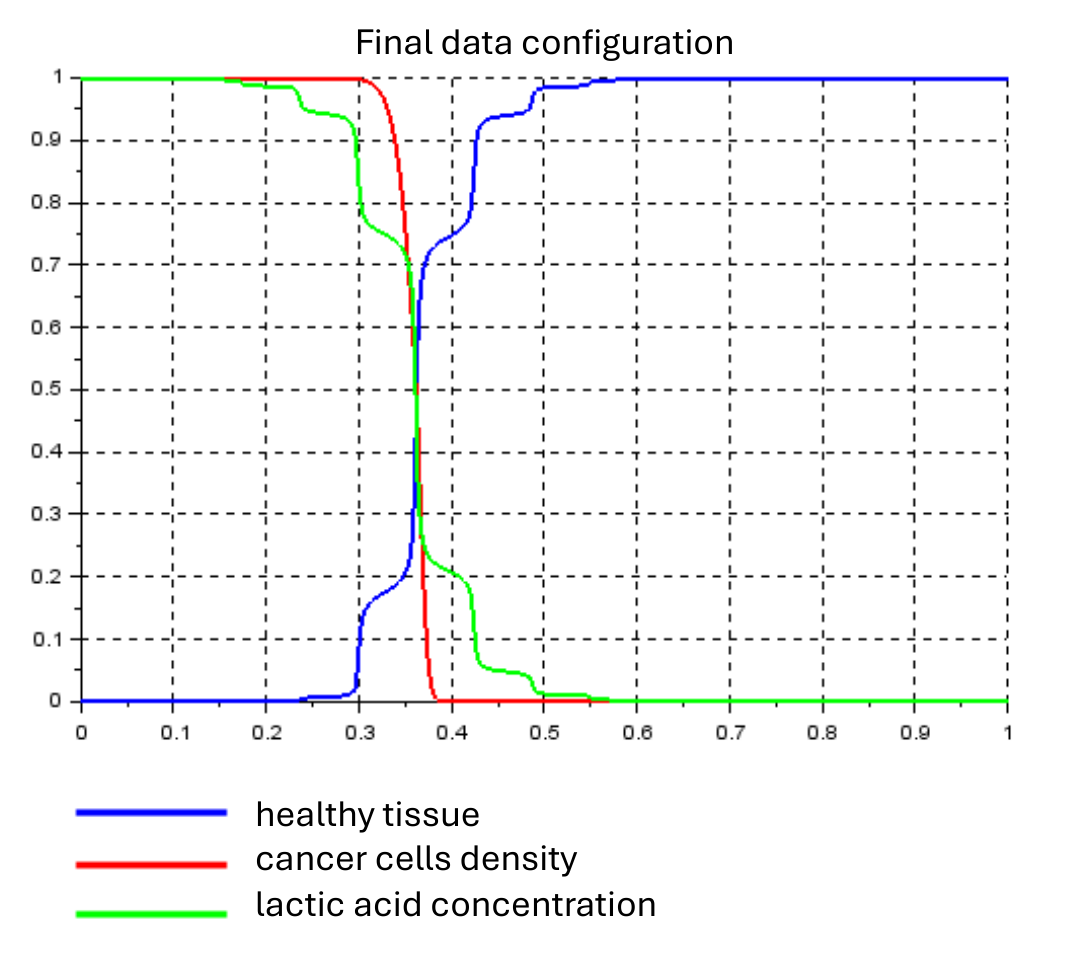}} \\
	\subfloat[][$d=30$]{\includegraphics[width=.32\textwidth]{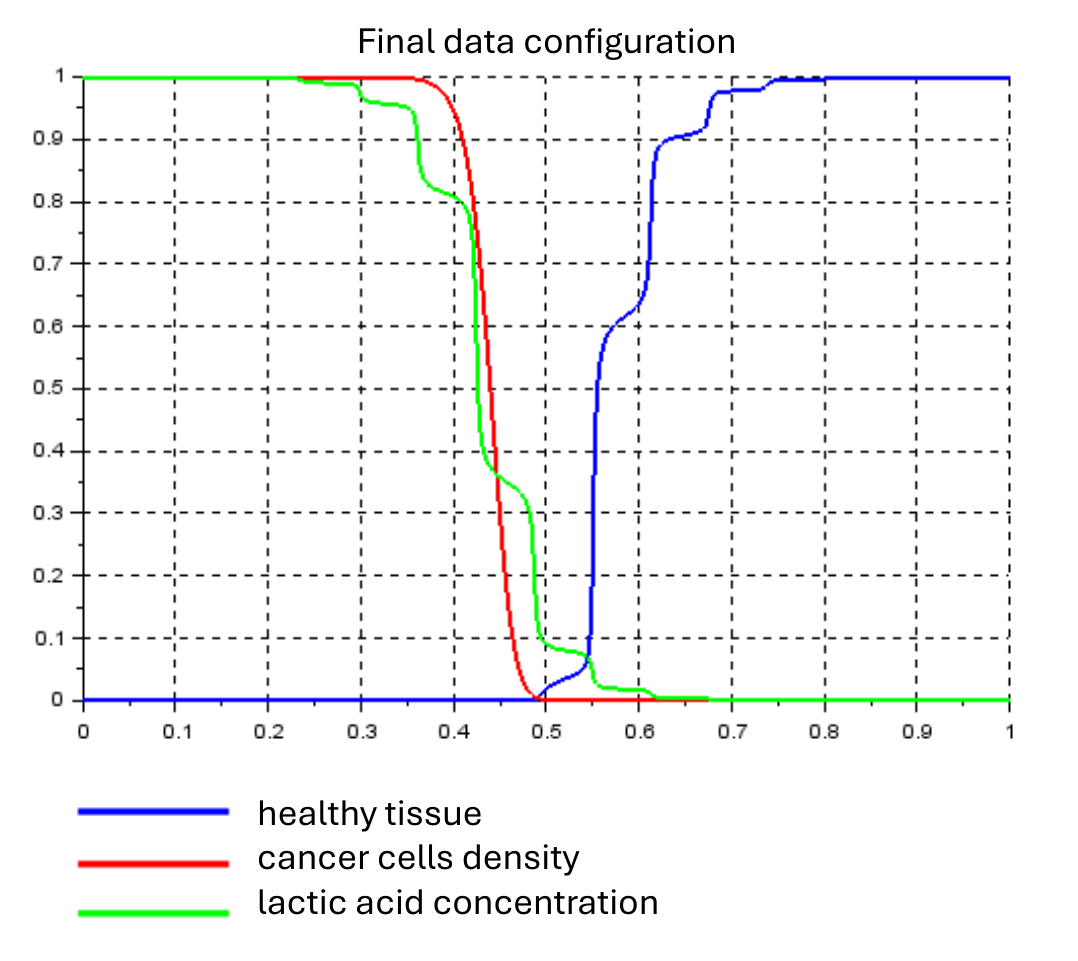}}\qquad
	\subfloat[][$d=60$]{\includegraphics[width=.32\textwidth]{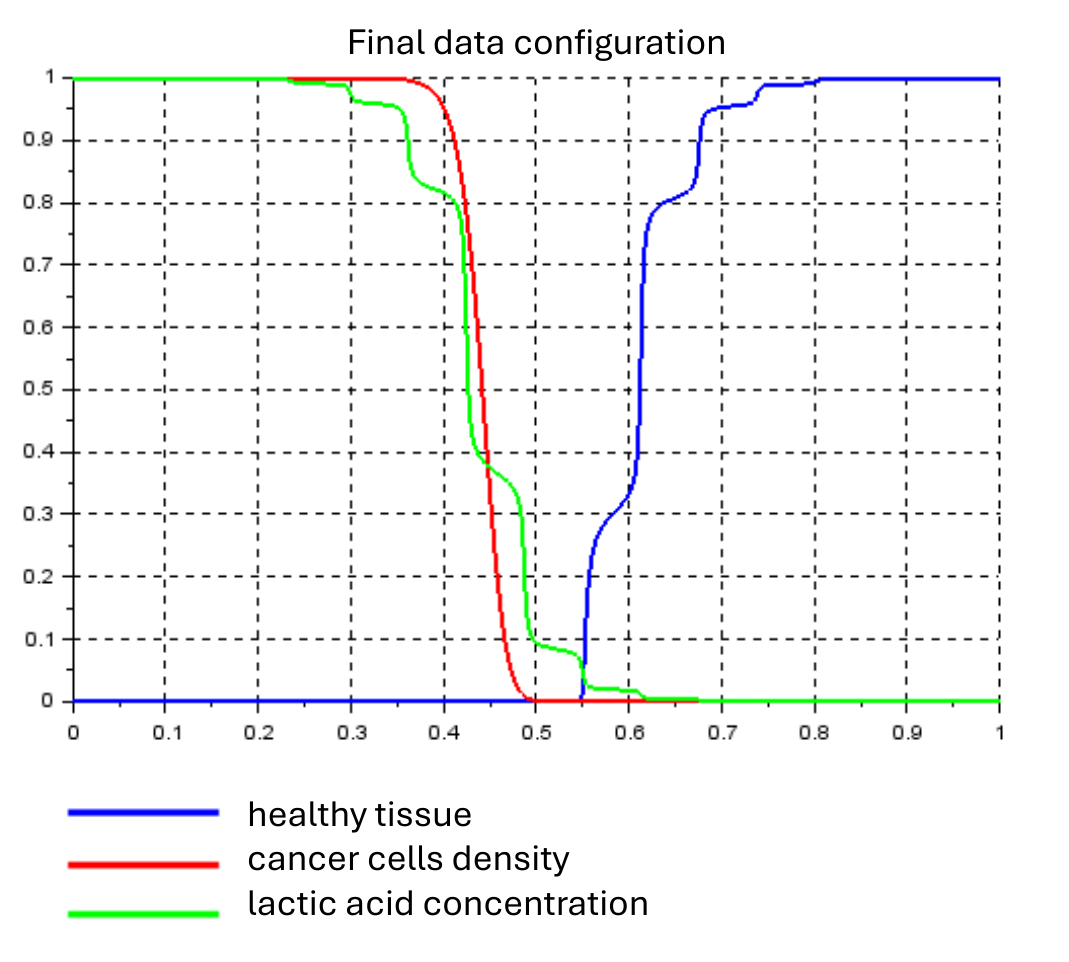}}
	\caption{\footnotesize{Different configurations of the numerical solution in presence of periodic diffusion $A$, 
	 {see formula \eqref{periodico1},}
	with frequency $\omega=100$ and existence of the spatial interstitial gap within the homogeneous invasion (c)-(d).}}
	\label{perconfig4}
\end{figure}

\begin{figure}[bht]\centering
	\subfloat[][$d=0.5$]{\includegraphics[width=.35\textwidth]{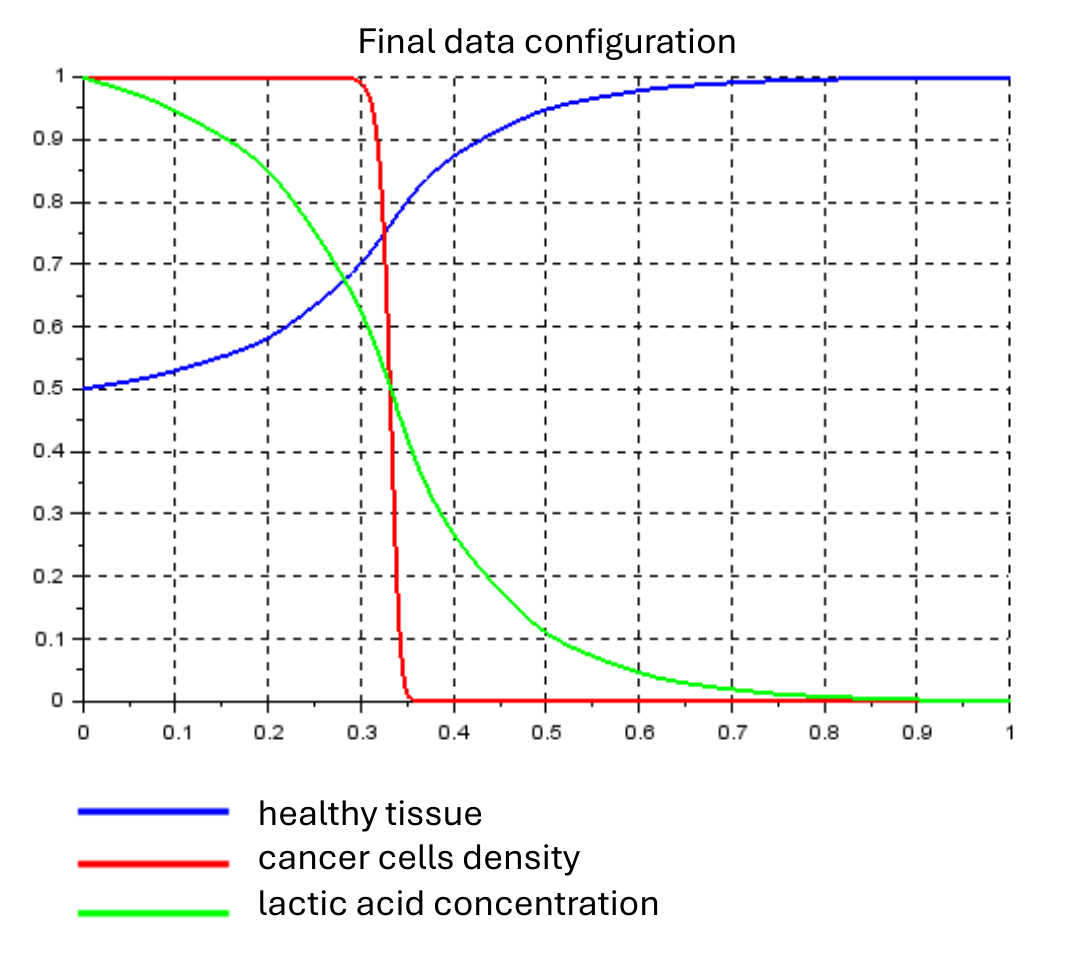}}\qquad
	\subfloat[][$d=20$]{\includegraphics[width=.35\textwidth]{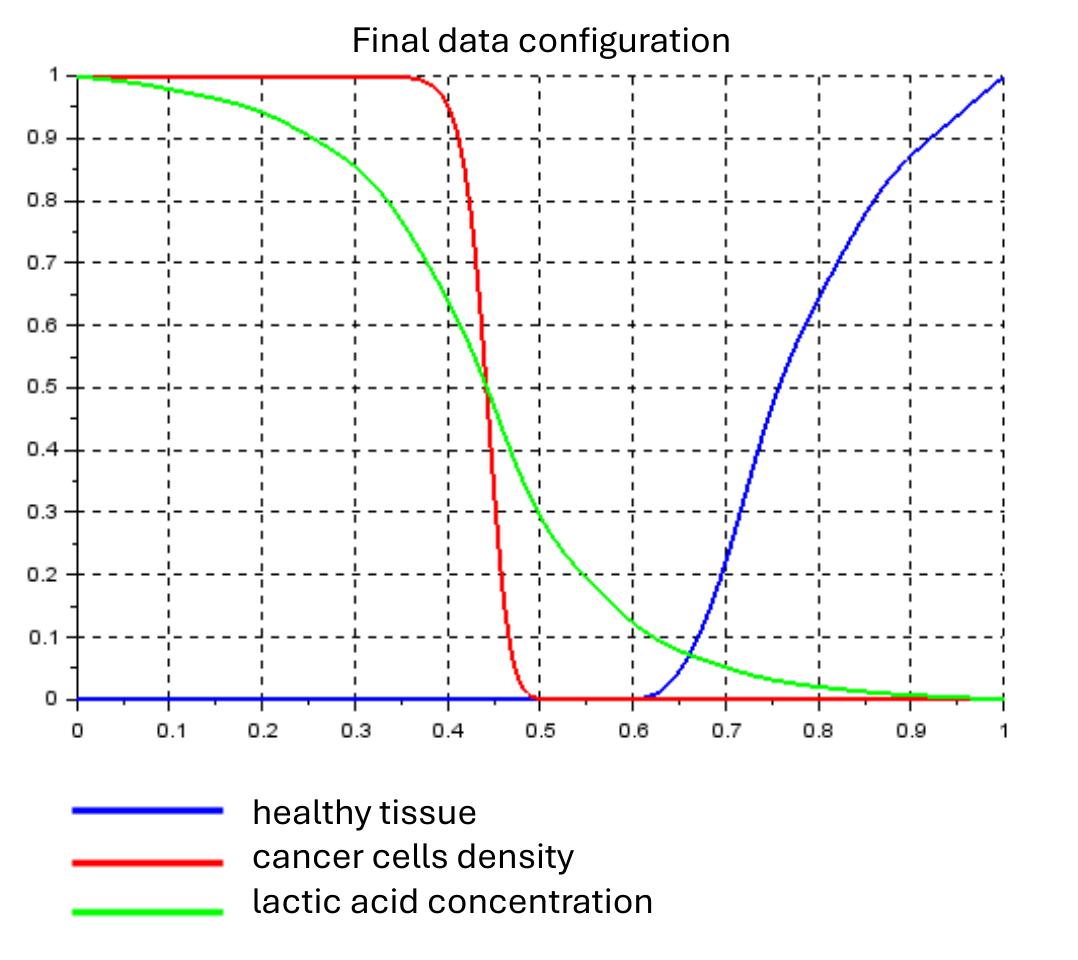}}
	\caption{\footnotesize{ {Case $\omega=50$ and $(\alpha_0,\alpha_1)=(0.8,1)$:
	heterogeneous (left) and homogeneous invasions (right).}}}
	\label{perconfig8}
\end{figure}

\begin{figure}[bht]\centering
	\subfloat[][$d=0.5$]{\includegraphics[width=.35\textwidth]{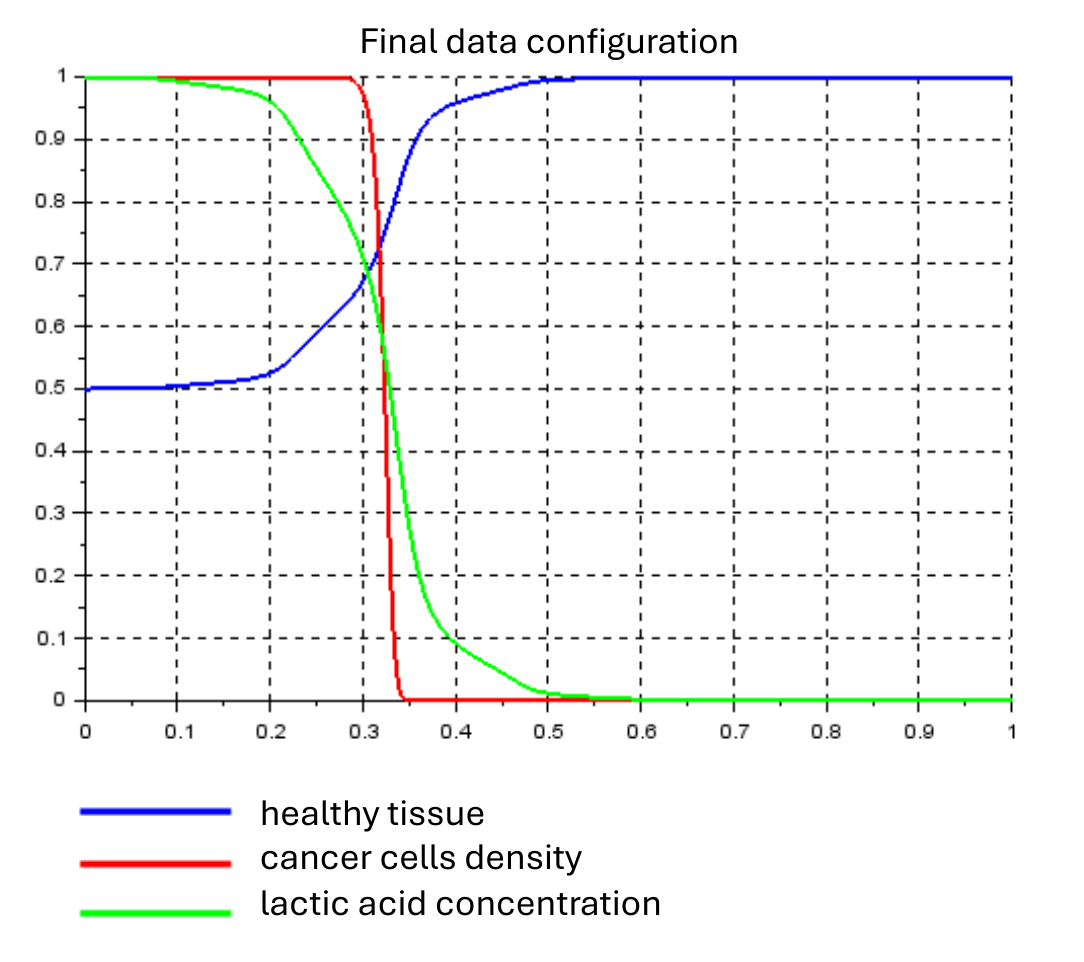}} \qquad
	\subfloat[][$d=20$]{\includegraphics[width=.35\textwidth]{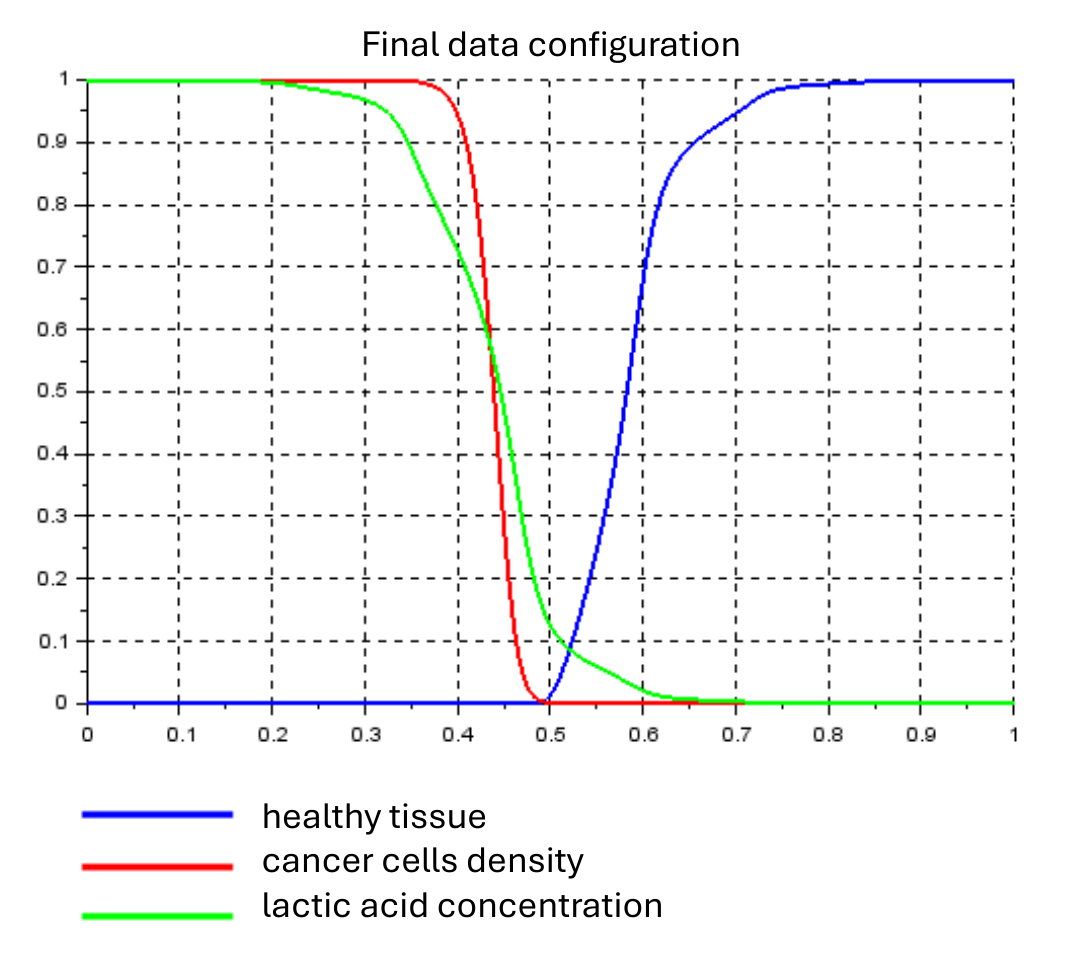}}
	\caption{\footnotesize{ {Case $\omega=50$ and $(\alpha_0,\alpha_1)=(0.1,0.3)$:
	heterogeneous (left) and homogeneous invasions (right).}}}
	\label{perconfig11}
\end{figure}

\begin{figure}[bht]\centering
	\subfloat[][$d=0.5$]{\includegraphics[width=.35\textwidth]{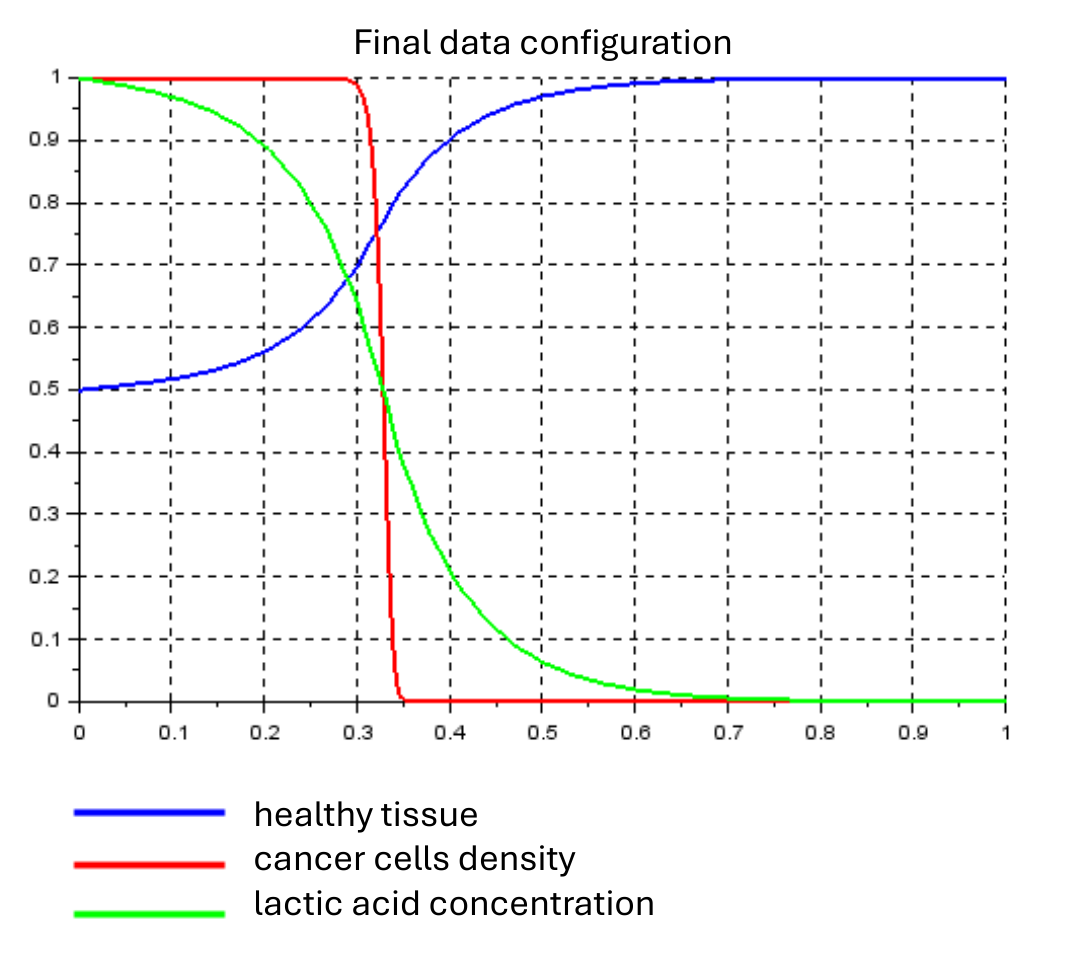}}\qquad
	\subfloat[][$d=20$]{\includegraphics[width=.35\textwidth]{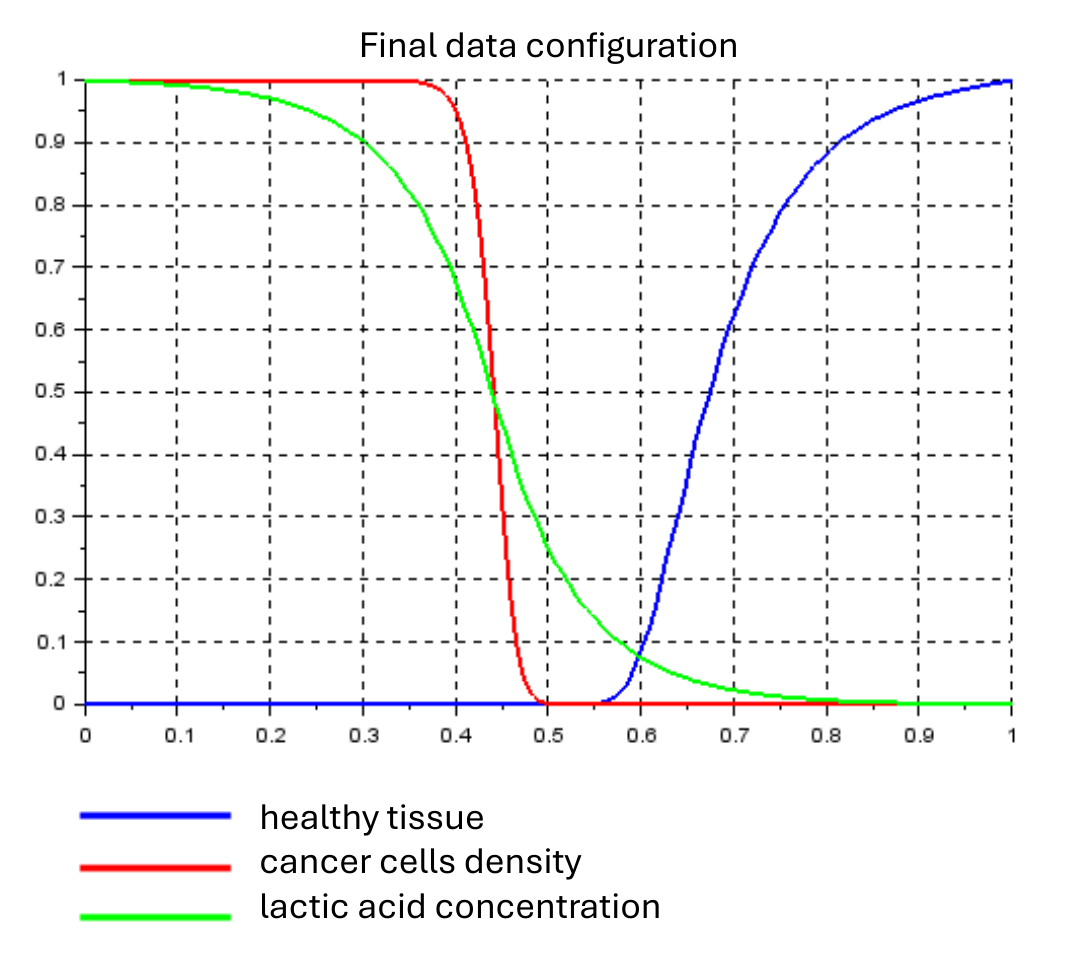}}
	\caption{\footnotesize{ {Case $\omega=50$ and $(\alpha_0,\alpha_1)=(0.95,1)$: 
	heterogeneous (left) and homogeneous invasions (right).}}}
	\label{perconfig12}
\end{figure}

\begin{figure}[bht]\centering
	\includegraphics[width=.5\textwidth]{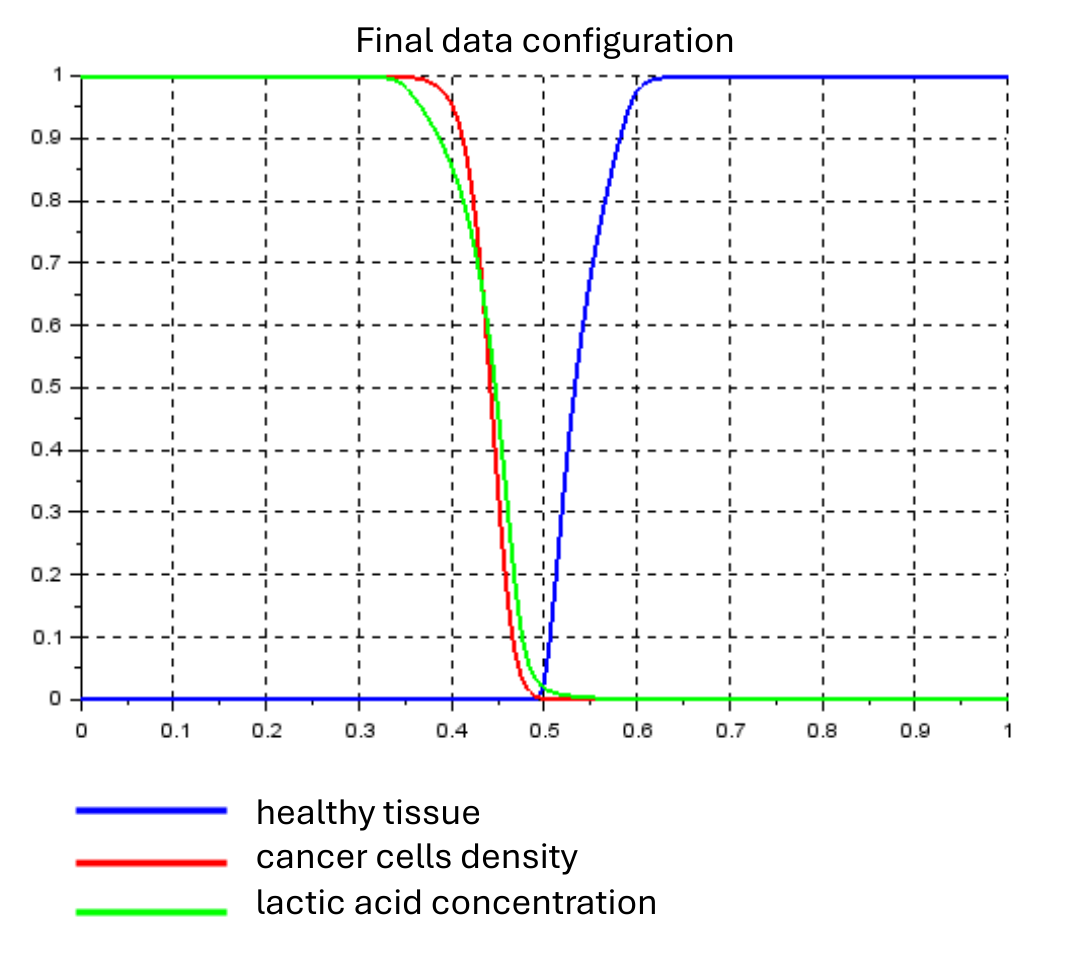}
	\caption{\footnotesize{ {The interstitial gap is absent at large regime of $d$ ($\omega=200$, 
	$(\alpha_0,\alpha_1)=(0.01,0.06)$ and $d=200$).}}}
	\label{perconfig13}
\end{figure}

\begin{figure}[bht]
	\centering
	\subfloat[][$r=1$ and $d=0.5$]{\includegraphics[width=.45\textwidth]{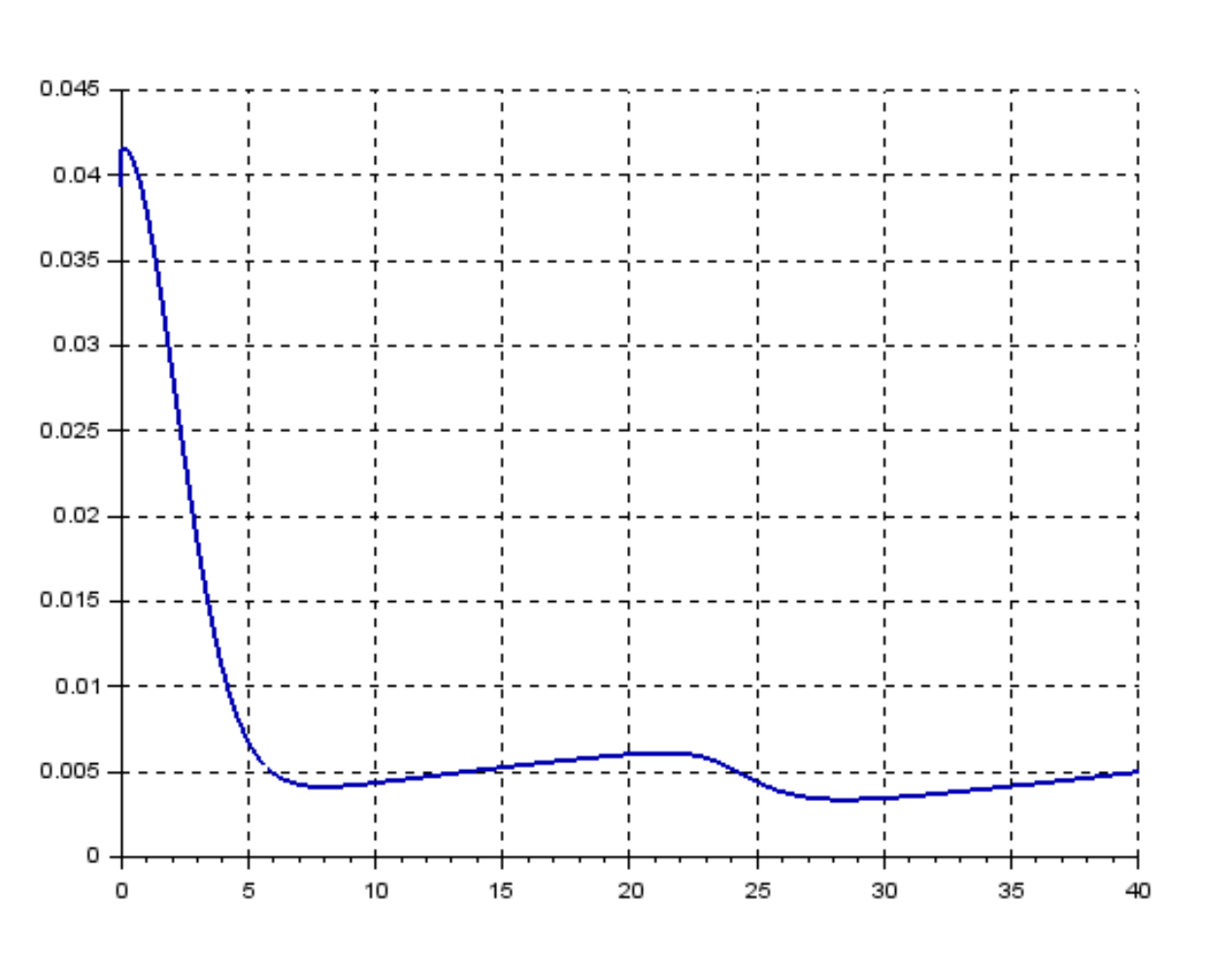}} \qquad
	\subfloat[][$r=10$ and $d=0.5$]{\includegraphics[width=.45\textwidth]{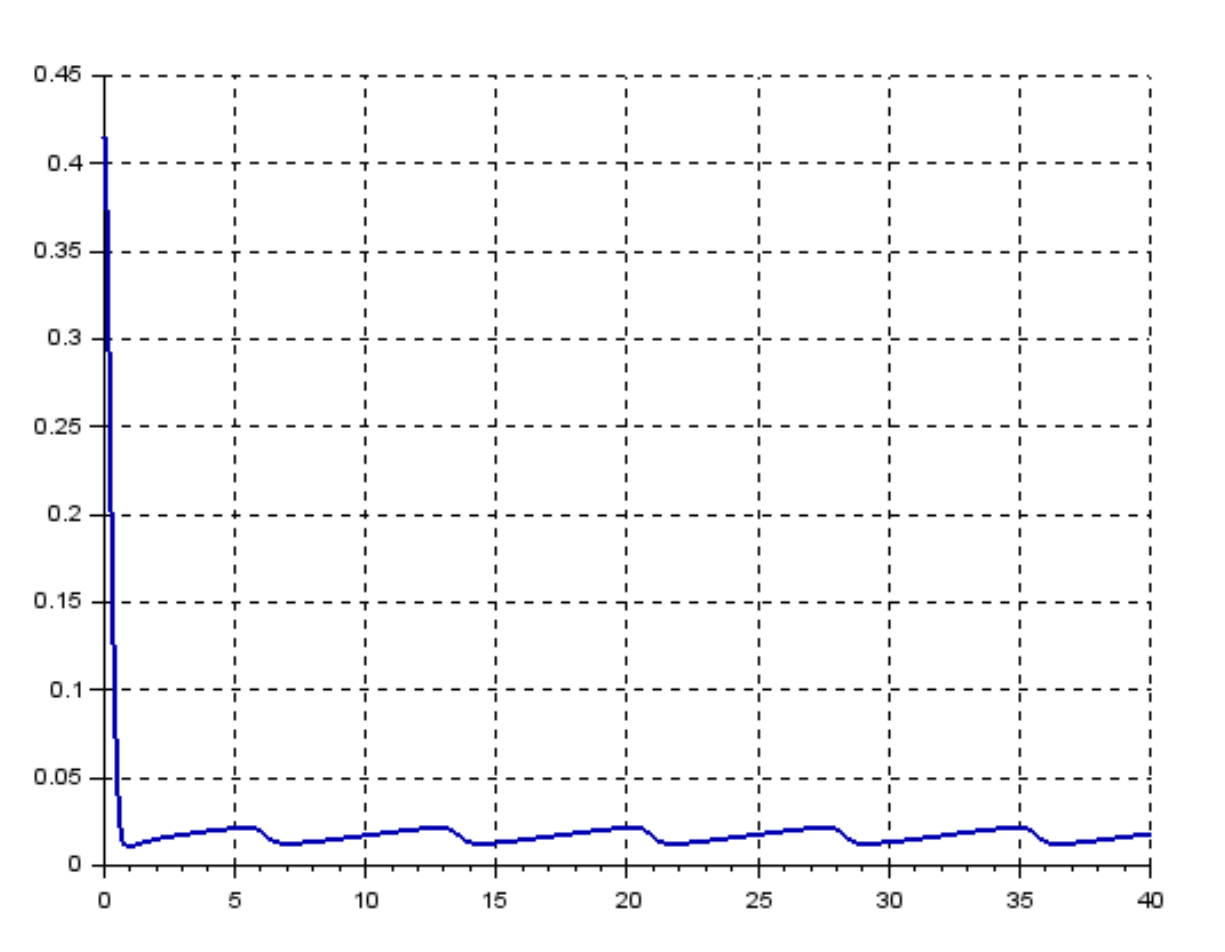}} \\
	\subfloat[][$r=1$ and $d=1.5$]{\includegraphics[width=.45\textwidth]{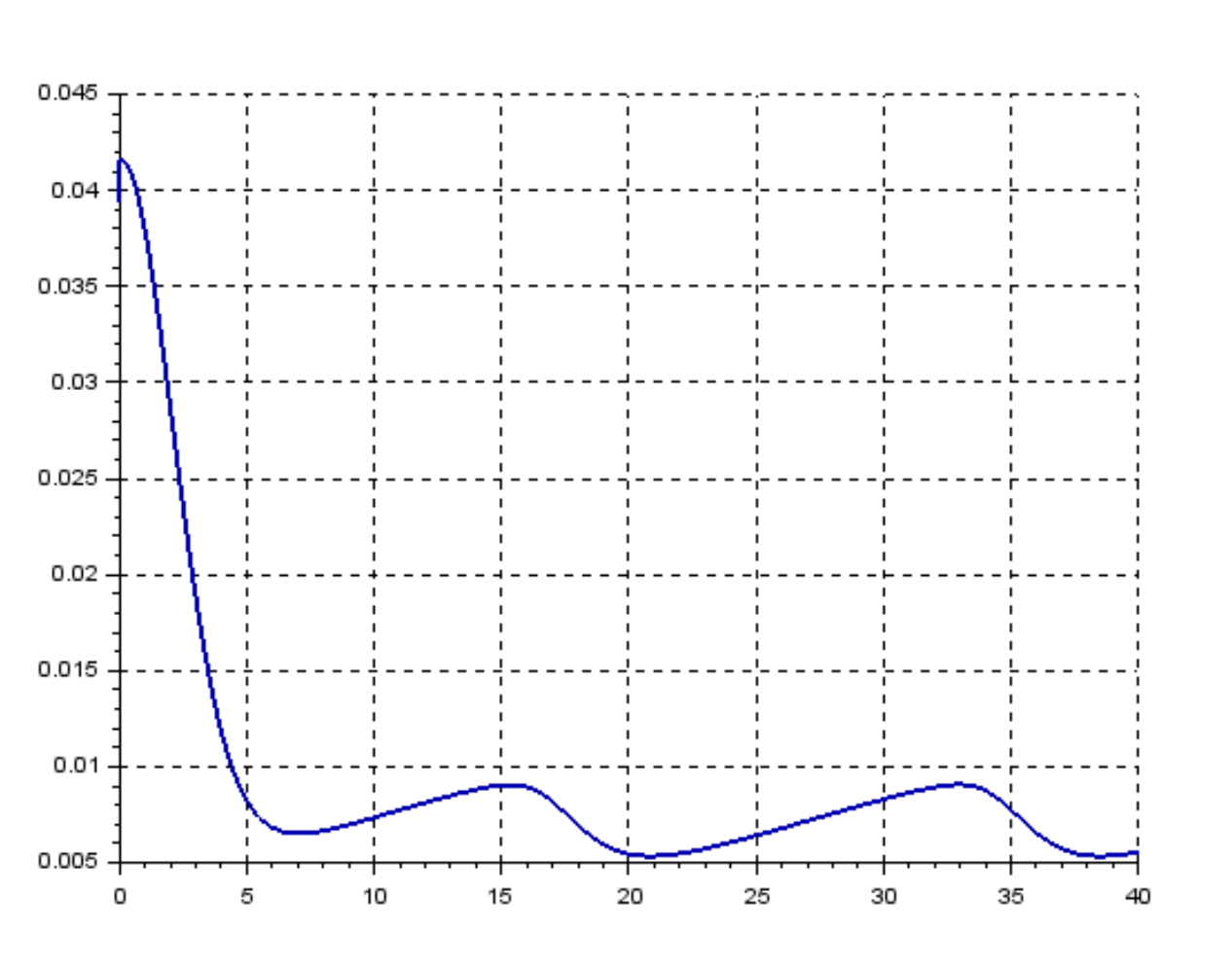}} \qquad
	\subfloat[][$r=10$ and $d=1.5$]{\includegraphics[width=.45\textwidth]{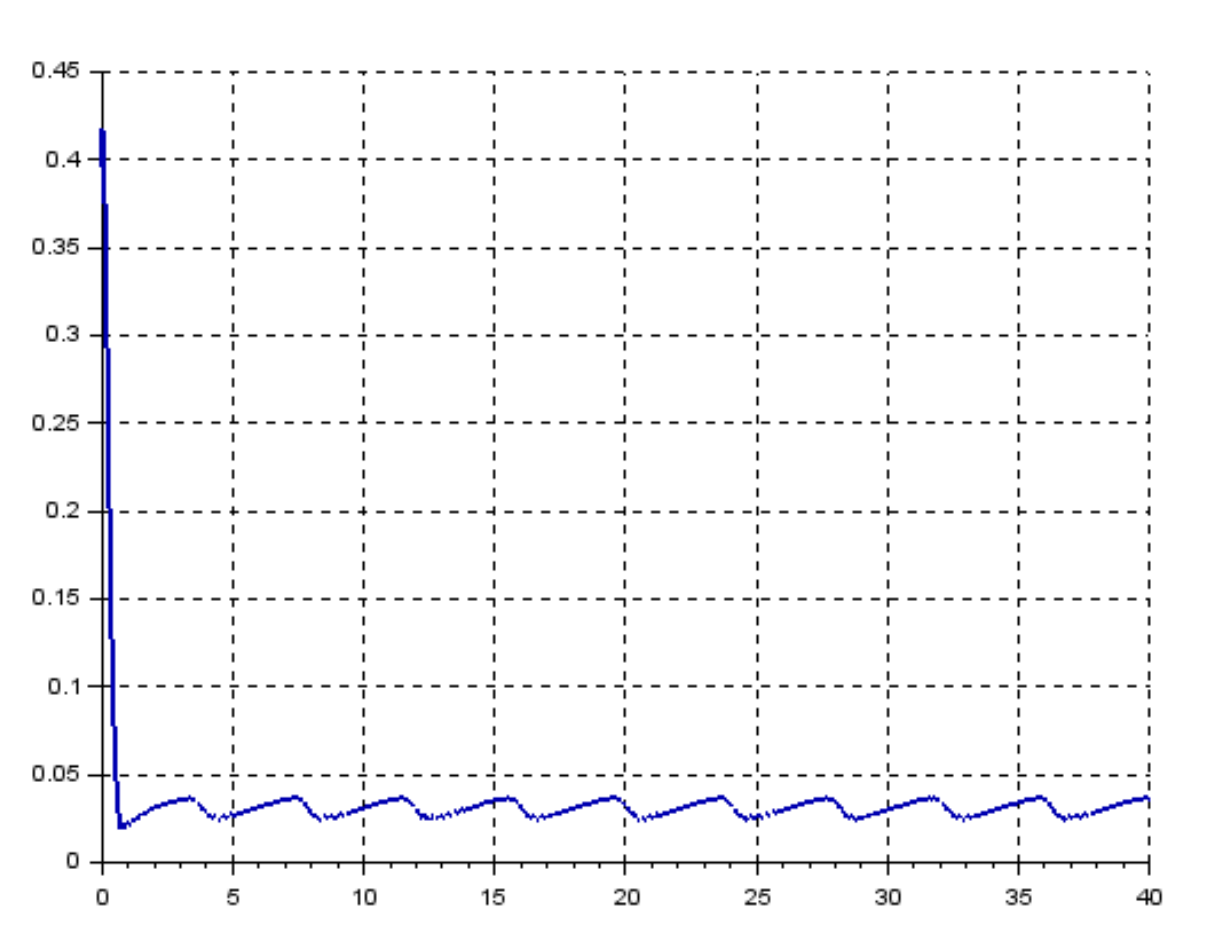}} \\
	\subfloat[][$r=1$ and $d=30$]{\includegraphics[width=.45\textwidth]{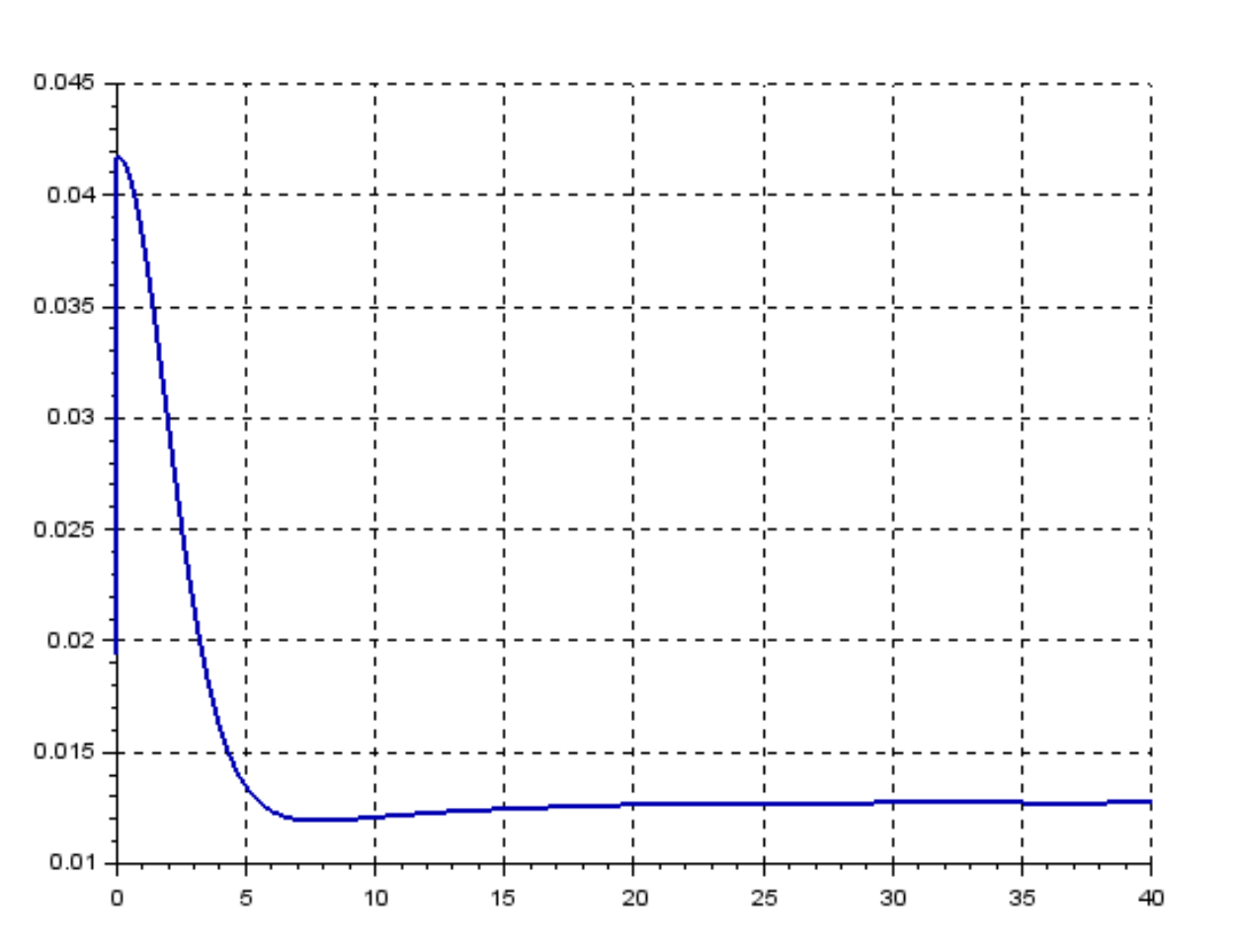}} \qquad
	\subfloat[][$r=10$ and $d=30$]{\includegraphics[width=.45\textwidth]{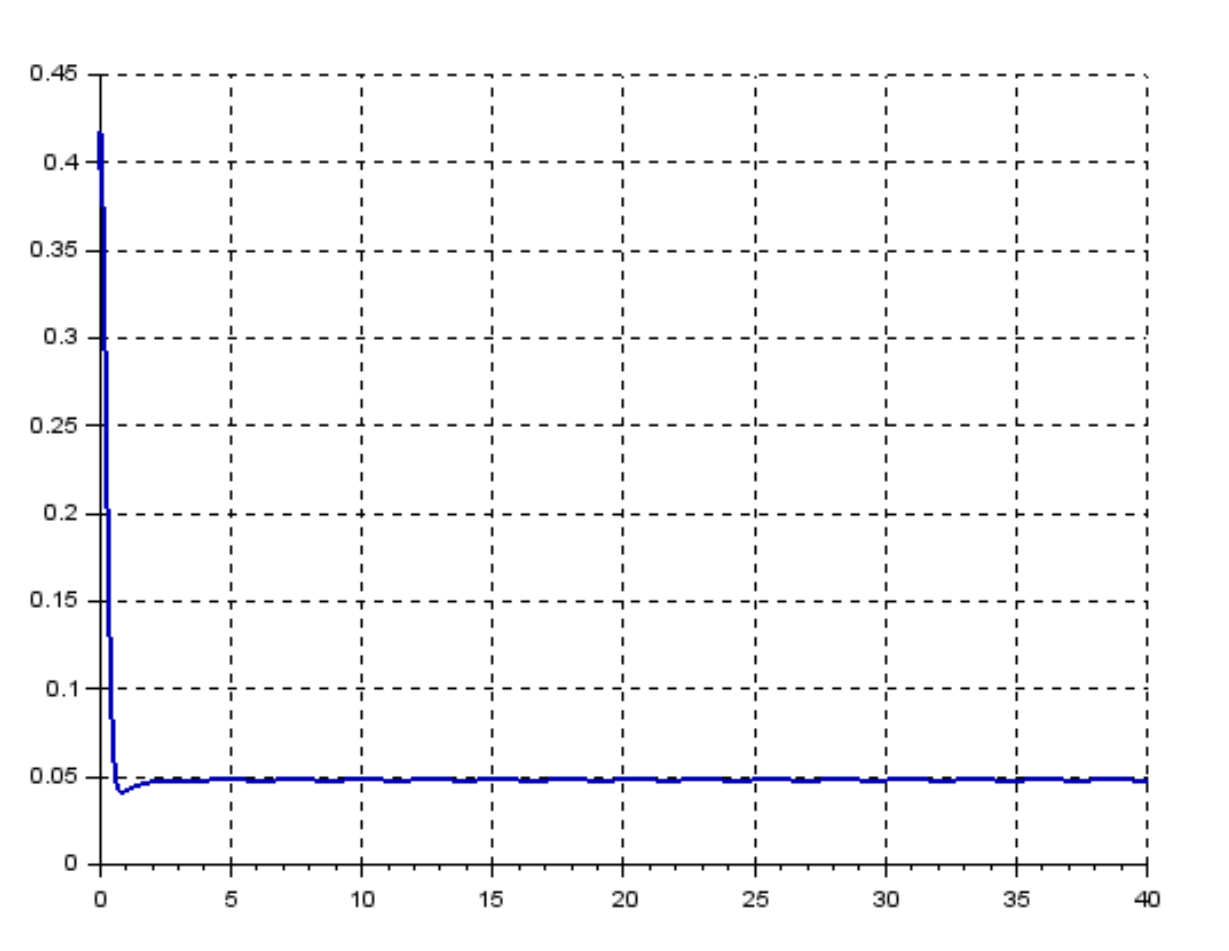}}
	\caption{\footnotesize{Wave speed approximation of the tumour front $v$ for different values of $r$ and $d$.}}
	\label{wavespeed5}
\end{figure}
}

\end{document}